\documentclass[11pt]{article}
\usepackage{graphicx, amsmath, amsfonts, amssymb, amsthm, fullpage, authblk, bbm}
\usepackage[dvipsnames]{xcolor}
\usepackage{algorithm}
\usepackage{algpseudocode}
\algnewcommand\algorithmicinput{\textbf{INPUT:}}
\algnewcommand\INPUT{\item[\algorithmicinput]}
\algnewcommand\algorithmicoutput{\textbf{OUTPUT:}}
\algnewcommand\OUTPUT{\item[\algorithmicoutput]}

\usepackage{hyperref}[]
\hypersetup{
    colorlinks=true,
    linkcolor=blue,
    filecolor=magenta,      
    urlcolor=cyan,
      citecolor=blue,
    }
\usepackage{cleveref}

\usepackage{mathtools}

\newtheorem{theorem}{Theorem}
\newtheorem{lemma}[theorem]{Lemma}
\newtheorem{proposition}[theorem]{Proposition}

\newtheorem{definition}{Definition}
\newtheorem{remark}{Remark}

\allowdisplaybreaks

\DeclareMathOperator*{\argmin}{arg\,min}
\DeclareMathOperator*{\argmax}{arg\,max}
\DeclareMathOperator*{\sign}{sign}

\DeclareMathOperator*{\TV}{TV}
\DeclareMathOperator*{\KL}{KL}
\newcommand*{\dd}{\mathrm{d}}
\DeclareMathOperator{\E}{\mathbb{E}}

\usepackage[round]{natbib}

\title{Federated Differential Privacy with Batch Contamination}
\author{Jason K.~C.~Ng, Thomas B.~Berrett and Yi Yu}
\affil{Department of Statistics, University of Warwick}

\begin{document}

\maketitle

\begin{abstract}
    We study mean estimation under federated differential privacy (FDP) constraints and batch contamination, considering both Gaussian and heavy-tailed distributions. In the homogeneous setting, where sample sizes and privacy budgets are constant across servers, we show that the costs of privacy and contamination disentangle.
    Under a bounded $v$th moment condition, with $v \in (1, \infty)$, the cost of batch contamination is of order $(\varepsilon/n)^{2-2/v}\vee \varepsilon^2/n$, where $\varepsilon$ is the contamination level and $n$ is the per-batch sample size.
    This reveals a phase transition at $\varepsilon\asymp n^{1-v/2}$, for $v > 2$, which is the first of its kind in the literature.   We develop a two-stage procedure that achieves the minimax rate up to poly-logarithmic factors.  A key component is a private smoothed voting procedure tailored to the FDP setting that provides optimal sensitivity control. For heterogeneous batch sample sizes and privacy budgets, we introduce a batch partitioning framework that groups batches with comparable per-batch sample sizes and privacy budgets. Surprisingly, the minimax rate is determined, up to poly-logarithmic factors, by the most informative group of batches. Consequently, in some regimes, it is minimax optimal to use only an arbitrarily small proportion of the data set.  Heterogeneity also introduces interactions between the costs of privacy and contamination, in contrast to their disentanglement in the homogeneous setting. For unknown contamination level, we propose an adaptive, near-minimax optimal Lepski-type extension for the heterogeneous Gaussian setting, designed specifically to prevent privacy leakage by ensuring no repeated use of data points.
\end{abstract}

\section{Introduction}\label{sec-intro}
In a world of ever increasing collection of large data sets, privacy concerns from individuals have prompted the introduction of regulations globally, such as the GDPR \citep{GDPR}, the AI Act \citep[][]{AIact}, and the CCPA in California, USA \citep[][]{CCPA}. Amidst the growing trend of regulation, differential privacy \citep{Dwork2006calibrating} has emerged as a rigorous framework for the research of privacy-preserving methodology. It has since been utilised in industry by technology companies such as Google \citep{erlingsson2014rappor}, and in the public sector \citep{USCensus}.

Partly driven by growing attention to the safety, security and privacy of machine learning and artificial intelligence \citep[e.g.][]{barreno2006can,amodei2016concrete}, there has been renewed interest in robustness in recent years. Since the discovery of nonrobustness in machine learning algorithms \citep[e.g.][]{szegedy2013intriguing} and the connections between privacy and robustness \citep[e.g.][]{dwork2009differential}, there has been a flurry of work on algorithms robust against data contamination \citep[see e.g.][for a review]{Diakonikolas2023algorithmic}, some of which are also simultaneously private \citep[see e.g.][]{Liu2021robust,liu2023label,Li2023on}.

Meanwhile, security risks and communication bottlenecks arising from the transfer of massive raw data sets pushes for decentralised learning. Federated learning \citep[e.g.][]{mcmahan2017communication,li2020federated, baduwal2026federated} is a paradigm addressing the challenges of training a global model using data from multiple local sources. Information exchange is facilitated through the transmission of summary statistics (e.g.~gradients, model parameter updates), such that each individual's raw data never leave the local database, enhancing the privacy protection of local users. However, summary statistics can still be used to infer information about each individual \citep[see e.g.][]{wang2019beyond,zhu2019deep,geiping2020inverting}, prompting a line of work on federated learning under various models of differential privacy \citep[e.g.][]{geyer2017differentially,cheu2019distributed,girgis2021shuffled, kairouz2021distributed}.   

To address the above challenges, we study robustness and privacy preservation under the setting of federated learning. The rest of this paper is organised as follows. \Cref{sec-general-setup} introduces the general setup of our problem, followed by a list of the main contributions of this paper in \Cref{sec-contributions}. \Cref{sec-related-work} gives a review of related work. We study the univariate mean estimation problem with homogeneous batch sample sizes and privacy budgets in \Cref{sec-mean-est}, where we consider the Gaussian setting in \Cref{sec-Gaussian-mean} and the heavy-tailed setting in \Cref{sec-finite-moments}. In \Cref{sec-mean-est-hetero}, we generalise to heterogeneous batch sample sizes and privacy budgets under the Gaussian setting. Details of the proofs are deferred to the Appendices.

\subsection{General setup}\label{sec-general-setup}
We consider a distributed data setting, where each server holds a single batch of data. Specifically, we have $K\in \mathbb{Z}_+$ local servers, and for every $k \in \{1,\dots,K\}$, the $k$th server holds a data set~$D_k$. Throughout this paper, we use the terms server and batch interchangeably. In the rest of this section, we will formulate the general problem, including the specific privacy constraints and contamination model under consideration.

\subsubsection{Privacy}
We consider the setting of federated differential privacy (FDP), where only summary statistics that are privatised can be communicated across servers. The privacy mechanisms used at each server must satisfy the central differential privacy (CDP) constraint. This notion of FDP we consider is most similar to that in the works of, for example,~\cite{Li2024federated, cai2024optimal, xue2024optimal}. Specifically, at each server $k \in \{1,\dots,K\}$, we output a privatised statistic $Z_k\in \mathcal{Z}_k$ generated by a privacy mechanism $Q_k$ that satisfies the following definition of FDP. 
\begin{definition}[Federated differential privacy, FDP]\label{FDP-def}
    Given a collection of privacy parameters $\{\alpha_k\}_{k=1}^K$, a privacy mechanism $Q = \{Q_k\}_{k=1}^K$ satisfies $\{\alpha_k\}_{k=1}^K$-FDP if, for each $k\in \{1,\dots,K\}$,~$Q_k$ is a conditional distribution on $\mathcal{Z}_k$ satisfying
    \[Q_k(Z_k \in A|D_k, Z_1,\dots, Z_{k-1}) \leq e^{\alpha_k} Q_k(Z_k \in A|D'_k, Z_1,\dots, Z_{k-1}),\]
    for all measurable sets $A \in \sigma(\mathcal{Z}_k)$, outputs $(Z_1,\dots,Z_{k-1}) \in \mathcal{Z}_1 \times \dots \times \mathcal{Z}_{k-1}$, and adjacent data sets $D_k, D'_k$ that differ by at most one entry.
\end{definition}
The parameter $\alpha_k>0$ specifies the level of privacy guarantee required at the $k$th server; a smaller value of $\alpha_k$ provides a stronger guarantee. This definition allows for sequential interactivity among servers. When each server holds only one sample and $\alpha_1 = \cdots = \alpha_K$, FDP is equivalent to an $\alpha_1$-local differential privacy (LDP) constraint. When there is only one server ($K=1$), it is equivalent to an $\alpha_1$-CDP constraint. 

Without ambiguity, in the sequel we let $\mathcal{Q}\equiv \mathcal{Q}_{\{\alpha_k\}_{k=1}^K}$ denote the set of privacy mechanisms satisfying the FDP constraint.

\subsubsection{Batch contamination}
For a class of distributions $\mathcal{P}$ of interest on the sample space $\mathcal{X}$ and contamination level $\varepsilon \in [0, 1/2)$, the Huber $\varepsilon$-contamination model \citep[e.g.][]{Huber1964robust} is defined by the collection of distributions
\begin{equation}\label{Huber-model}\{P_\varepsilon = (1-\varepsilon)P+\varepsilon G: P\in \mathcal{P}, G \in \mathcal{G}\},\end{equation}
where $\mathcal{G}$ denotes the set of all distributions on $\mathcal{X}$. We consider a batch contamination variant, i.e.
\begin{equation}\label{general-data-distribution}
    D_k \overset{\mathrm{ind.}}{\sim} P_{\varepsilon,k} = (1-\varepsilon)P^{\otimes n_k}+\varepsilon G_k
\end{equation}
for some $\varepsilon\in [0,1/2)$, $P\in \mathcal{P}$ and $G_k\in \mathcal{G}_k$, where $\mathcal{G}_k$ denotes the set of all distributions on $\mathcal{X}^{\otimes n_k}$. The \textit{inlier distribution} $P^{\otimes n_k}$ of server $k$ is an $n_k$-fold product distribution, with $P$ being identical across all servers. With probability $\varepsilon$, independently of all other servers, each server $k\in \{1,\dots,K\}$ is contaminated by an arbitrary \textit{outlier distribution} $G_k$ on $\mathcal{X}^{\otimes n_k}$. 
\begin{remark}
    Unlike the inlier distribution, the outlier distribution $G_k$ is arbitrary and not necessarily a product distribution. In addition, unlike the standard Huber model, outlier distributions of different servers are allowed to be arbitrarily different.
\end{remark}
With \eqref{general-data-distribution}, the class of distributions considered is
\begin{equation}\label{general-contamination-model}\mathcal{P}_\varepsilon(\mathcal{P}) = \{\{P_{\varepsilon,k}\}_{k=1}^K:P\in \mathcal{P}, G_k \in \mathcal{G}_k \text{ for each } k\in \{1,\dots,K\}\}.\end{equation}

\subsubsection{Minimax risk}
To quantify the fundamental difficulty of a statistical problem, we consider the notion of private minimax risk $\mathcal{R}$ under the squared error loss defined by
\begin{equation}\label{general-minimax-risk}\mathcal{R} = \inf_{Q\in \mathcal{Q}}\inf_{\hat{\theta}}\sup_{\{P_{\varepsilon,k}\}_{k = 1}^K\in\mathcal{P}_\varepsilon(\mathcal{P})}\mathbb{E}_{\{P_{\varepsilon,k}\}_{k=1}^K,Q}[\{\hat{\theta}-\theta(P)\}^2],\end{equation}
where $\mathcal{P}_{\varepsilon}(\mathcal{P})$ denotes the class of data generating mechanisms under the batch contamination model defined in~\eqref{general-contamination-model}, $\theta(P)$ is a population functional of interest, the expectation is with respect to both the data generating mechanism $\{P_{\varepsilon,k}\}_{k=1}^K$ and the privacy mechanism $Q$, and the second infimum is over all measurable functions $\hat{\theta}(Z_1,\dots,Z_K)$ of the privatised statistics $(Z_1,\dots,Z_K)$ generated from the privacy mechanism $Q$. In this paper, we consider $\theta(P)$ to be the mean of a univariate distribution $P$. The study of the univariate inlier mean estimation problem serves as the blueprint to understand more complex problems under batch contamination and FDP.

\subsection{Contributions}\label{sec-contributions}
We summarise below the key contributions of this paper.

\begin{itemize}
    \item To the best of our knowledge, this is the first study of batch contamination under FDP, providing a framework for developing robust methodology with rigorous privacy guarantees in federated learning.  
    In \Cref{sec-Gaussian-mean}, we establish the minimax rate for Gaussian mean estimation with homogeneous batch sizes and privacy budgets, 
    and show that costs of privacy and contamination disentangle, extending the findings of \cite{Liu2021robust} and    \cite{Li2023on} to the federated setting.  
    Our two-stage estimator attains this rate up to poly-logarithmic factors.  A key component is a  
    novel smoothed voting procedure tailored to FDP constraints.  By replacing binary votes with smoothed votes, we avoid abrupt changes in a server's vote caused by a single observation and achieve optimal sensitivity control. 
    \item In \Cref{sec-finite-moments}, we extend the theory to heavy-tailed inlier distributions satisfying only a bounded $v$th-moment condition,  for any fixed $v\in (1,\infty)$. 
    To the best of our knowledge, a minimax characterisation covering this full range of moment assumptions is new both for batch contamination without privacy constraints and for FDP without contamination.  
    We establish matching upper and lower bounds up to poly-logarithmic factors, which show that the costs of privacy and contamination remain disentangled. 
    In particular, the cost of batch contamination is of order $(\varepsilon/n)^{2-2/v}\vee\varepsilon^2/n$, revealing a new phase transition for $v\in (2,\infty)$ at $\varepsilon \asymp n^{1-v/2}$.  The two terms reflect the effects of contamination in the tails and the main bulk of the inlier distribution, respectively. 
    
    \item In \Cref{sec-mean-est-hetero}, we characterise the minimax rate of Gaussian mean estimation under batch contamination with fully  
    heterogeneous batch sample sizes and privacy budgets. 
    We introduce a batch partitioning framework that groups servers with comparable batch sample sizes and privacy budgets.  
    Surprisingly, the most informative group of batches alone determines the minimax rate of the problem up to poly-logarithmic factors.  
    Consequently, in some regimes, near-optimal estimation requires only an 
    arbitrarily small proportion of the data. 
    Heterogeneity  
    also introduces interactions between the costs of privacy and contamination, in contrast to the homogeneous setting and the findings in \cite{Liu2021robust} and \cite{Li2023on}. Specifically, the group that is optimal under privacy alone may differ from that under contamination alone, and the joint minimax risk can be strictly larger than the maximum of the risks arising from the two requirements considered separately.
    \item Our two-stage procedure adapts automatically to the unknown contamination level $\varepsilon$   
    in the homogeneous Gaussian setting.  For the heterogeneous Gaussian setting, \Cref{sec-adaptation} introduces a Lepski-type procedure that adapts the choice of batch group to the contamination level.  The procedure exploits our batch partitioning framework.  Each group estimator is computed only once and reused whenever the same group is selected for another candidate contamination level. As a result, each data point is used at most once, preserving FDP guarantees without composing privacy budgets across iterations. 
\end{itemize}

\subsection{Related work}\label{sec-related-work}

\textbf{Robust statistics. }The field of robust statistics emerged in the 1960s through the work of \cite{tukey1960survey} and \cite{Huber1964robust} on contaminated distributions, prompting the development of classical methods for robust estimation, such as the trimmed mean and M-estimators. Hampel formalised notions of the influence function and breakdown point \citep[e.g.][]{hampel1968contributions, hampel1974influence}, providing a framework for the analysis of  stability and breakdown properties of statistical procedures. These foundational concepts inspired a great amount of research in classical robust statistics; see a number of textbooks \citep[e.g.][]{huber2009robust, hampel2005robust, maronna2006robust} for an overview. Modern developments have seen a shift in emphasis towards non-asymptotic analysis of various statistical tasks under contamination models, such as robust estimation \citep[e.g.][]{Chen2018robust, prasad2020robust}, robust testing \citep[e.g.][]{Chen2016a,canonne2026full}, robust regression \citep[e.g.][]{dalalyan2019outlier,gao2020robust} and density estimation \citep{Liu2019density,uppal2020robust}. Some of the proposed algorithms, despite their statistical optimality, are computationally intractable in high dimensions, a classical example being the depth-based Tukey median \citep{tukey1975mathematics}. Recent breakthroughs have seen tractable algorithms \citep[e.g.][]{lai2016agnostic,diakonikolas2019robust} that are robust against contamination even in high dimensions; see \cite{diakonikolas2019recent} for a survey. As interest in practical algorithms grows, one line of work is on developing algorithms that do not require knowledge of the contamination level. Recent work include mean estimation \citep[e.g.][]{Dalalyan2022all-in-one}, density estimation \citep{Liu2019density} and adaptive confidence sets \citep{luo2026adaptive}; \cite{Jain2022robust} suggests applying a Lepski-type method \citep{Lepskii1991on,Lepskii1992asymptotically, lepskii1993asymptotically} as a general approach for adapting algorithms to unknown contamination level. 

\medskip
\noindent\textbf{Heavy-tailed distributions and robust statistics. }Even without the presence of contamination, standard methods, such as the empirical mean, exhibit poor finite-sample performance under heavy-tailed distributions.  However, assuming only finite variance, it was shown that the optimal finite-sample sub-Gaussian error bound can be attained by robust mean estimators such as Catoni's M-estimator \citep{catoni2012challenging}, the median-of-means \citep[e.g.][]{Devroye2016sub-gaussian} and the trimmed mean (or truncated mean) \citep[e.g.][]{lugosi2021robust}; see \cite{lugosi2019mean} for a survey. Beyond mean estimation, regression \citep[e.g.][]{hsu2014heavy}, covariance estimation \citep[e.g.][]{wei2017estimation,minsker2018sub}, bandit problems \citep[e.g.][]{bubeck2013bandits, agrawal2021regret}, among other statistical problems, have also been studied under heavy-tailed distributions.  In fact, analysis of the heavy-tailed model is deeply connected to robust statistics.  \cite{prasad2019unified} developed optimal algorithms that are simultaneously robust under the heavy-tailed model and Huber's $\varepsilon$-contamination model. Algorithms robust against both adversarial contamination and heavy tails have also been studied \citep[e.g.][]{diakonikolas2020outlier,pensia2025robust}. More generally, \cite{cherapanamjeri2025heavy} showed that an optimal mean estimation algorithm robust against adversarial contamination is also optimal under the heavy-tailed model.

\medskip
\noindent\textbf{Batch contamination. }Introduced by \cite{qiao2017learning}, a related line of work studies robust learning from data collected in batches, a fraction of which may be adversarially contaminated. \cite{qiao2017learning} initiated the study of the problem of learning discrete distributions to within a small total variation (TV) distance, which is built upon by \cite{chen2020efficiently} and \cite{Jain2020optimal}. \cite{Jain2020a} extends the framework to include estimation of certain structured continuous distributions and robust classification. While these earlier works focus on distribution estimation, more recent work by \cite{aliakbarpour2026high} considers multivariate mean estimation assuming finite covariance, with adversarial corruption at the batch-level, while also incorporating heterogeneity by allowing mean shifts and adversarial contamination among inlier batches. \cite{tian2026contaminated} integrates contaminated batches into a transfer and multi-task learning framework. In contrast to previous work, we consider the extensions of general bounded-moment assumptions on the inlier distribution or general batch sample size heterogeneity.

\medskip
\noindent\textbf{Privacy. }Since the introduction of the framework of differential privacy \cite[e.g.][]{Dwork2006calibrating}, privacy has been incorporated in mean and covariance matrix estimation \cite[e.g.][]{biswas2020coinpress,amin2019differentially,wang2021differentially}, density estimation \citep[e.g.][]{lalanne2023cost}, hypothesis testing \citep[e.g.][]{acharya2018differentially,canonne2019structure,berrett2020locally}, regression \citep[e.g.][]{alabi2020differentially,bernstein2019differentially, milionis2022differentially}, classification \citep[e.g.][]{berrett:hal-03603855,auddy2026minimax}, change point detection \citep[e.g.][]{cummings2018differentially,berrett2021locally}, and many other statistical tasks. Due to their popularity, information-theoretic tools for studying statistical optimality have been well-established for both the central model \citep[e.g.][]{wasserman2010statistical,Bun2014fingerprinting,Cai2021cost,cai2023score} and the local model \citep[e.g.][]{Duchi2018minimax} of differential privacy. The model of differential privacy that we consider under the setting of federated learning most closely resembles that of \cite{Li2024federated} and  \cite{cai2024optimal}, which can be seen to blend the features of the central and local models of privacy. Under this framework, problems on mean estimation and regression \citep{Li2024federated}, nonparametric regression \citep{cai2024optimal}, functional mean estimation \citep{xue2024optimal,cai2026optimal}, classification \citep{auddy2026minimax, xue2026federated}, survival analysis \citep{Hung2025optimal}, among others, have been studied.

\medskip
\noindent\textbf{Privacy and robust statistics. }The connections between privacy and robust statistics have been studied since the early developments of differential privacy. Recognising the common goal of being insensitive to perturbations, \cite{dwork2009differential} studies the use of robust estimators as a starting point for the construction of differentially private estimators. Subsequent work developed this idea, drawing upon classical robust estimators such as the trimmed mean \citep{smith2011privacy} and M-estimators \citep{lei2011differentially}, while \cite{chaudhuri2012convergence} links privacy with the notions of the influence function and gross error sensitivity. \cite{avella2021privacy} builds upon this line of work through the study of differentially private testing. General approaches of converting between robust and private algorithms have been studied by \cite{liu2022differential},  \cite{hopkins2023robustness} and \cite{asi2023robustness}.

More recently, there has been growing interest on methods that are simultaneously private and robust against contamination or heavy tails, ranging from mean estimation \citep[e.g.][]{Kamath20private,Liu2021robust}, regression \citep[e.g.][]{liu2023label}, and principal component analysis \citep[e.g.][]{kim2025robust}, to name but a few. Beyond the CDP setting, robustness is also studied under LDP constraints as interest in distributed data grows. \cite{Li2023on} studied Huber's contamination model under LDP through the problems of testing, mean estimation and density estimation. Other works,  \citep[e.g.][]{Cheu2021manipulation,acharya2021robust, Chhor2023robust} study a stronger adversarial model with LDP, where privatised data can be attacked or contaminated by an adversary. In particular, \cite{Chhor2023robust} considers contamination in batches. Our work aims to bridge the gap between existing work on robustness under CDP and LDP using the less studied frameworks of FDP and batch contamination, which are particularly suited to distributed and federated learning.

\subsection{Notation}
For $n\in \mathbb{N}$, denote $[n] = \{1,\dots,n\}$. For $a,b,x\in \mathbb{R}$, let $a\wedge b = \min\{a,b\}$, $a\vee b = \max\{a,b\}$, $a_+ = a\vee 0$, $[x]_a^b = (x\vee a)\wedge b$, $\lfloor a\rfloor$ be the largest integer less than or equal to $a$, and $\lceil a\rceil$ be the smallest integer greater than or equal to $a$. For positive sequences $\{a_m\}$ and $\{b_m\}$, we write $a_m \lesssim b_m$ if there exists a constant $C>0$ such that $\sup_{m}a_m/b_m\leq C$, with $a_m \gtrsim b_m$ denoting $b_m \lesssim a_m$, and $a_m \asymp b_m$ denotes $a_m \lesssim b_m \lesssim a_m$. Further, write $a_m\lesssim_{\log} b_m$, $a_m\gtrsim_{\log} b_m$, and $a_m \asymp_{\log} b_m$, if $a_m\lesssim b_m$, $a_m\gtrsim b_m$, and $a_m \asymp b_m$ up to poly-logarithmic factors in $m$. 
Denote by $\mathbf{1}\{A\}$ the indicator function of the event $A$. Given a distribution $P$, the distribution $P^{\otimes n}$ denotes the $n$-fold product distribution arising from $n$ i.i.d.~observations drawn from $P$. For distributions $P$ and $Q$, their TV distance is $\TV(P,Q) = \sup_A|P(A)-Q(A)|$, where the supremum is taken over all measurable sets~$A$; their Kullback-Leibler divergence is $\KL(P,Q) = \int \log(\dd P /\dd Q) \,\dd P$ if $P$ is absolutely continuous with respect to $Q$. We say a random variable $X$ is sub-Gaussian with variance parameter $\sigma^2$ if $\E \{e^{\lambda (X-\E X)}\}\leq e^{\lambda^2\sigma^2/2}$ for all $\lambda \in \mathbb{R}$. We say a random variable $X$ is sub-exponential with parameters $(\nu,b)$ if $\E \{e^{\lambda(X-\E X)}\} \leq e^{\nu^2 \lambda^2/2}$ for all $|\lambda|<1/b$. We define the $\mathrm{Laplace}(\mu,\lambda)$ distribution to be the distribution with density $(2\lambda)^{-1}e^{-|x-\mu|/\lambda}$, and the standard Laplace distribution to be the $\mathrm{Laplace}(0,1)$ distribution.

\section{Mean estimation with homogeneous servers}\label{sec-mean-est}
In this section, we investigate the homogeneous mean estimation problem. Specifically, under the general setup described in \Cref{sec-general-setup}, we restrict to a common batch size $n = n_1 = \dots = n_K$ and a common privacy parameter $\alpha = \alpha_1=\dots =\alpha_K$. We first study the Gaussian setting in \Cref{sec-Gaussian-mean}, characterising its fundamental statistical difficulty and detailing a two-step estimation procedure that is minimax optimal up to poly-logarithmic factors. We then extend to the heavy-tailed setting in \Cref{sec-finite-moments}, proving matching lower and upper bounds, up to poly-logarithmic factors, of the minimax risk using a version of the two-step procedure proposed in \Cref{sec-Gaussian-mean}.
\subsection{Gaussian mean estimation}\label{sec-Gaussian-mean}

In this subsection, we consider the class of univariate Gaussian distributions
\begin{equation}\label{Gaussian-distribution}\mathcal{P} = \{\mathcal{N}(\theta,1): \theta \in [-1,1]\}.\end{equation} 
For $k \in [K]$, denoting the data set stored in the $k$th server by $D_k = \{X_{k,i}\}_{i = 1}^{n}$, model (\ref{general-data-distribution}) reduces to
\begin{equation}\label{Gaussian-model}
    X_k = (X_{k, 1}, \ldots, X_{k, n}) \overset{\mathrm{ind.}}{\sim} P_{\varepsilon,k} = (1 - \varepsilon) \mathcal{N}^{\otimes n}(\theta, 1) + \varepsilon G_k,
\end{equation}
where $G_k$ is any distribution on $\mathbb{R}^n$.

Considering the mean estimation problem under the squared loss, the minimax risk in (\ref{general-minimax-risk}) takes the form
\begin{equation}\label{Gaussian-minimax-risk}
\mathcal{R}_\alpha(\varepsilon) = \inf_{Q\in \mathcal{Q}}\inf_{\hat{\theta}}\sup_{\{P_{\varepsilon,k}\}_{k=1}^K\in\mathcal{P}_\varepsilon(\mathcal{P})} \mathbb{E}_{\{P_{\varepsilon,k}\}_{k=1}^K,Q}[\{\hat{\theta}-\theta(P)\}^2],
\end{equation}
where $\theta(P)$ denotes the mean of the distribution $P$ and, for brevity, we write $\mathcal{Q}$ to denote the set of privacy mechanisms with the common privacy parameter $\alpha$ in this homogeneous setting. We establish the minimax rate of this problem, up to poly-logarithmic factors, in the following theorem.
\begin{theorem}\label{thm-Gaussian-mean-minimax-rate}
    Consider the homogeneous Gaussian mean estimation problem defined by \eqref{Gaussian-distribution} and~\eqref{Gaussian-model}. There exists an absolute constant $c>0$ such that when the condition 
    \begin{equation}\label{UB-condition-thm1}
    \varepsilon+ \sqrt{\frac{\log(K(n\wedge n^2\alpha^2))}{K}}+\frac{1}{\sqrt{K}(\sqrt{n}\wedge n\alpha)}< c
    \end{equation}
    holds, the minimax risk defined in (\ref{Gaussian-minimax-risk}) satisfies
    \begin{equation} \label{Eq:HomogeneousGaussianRate} \mathcal{R}_\alpha(\varepsilon) \asymp_{\log} \frac{1}{nK}\vee\frac{1}{n^2 K \alpha^2}\vee \frac{\varepsilon^2}{n}. \end{equation}
\end{theorem}

\Cref{thm-Gaussian-mean-minimax-rate} characterises the minimax rate with the lower bound in \Cref{LB-Gaussian} and the upper bound in \Cref{UB-thm} respectively.  The condition \eqref{UB-condition-thm1} is used in the upper bound only. This is nevertheless a mild assumption.  It is common in the robust statistics literature to assume a small constant upper bound on $\varepsilon$ \citep[e.g.][]{Diakonikolas2023algorithmic}, whereas the violation of the constant lower bound on $\sqrt{K}(\sqrt{n}\wedge n\alpha)$ results in a regime where consistent estimation is impossible, as is evident from \Cref{LB-Gaussian}. The bound on the remaining term is more subtle; it requires the number of servers $K$ to be at least of order $\log(K(n\wedge n^2\alpha^2))$. In view of \Cref{LB-Gaussian} below, it is impossible to upper bound the minimax risk by the rate~\eqref{Eq:HomogeneousGaussianRate} without assuming a lower bound on $K\log(1/\varepsilon)$ of this same order. Up to the $\log(1/\varepsilon)$ factor, we therefore see that this condition is necessary for~\eqref{Eq:HomogeneousGaussianRate} to hold.  

In the regime of interest where the condition \eqref{UB-condition-thm1} is satisfied, the minimax rate is comprised of three terms, which can be interpreted as the fundamental difficulty of the problem due to the usual statistical sampling error, the cost of privacy and the effect of contamination, respectively.
In the non-private and non-contaminated setting ($\alpha = \infty$ and $\varepsilon=0$), the recovered parametric rate of $1/(nK)$ due to sampling error is a standard result. 

In the presence of batch contamination, the non-private rate is $1/(nK) \vee \varepsilon^2/n$. Assuming that the cost of contamination dominates, it is well-known that univariate Gaussian mean estimation with known variance $\sigma^2$ under the standard Huber $\varepsilon$-contamination model \eqref{Huber-model} is characterised by the rate $\varepsilon^2 \sigma^2$; this can be seen in, for example, Theorem 2.2 in \cite{Chen2018robust} after scaling to the variance $\sigma^2$.  
Since the sample mean is a sufficient statistic, one can draw a connection to our results by considering the sample mean within servers $\bar{X}_{k} = n^{-1}\sum_{i=1}^n X_{k,i}$, which satisfies
\[\bar{X}_k \overset{\mathrm{ind.}}{\sim}(1 - \varepsilon) \mathcal{N}(\theta, 1/n) + \varepsilon \bar{G}_k,\]
where $\bar{G}_k$ is a distribution on $\mathbb{R}$ for each $k\in [K]$ and the variance is set to $\sigma^2 = 1/n$. \Cref{thm-Gaussian-mean-minimax-rate} demonstrates that, with the additional structure of the batch contamination model, the effect of the contamination term $\varepsilon^2/n$ decreases as the server size $n$ increases, even as the average number of corrupted samples grows as $\varepsilon nK$ and the proportion of corrupted samples remains at $\varepsilon$. This is consistent with findings in previous work on learning inlier distributions under adversarial batch contamination \citep[e.g.][]{Jain2020a, Jain2020optimal}, and in contrast to the unavoidable cost of $\varepsilon^2$ under the usual Huber $\varepsilon$-contamination model \eqref{Huber-model}. 

In the private setting without contamination ($\varepsilon=0$), we recover the FDP rate $K^{-1}(n\wedge n^2\alpha^2)^{-1}$ consistent with previous work \citep[e.g.][Theorem 2]{Li2024federated}. The quantity $K(n\wedge n^2\alpha^2)$ can be viewed as the effective sample size for the FDP problem, with a part due to the usual statistical error, and a part due to privatisation variance. The rate interpolates the rates $(n\wedge n^2\alpha^2)^{-1}$ of the CDP ($K=1$) regime \citep[e.g.~Theorem 3.1 in][]{Cai2021cost} and $K^{-1}(1\wedge \alpha^2)^{-1}$ of the LDP ($n= 1$) regime \citep[e.g.~Corollary 5 in][]{duchi2019lower}. 

Under the presence of both batch contamination and the FDP constraint, note that the minimax rate~\eqref{Eq:HomogeneousGaussianRate} contains no interaction between $\varepsilon$ and $\alpha$, such that the costs of contamination and privacy are separated. Consequently, there is no additional cost in combining robustness and privacy. To be specific, modifying a robust procedure to satisfy FDP constraints need not incur additional statistical costs when $\alpha\gtrsim n^{-1/2} \wedge ( \varepsilon\sqrt{nK})^{-1}$, and the same is true for incorporating robustness to a FDP procedure when $\varepsilon \lesssim K^{-1/2} \vee (\alpha\sqrt{nK})^{-1}$. The same phenomenon was observed by previous work on Huber contamination under LDP \citep{Li2023on} and adversarial contamination under CDP \citep[e.g.][]{Liu2021robust}; \Cref{thm-Gaussian-mean-minimax-rate} shows that this extends to the batch contamination model under FDP. In contrast, previous work on LDP with adversarial contamination, where data contamination can take place after data privatisation, have found that an interaction term depending on both $\varepsilon$ and $\alpha$ is inevitable. This is the case for both standard adversarial contamination \citep[e.g.][]{Cheu2021manipulation} and batch adversarial contamination \citep{Chhor2023robust}. The key difference is that we assume no contamination takes place after the privatisation process. 

\subsubsection{Lower bound}\label{sec-uni-mean-under-priv-LB}
We now establish the lower bound on the minimax risk in detail. 
\begin{proposition}\label{LB-Gaussian} Consider the homogeneous Gaussian mean estimation problem detailed in~\eqref{Gaussian-model}. The minimax risk defined in (\ref{Gaussian-minimax-risk}) satisfies

\begin{equation}\label{LB-Gaussian-eq}\mathcal{R}_\alpha(\varepsilon)\geq c\left(\frac{1}{nK}\vee\frac{1}{1\vee n^2 K \alpha^2}\vee \frac{\varepsilon^2}{n}\vee \varepsilon^{\lceil K/2 \rceil}\right)\end{equation}
    for some absolute constant $c>0$. 
\end{proposition}
The proof of \Cref{LB-Gaussian} is deferred to \Cref{sec:homogeneous-LB-proofs}. Before outlining the key elements of the proof, we briefly compare \Cref{LB-Gaussian} with \Cref{thm-Gaussian-mean-minimax-rate}. Due to the assumption in \eqref{Gaussian-distribution} that $\theta \in [-1,1]$, the minimax risk is of constant order even when the privacy parameter $\alpha\leq n^{-1}K^{-1/2}$. The additional term $\varepsilon^{\lceil K/2 \rceil}$ follows from considering the event that at least half of the servers are contaminated. This is dominated by the other terms, i.e. 
\[\varepsilon^{\lceil K/2 \rceil}\leq \frac{1}{nK}\vee\frac{1}{n^2 K \alpha^2}\vee \frac{\varepsilon^2}{n}\]
only if $K\log(1/\varepsilon)\geq c'\log(K(n\wedge n^2\alpha^2))$, where $c'>0$ is an absolute constant. Therefore, we require that the number of batches $K$ is not too small, as enforced by condition \eqref{UB-condition-thm1} (up to the $\log(1/\varepsilon)$ factor), in order for the rate stated in \eqref{LB-Gaussian-eq} to be equal to that of \eqref{Eq:HomogeneousGaussianRate} in \Cref{thm-Gaussian-mean-minimax-rate}.

At a high level, the proof is based on the observation that $\mathcal{R}_\alpha(\varepsilon) \geq \mathcal{R}_\alpha(0) \vee \mathcal{R}_{\infty}(\varepsilon)\geq \mathcal{R}_\infty(0)$, such that minimax risk satisfies
\[
    \mathcal{R}_\alpha(\varepsilon) \geq \mathcal{R}_\infty(0)\vee \mathcal{R}_\alpha(0) \vee \mathcal{R}_{\infty}(\varepsilon),
\]
where $\alpha = \infty$ in the subscript denotes the minimax risk without any privacy constraints. In the absence of contamination and privacy constraints, the problem can be viewed as standard mean estimation with $nK$ i.i.d.~samples, corresponding to the standard lower bound of $1/(nK)$ for $\mathcal{R}_\infty(0)$.

To prove a lower bound on $\mathcal{R}_\infty(\varepsilon)$, we note that two different inlier distributions with TV distance less than $\varepsilon$ are indistinguishable under the contamination model. We establish the lower bound by constructing two $n$-fold product distributions of Gaussians with mean separation of order at least $\varepsilon/\sqrt{n}$, combined with an application of Le Cam's two point method \citep[e.g.][]{LeCam1973convergence,Yu1997assouad}. This involves bounding the TV modulus of continuity, as is standard in robust statistics \citep[e.g.][]{Chen2016a,Chen2018robust,Li2023on}. 

The lower bound on $\mathcal{R}_{\alpha}(0)$ is again proved using Le Cam's two point method. 
We combine the coupling method \citep[e.g.~Lemma 6.1 in][]{karwa2017finite} for CDP with derivation techniques for LDP data processing inequalities \citep[e.g.~the proof of Corollary 3 in][]{Duchi2018minimax}. This results in a tight upper bound on the TV distance between the privatised distributions. 

\subsubsection{Upper bound}\label{sec-Gaussian-mean-estimation-UB}

We now propose a two-step estimator in \Cref{alg:homogeneous-Gaussian-mean-estimation} that will later be shown to be minimax rate optimal up to poly-logarithmic factors in \Cref{UB-thm}. One notable feature of \Cref{alg:homogeneous-Gaussian-mean-estimation} is its adaptivity to the contamination level $\varepsilon$.  
At a high level,
we split the batches into two halves, and perform a localisation step using the first half, followed by a refinement step using the remaining half.  In more detail, it proceeds as follows.
\begin{enumerate}
    \item (Localisation.) Partition the interval $[-1,1]$ into bins of width $r$. Each batch in the first half assigns a smoothed vote taking a value in $[0,2r]$ to each bin, where each vote is subsequently privatised using the Laplace mechanism. We then select the bin with the largest sum of privatised votes and use its location $\mathfrak{A}^*$ as an initial estimator of the mean. The smoothed votes are robust to outliers, and, with high probability, the distance between the true mean $\theta$ and the initial estimator $\mathfrak{A}^*$ is of order $r$. See \Cref{alg:mean-localisation}.
    \item (Refinement.) Given a truncation parameter $M$, we use the remaining half to refine the estimator within the interval $[\mathfrak{A}^*-M, \mathfrak{A}^*+M]$, leading to the final estimate $\hat{\theta}$. See \Cref{alg:mean-refinement}.
\end{enumerate}
Notably, our privacy mechanism is interactive only through the output $\mathfrak{A}^*$ of \Cref{alg:mean-localisation}; communication between servers therefore remains limited.

In the formal statements of our algorithms, we recall that $\{X_{k,i}\}_{i=1}^{n}$ denotes the $k$th batch of data. For simplicity, we only consider the case that $K$ is even and use the sets $B=[K/2]$ and $B'=[K]\setminus [K/2]$ to denote the collections of batches used in the localisation step and the refinement step respectively. When $K$ is odd, we may discard one batch at random before proceeding with \Cref{alg:homogeneous-Gaussian-mean-estimation}, without affecting the rate in \Cref{UB-thm}. Although only the homogeneous setting is considered in this section, we present Algorithms~\ref{alg:mean-localisation} and \ref{alg:mean-refinement} to allow heterogeneous batch sample sizes and privacy budgets for later generalisation.

\begin{algorithm}
    \caption{FDP mean estimation with homogeneous servers}
    \label{alg:homogeneous-Gaussian-mean-estimation}
    \begin{algorithmic}[1]
        \INPUT {Data $\{\{X_{k,i}\}_{i=1}^{n}\}_{k=1}^{K}$, privacy parameter $\alpha>0$, bin width $r>0$, tuning parameters $T>0$, $M>0$}
        \State Set $\alpha_k \leftarrow \alpha$ for each $k \in \{1,\dots,K\}$.
        \State Set $\mathfrak{A}^* \leftarrow$ {\sc MeanLocalisation}($\{\{X_{k,i}\}_{i=1}^{n}\}_{k=1}^{K/2}$, $\{\alpha_k\}_{k=1}^{K/2}$, $r$, $T$). \Comment{\Cref{alg:mean-localisation}}
        \State Set $\hat{\theta}\leftarrow$ {\sc MeanRefinement}($\{\{X_{k,i}\}_{i=1}^{n}\}_{k=K/2+1}^{K}$, $\mathfrak{A}^*$, $\{\alpha_k\}_{k=K/2+1}^{K}$, $M$, $T$) \Comment{\Cref{alg:mean-refinement}}
        \State \Return $\hat{\theta}$
    \end{algorithmic}
\end{algorithm}

\begin{algorithm}
    \caption{\textsc{MeanLocalisation}($\{\{X_{k,i}\}_{i=1}^{n_k}\}_{k\in B}$, $\{\alpha_k\}_{k\in B}$, $r$, $T$)}
    \label{alg:mean-localisation}
    \begin{algorithmic}[1]
        \INPUT{Data $\{\{X_{k,i}\}_{i=1}^{n_k}\}_{k\in B}$, privacy parameters $\{\alpha_k\}_{k\in B}$}, bin width $r>0$, tuning parameter $T>0$
        
        \For{$j \in [\lceil 2r^{-1} \rceil ]$}
        \State Generate $\{W_{k,j}\}_{k\in B}$ with $W_{k,j}\overset{\mathrm{i.i.d.}}{\sim} \mathrm{Laplace}(0,1)$.
        \State Set
        \[Z_j \leftarrow  \sum_{k \in B}\left[\left(2r - \left|\frac{1}{n_k}\sum_{i=1}^{n_k} [X_{k,i}]_{-T}^T - \left\{-1+(j-1/2)r\right\}\right|\right)_+ + \frac{16 T}{\alpha_k n_k} W_{k,j}\right] \]
        \EndFor
        \State Set $j^* \leftarrow \min\{\argmax_{j\in [\lceil 2r^{-1} \rceil ]} Z_j\}$. 
        \State \Return $\mathfrak{A}^*\leftarrow -1 + (j^*-1/2)r$ 

    \end{algorithmic}
\end{algorithm}

\begin{algorithm}
    \caption{\textsc{MeanRefinement}($\{\{X_{k,i}\}_{i=1}^{n_k}\}_{k\in B'}$, $\mathfrak{A}^*$, $\{\alpha_k\}_{k\in B'}$, $M$, $T$)}
    \label{alg:mean-refinement}
    \begin{algorithmic}[1]
        \INPUT{Data \{$\{X_{k,i}\}_{i=1}^{n_k}\}_{k\in B'}$, initial estimator $\mathfrak{A}^*$, privacy parameters $\{\alpha_k\}_{k\in B'}$, tuning parameters $M>0$, $T>0$}
        \State Generate $\{L_k\}_{k\in B'}$ with $L_k \overset{\mathrm{i.i.d.}}{\sim}\mathrm{Laplace}(0,1)$.
        \State Set \[\hat{\theta} \leftarrow \frac{1}{|B'|}\sum_{k\in B'} \left(\mathfrak{A}^*+\left[\frac{1}{n_k}\sum_{i = 1}^{n_k} [X_{k,i}]_{-T}^T -\mathfrak{A}^*\right]_{-M}^M + \frac{2T}{\alpha_k n_k}L_k\right)\]
        \State \Return $\hat{\theta}$
    \end{algorithmic}
\end{algorithm}

In both Algorithms~\ref{alg:mean-localisation} and~\ref{alg:mean-refinement}, the quantity
\begin{equation}\label{key-summary-stat}\frac{1}{n_k}\sum_{i = 1}^{n_k} [X_{k,i}]_{-T}^T\end{equation}
with truncation parameter $T$ can be viewed as the key summary statistic of the $k$th batch used to build our estimator of $\theta$. This is in the spirit of classical robust mean estimation methods such as the trimmed mean or truncated mean \citep[e.g.][]{bickel1965some} and their recent variants \citep[e.g.][]{bubeck2013bandits,Devroye2016sub-gaussian,lugosi2021robust}; versions of these statistics were shown to be optimal for both heavy-tailed mean estimation and robustness against strong adversarial contamination \citep[see e.g.][]{lugosi2021robust}. While the truncated mean is rooted in classical robust statistics, the key purpose of truncation by parameter $T$ of our procedure in the Gaussian setting is instead to control the sensitivity of the truncated mean and limit the magnitude of privatisation noise added by the Laplace mechanism. The robustness aspect will be explored as we extend to heavy-tailed distributions in \Cref{sec-finite-moments}.

Crucial in obtaining our preliminary estimator $\mathfrak{A}^*$, the localisation step of our procedure extends the method of private histograms \citep[see e.g.][]{chawla2012privacy,karwa2017finite, Li2023on}. While using binary votes is common in the privacy literature, we consider a novel version tailored to FDP constraints, where each batch casts smoothed votes taking values in $[0,2r]$. The value of the smoothed vote is continuous to changes in the truncated mean in \eqref{key-summary-stat}; this prevents the discrete jumps from a single data point change that lead to an $\ell_1$-sensitivity of constant order in binary voting. Consequently, the smoothed votes' $\ell_1$-sensitivity of order $T/n$ can be optimally tuned both using the parameter $T$ and according to the batch sample size $n$. Meanwhile, batch-level robustness is ensured by limiting both the number and magnitude of non-zero votes. 

Overall, the final estimator $\hat{\theta}$ can be viewed as a private truncated mean of truncated means (see \Cref{alg:mean-refinement}), where the truncation interval $[\mathfrak{A}^*-M, \mathfrak{A}^*+M]$ is calibrated in the preliminary localisation step (\Cref{alg:mean-localisation}). This truncated mean approach involving a preliminary calibration of the truncation interval resembles that of the heavy-tailed estimation procedure detailed in Section~7 in \cite{Devroye2016sub-gaussian}. More broadly, this localisation-refinement approach is an example of multi-step strategies that are common in both robust statistics \citep[e.g.~outlier filtering algorithms in][]{diakonikolas2019robust, prasad2019unified, diakonikolas2019recent} and privacy literature \citep[e.g.][]{karwa2017finite}, as well as recent work on simultaneously robust and private algorithms under CDP \citep[e.g.][]{Kamath20private, Liu2021robust} and LDP \citep[e.g.][]{Li2023on}.

We proceed to present the theoretical guarantees of \Cref{alg:homogeneous-Gaussian-mean-estimation} in \Cref{UB-thm} below.

\begin{proposition}\label{UB-thm}
    Consider the homogeneous Gaussian mean estimation problem detailed in \eqref{Gaussian-distribution} and~\eqref{Gaussian-model}. There exists an absolute constant $c>0$ such that when the condition~\eqref{UB-condition-thm1} holds,
    the output estimator~$\hat{\theta}$ of \Cref{alg:homogeneous-Gaussian-mean-estimation} satisfies 
    \[\mathbb{E}[|\hat{\theta}-\theta|^2]
            \leq C\left\{\frac{1}{nK}+\frac{\varepsilon^2 T^2}{n}+\frac{T^2}{n^2 K\alpha^2}(1+\varepsilon^2T^2)\right\}\]
    with input tuning parameters \begin{equation}\label{UB-thm-tuning-params}
        T =  C_1 \left\{1\vee \sqrt{\log(K(n \wedge n^2 \alpha^2))}\right\}, \quad r = C_2 T\left(\frac{1}{n}\vee \frac{T^2}{n^2K\alpha^2}\right)^{1/2} \quad \text{and} \quad M=3r,
    \end{equation}
    where $C, C_1, C_2>0$ are absolute constants. 
    
\end{proposition}

Comparing the rate achieved by \Cref{alg:homogeneous-Gaussian-mean-estimation} with \Cref{LB-Gaussian}, there is a poly-logarithmic gap of at most $T^4 \lesssim 1\vee \log^2(nK)$, thus our rate is optimal up to poly-logarithmic factors in the total sample size $nK$. We note that this near optimality is obtained without any knowledge of the contamination level $\varepsilon$.

For the choice of tuning parameters, the truncation parameter $T$ is chosen to balance the bias introduced by truncation against the variance arising from the Laplace privatisation mechanism. Since the inlier distribution is assumed to be Gaussian, taking $T$ to scale logarithmically ensures that, with high probability, no truncation occurs for the inlier observations, limiting its contribution to the bias.  The resulting inflation in the variance term remains logarithmic.  
As for the bin width, we require $r\asymp_{\log}\{n^{-1/2}\vee (n^2K\alpha^2)^{-1/2}\}$, where the two terms correspond, respectively, to the high-probability concentration rate of $n^{-1}\sum_{i=1}^n[X_{k,i}]_{-T}^T$ for inlier observations, and to the scale of the Laplace noise added for privatisation.   Finally, the truncation level $M$ in \Cref{alg:mean-refinement} is chosen so that the interval $[\mathfrak{A}^*-M,\mathfrak{A}^*+M]$ contains the true mean with high probability, taking into account the estimation accuracy of $\mathfrak{A}^*$.  Since the condition~\eqref{UB-condition-thm1} ensures that the $\mathfrak{A}^*$ remains within a small number of bins of $\theta$ with high probability, $M$ is chosen as a constant multiple of the bin width $r$. 

To see how the choice of tuning parameters lead to the stated rate of convergence, we note that the localisation step outputs $\mathfrak{A}^*$ with a squared error rate of convergence of 
\[r^2 \asymp_{\log} \frac{1}{n}\vee \frac{1}{n^2K\alpha^2}.\]
The rate of convergence of the final estimator $\hat{\theta}$ is then refined to 
\[r^2 \left(\varepsilon^2\vee \frac{1}{K}\right) +\frac{T^2}{n^2K\alpha^2}\asymp_{\log} \frac{1}{nK}\vee\frac{\varepsilon^2}{n} \vee \frac{1}{n^2K\alpha^2}.\]

\subsection{Heavy-tailed mean estimation}\label{sec-finite-moments}
In this section, we consider the homogeneous estimation problem with a heavy-tailed inlier distribution. Specifically, we consider the class of distributions with a bounded $v$th central moment
\begin{equation}\label{finite-moments-class}
    \mathcal{P}=\mathcal{P}_v= \{P:\theta = \mathbb{E}_{Y\sim P}[Y] \in [-1,1], \,\mathbb{E}_{Y\sim P}[|Y-\theta|^v]\leq 1\}.
\end{equation}
We view the parameter $v>1$ as a fixed but arbitrary constant.
We denote the data set $D_k$ of size $n$ stored within the $k$th server by $D_k = \{X_{k,i}\}_{i=1}^n$. Then, under our contamination model  $\mathcal{P}_\varepsilon(\mathcal{P})$, 
\begin{equation}\label{finite-moments-model}
    X_{k}= (X_{k,1},\dots,X_{k,n})\overset{\mathrm{ind.}}{\sim} (1-\varepsilon)P^{\otimes n}+\varepsilon G_k,
\end{equation}
where $P\in \mathcal{P}$ and $G_k$ is any distribution on $\mathbb{R}^n$. We establish the minimax rate of the homogeneous heavy-tailed mean estimation problem in the following theorem.
\begin{theorem}\label{thm-heavy-tail-minimax-rate}
    Consider the homogeneous heavy-tailed mean estimation problem detailed in \eqref{finite-moments-class} and~\eqref{finite-moments-model}. There exists an absolute constant $c>0$ such that when the condition \eqref{UB-condition-thm1} holds, the minimax risk defined in (\ref{Gaussian-minimax-risk}) satisfies
    \[\mathcal{R}_\alpha(\varepsilon) \asymp_{\log} \frac{1}{nK}\vee \frac{1}{(nK)^{2-2/v}}\vee\frac{1}{(n^2K\alpha^2)^{1-1/v}}\vee\frac{\varepsilon^2}{n}\vee\frac{\varepsilon^{2-2/v}}{n^{2-2/v}}. \]
    In particular, when $v\in (1,2]$,
    \[\mathcal{R}_\alpha(\varepsilon)\asymp_{\log}  \frac{1}{(nK)^{2-2/v}}\vee\frac{1}{(n^2K\alpha^2)^{1-1/v}}\vee\frac{\varepsilon^{2-2/v}}{n^{2-2/v}},\]
    and when $v\geq 2$,
    \[\mathcal{R}_\alpha(\varepsilon)\asymp_{\log}\frac{1}{nK}\vee\frac{1}{(n^2K\alpha^2)^{1-1/v}}\vee\frac{\varepsilon^{2-2/v}}{n^{2-2/v}}\vee\frac{\varepsilon^2}{n}.\]
\end{theorem}
\Cref{thm-heavy-tail-minimax-rate} characterises the minimax rate with the lower bound in \Cref{finite-moments-LB} and the upper bounds in \Cref{finite-moment-UB,infinite-variance-UB}, respectively. Similarly to \Cref{thm-Gaussian-mean-minimax-rate}, the upper bounds require condition \eqref{UB-condition-thm1} to hold. Under the heavy-tailed setting of \eqref{finite-moments-class}, we observe a phase transition boundary at $v=2$, separating the rates in the settings of finite variance $(v\geq 2)$ and infinite variance $(1<v<2)$, respectively. The minimax rate under the Gaussian setting is recovered as $v\rightarrow \infty$. In all regimes, we again observe the separation of the costs of privacy and contamination found in the Gaussian setting. We proceed to comment on the three regimes corresponding to absence and presence of privacy and contamination.  

In the absence of both privacy constraints and contamination, we recover the rate $(nK)^{-1} \vee (nK)^{-(2-2/v)}$. It is well-established that, under the setting of bounded variance ($v\geq 2$), the minimax rate is characterised by the parametric rate. Consistent with our result, a version of the lower bound for $v\in (1,2]$ is established in \cite{Devroye2016sub-gaussian}. Results from \cite{bubeck2013bandits} show that the truncated mean and the median of means attain the optimal rate for $v\in (1,2]$.  

Under the FDP constraint, the cost of privacy is characterised by the term $(n^2K\alpha^2)^{-(1-1/v)}$ for all $v\in(1,\infty)$. Considering the special cases of $n=1$ and $K=1$, our result interpolates the LDP rate $(K\alpha^2)^{-(1-1/v)}$ \citep[Proposition 4 and Theorem 5 in][]{Li2023on} and the CDP rate $(n^2\alpha^2)^{-(1-1/v)}$ \cite[Theorem 1.1 in][proved for $v \geq 2$]{Kamath20private}. Although heavy-tailedness does not affect the rate in the non-private and non-contaminated setting when $v\geq 2$, the factor of $1-1/v$ present in the exponent of both $n$ and $K$ shows that the cost of privacy is materially altered by heavy-tailedness even for $v\geq 2$. 

The effect of batch contamination is characterised by the rate $(\varepsilon/n)^{2-2/v}\vee \varepsilon^2/n$ and we see that this grows progressively with heavy-tailedness. Recent work \citep{aliakbarpour2026high} studying multivariate mean estimation with untrusted batches is closely related to ours, where heavy-tailedness is studied under a bounded covariance assumption (corresponding to $v=2$). We extend this to a bounded $v$th moment assumption with arbitrary $v\in (1,\infty)$ in the univariate case, and the rate $\varepsilon/n$ that we recover for the case $v =2$ is consistent with their findings. In the case $n = 1$, the problem reduces to heavy-tailed mean estimation \eqref{finite-moments-class} under Huber contamination \eqref{Huber-model}, where the corresponding cost of contamination is $\varepsilon^{2-2/v}$ (see e.g.~Lemma 3 in \citealp{prasad2019unified}; Lemma 1.12 in \citealp{Diakonikolas2023algorithmic}). Our results show that when $v \geq 2$ there are two regimes bridging the Gaussian rate of $\varepsilon^2/n$ and the infinite-variance rate of $(\varepsilon/n)^{2-2/v}$, with the dominant rate transitioning between the two at $\varepsilon\asymp n^{1-v/2}$. To the best of our knowledge, \Cref{thm-heavy-tail-minimax-rate} is the first result for heavy-tailed mean estimation under batch contamination with the full range of $v\in (1,\infty)$. As seen in the proofs of our upper and lower bounds, we may view the term $\varepsilon^2/n$ as the effect of contamination in the main bulk of the data distribution and $(\varepsilon/n)^{2-2/v}$ as the effect of contamination in the tails, which worsens as $v$ decreases.

\subsubsection{Lower bound}
We proceed to present the following lower bound on the minimax risk.
\begin{proposition}\label{finite-moments-LB}
    Consider the homogeneous heavy-tailed mean estimation problem detailed in \eqref{finite-moments-class} and \eqref{finite-moments-model}. The minimax risk defined in \eqref{Gaussian-minimax-risk} satisfies
    \[\mathcal{R}_\alpha(\varepsilon)\geq c_v\left(\frac{1}{nK}\vee \frac{1}{(nK)^{2-2/v}}\vee\frac{1}{1\vee(n^2K\alpha^2)^{1-1/v}}\vee\frac{\varepsilon^2}{n}\vee\frac{\varepsilon^{2-2/v}}{n^{2-2/v}}\vee \varepsilon^{\lceil K/2 \rceil}\right)\]
    for some constant $c_v>0$ depending only on $v$.
\end{proposition}
The proof of \Cref{finite-moments-LB} is deferred to \Cref{sec:homogeneous-LB-proofs}. As discussed in \Cref{LB-Gaussian}, we focus on the regime that $\varepsilon^{\lceil K/2 \rceil}$ is not the dominant term in the rate. Note that we have
\[\varepsilon^{\lceil K/2 \rceil}\leq \frac{1}{nK}\vee \frac{1}{(nK)^{2-2/v}}\vee\frac{1}{(n^2K\alpha^2)^{1-1/v}}\vee\frac{\varepsilon^2}{n}\vee\frac{\varepsilon^{2-2/v}}{n^{2-2/v}}\]
only if the same condition $K\log(1/\varepsilon)\geq c_v'\log(K(n\wedge n^2\alpha^2))$ holds, where $c'_v>0$ is a constant depending only on $v$. The lower bound proof follows the general strategy outlined in the discussion following \Cref{LB-Gaussian}. The $\varepsilon^{\lceil K/2 \rceil}$ term follows by the same argument, and the construction of $n$-fold product distributions of Gaussians (after a constant rescaling of the variance to satisfy \eqref{finite-moments-class}) is used to prove the rates $(nK)^{-1}$ and $\varepsilon^2/n$ that also appear in \Cref{LB-Gaussian}. The lower bound proof strategy of constructing two-point distributions has been used in the robust statistics and differential privacy literature studying heavy-tailed distributions \citep[e.g.][]{Devroye2016sub-gaussian, Kamath20private, Li2023on}. Adapting this strategy, we establish the lower bounds $(nK)^{-(2-2/v)}$, $(n^2K\alpha^2)^{-(1-1/v)}$ and~$(\varepsilon/n)^{2-2/v}$ that capture the effect of tail events under the settings of FDP and batch contamination.

\subsubsection{Upper bound}
To prove upper bounds on the minimax risk, we show that, with the input tuning parameters adapted to the heavy-tailed setting, \Cref{alg:homogeneous-Gaussian-mean-estimation} can be directly applied and achieves optimal rates up to poly-logarithmic factors. \Cref{finite-moment-UB,infinite-variance-UB} below establish upper bounds on the minimax risk for the finite variance setting ($v\geq 2$) and the infinite variance setting ($v\in (1,2)$) respectively.

\begin{proposition}\label{finite-moment-UB}
    Consider the homogeneous heavy-tailed mean estimation problem detailed in \eqref{finite-moments-class} and~\eqref{finite-moments-model} with $v\geq 2$. There exists an absolute constant $c>0$ such that when the condition~\eqref{UB-condition-thm1} holds,
    the output estimator $\hat{\theta}$ of \Cref{alg:homogeneous-Gaussian-mean-estimation} satisfies
\[\E [|\hat{\theta}-\theta|^2 ]\leq C_v\left[\frac{1}{nK}+ \frac{\varepsilon^{2-2/v}}{n^{2-2/v}}+\left(\frac{\varepsilon^2}{n} +\frac{1}{(n^2 K\alpha^2)^{1-1/v}}\right)\{1\vee\log(K(n\wedge n^2\alpha^2))\}\right]\]
with input tuning parameters 
    \begin{equation}\label{finite-var-ub-t-choice}T = A_1 (n^2K\alpha^2)^{1/(2v)}
    , \quad 
    r = A_2\left(\left[1\vee\sqrt{\log(K(n\wedge n^2\alpha^2))}\right]\left\{\frac{1}{\sqrt{n}}\vee T^{1-v}\right\}\vee \frac{\varepsilon^{-1/v}\wedge K^{1/v}}{n^{1-1/v}}\right),\end{equation}
and 
\begin{equation}\label{finite-var-ub-m-choice}M=A_3 r,\end{equation}where $C_v, A_1, A_2, A_3>0$ are constants that depend only on $v$. 
\end{proposition}
 
\begin{proposition}\label{infinite-variance-UB}
    Consider the homogeneous heavy-tailed mean estimation problem detailed in \eqref{finite-moments-class} and~\eqref{finite-moments-model} with $v\in(1,2)$. There exists an absolute constant $c>0$ such that when the condition~\eqref{UB-condition-thm1} holds,
    the output estimator $\hat{\theta}$ of \Cref{alg:homogeneous-Gaussian-mean-estimation} satisfies
    \[\E [|\hat{\theta}-\theta|^2 ]\leq C_v\log^2(nK) \left\{ \frac{1}{(nK)^{2-2/v}}\vee\frac{1}{(n^2K\alpha^2)^{1-1/v}}\vee\frac{\varepsilon^{2-2/v}}{n^{2-2/v}}\right\}\] with input tuning parameters 
    \begin{equation}\label{infinite-var-ub-t-choice}
        T = \left(nK\wedge n\alpha\sqrt{K}\wedge  n/\varepsilon\right)^{1/v},
    \end{equation} 
    \begin{equation}\label{infinite-var-ub-r-choice}r =  B_1\left(T^{1-v}\left\{1\vee\sqrt{\log(nK\wedge n\alpha\sqrt{K})}\right\}\vee \frac{\varepsilon^{-1/v}\wedge K^{1/v}}{n^{1-1/v}}\log(K/\delta)\right),\end{equation}
    \begin{equation}\label{infinite-var-m-choice}M = B_2 r, \quad \text{where} \quad \delta = \frac{12}{B_1} \left(\frac{1}{(nK)^{2-2/v}} \vee  \frac{1}{(n^2 K\alpha^2)^{1-1/v}}\right),\end{equation}
    where $C_v,B_1,B_2>0$ are constants that depend only on $v$.

\end{proposition}

The proofs of \Cref{finite-moment-UB,infinite-variance-UB} are deferred to \Cref{appendix:homogeneous-UB-proofs}. Comparing with \Cref{finite-moments-LB}, \Cref{alg:homogeneous-Gaussian-mean-estimation} is rate-optimal up to a factor of at most $\log^2(nK)$ with appropriate tuning parameters. We proceed to discuss the differences in the choices of tuning parameters between the Gaussian setting and the heavy-tailed setting below, and provide an overview of how the tuning parameters affect the rate achieved by \Cref{alg:homogeneous-Gaussian-mean-estimation}.

Our upper bound can be interpreted as
\[\underbrace{\vphantom{\frac{\varepsilon M^2}{K}}T^{2-2v}+(nM)^{2-2v}}_{\text{Bias}^2}+\underbrace{\frac{1 \vee T^{2-v}}{nK}}_{\text{Variance}}+\underbrace{\varepsilon^2M^2+\frac{\varepsilon M^2}{K}}_{\text{Contamination}}+\underbrace{\frac{T^2}{n^2K\alpha^2}}_{\text{Privatisation}}.\]
In the finite variance setting ($v\geq 2$), as in the earlier Gaussian setting, the truncation parameter $T\geq 1$ is only required to balance the variance due to privatisation and the bias induced by truncating inlier samples, and we see that $T\asymp (n^2K\alpha^2)^{1/(2v)}$ suffices. Note that, due to the heavy tails of the inlier distribution, this choice of $T$ is considerably larger than that in the Gaussian setting to avoid truncating an excessive number of inlier samples. As discussed in the Gaussian setting, the parameter $r^2$ again must be at least the variance and squared bias of the summary statistic $n^{-1}\sum_{i=1}^n[X_{k,i}]_{-T}^T$ of each server, as well as the variance due to privatisation. \Cref{trunc-lemma-finite-moments} shows that the bias contribution due to truncation by $T$ is upper bounded by $T^{1-v}$. In the $v\geq 2$ case, this enforces $M^2 \gtrsim r^2 \gtrsim n^{-1}\vee T^{2-2v}$. Further, we use a version of the Fuk--Nagaev inequality \citep[e.g.~Equation (1) in][]{jirak2025concentration} to bound the probability of truncation by the parameter $M$, leading to the squared bias term of $(nM)^{2-2v}$. Finally, the terms $(nM)^{2-2v}$, $\varepsilon^2M^2$ and $\varepsilon M^2/K$ need to be balanced appropriately, for which the choice of \[M\asymp r \asymp_{\log} \frac{1}{\sqrt{n}} \vee T^{1-v} \vee \frac{\varepsilon^{-1/v}\wedge K^{1/v}}{n^{1-1/v}}\]
is one that leads to optimal rates, up to poly-logarithmic factors, for the finite variance ($v\geq 2$) setting.

In the infinite variance setting ($v\in (1,2)$), the $v$th moment bound on the inlier distribution is no longer sufficient in upper bounding the variance and tail probabilities of the truncated mean $n^{-1}\sum_{i=1}^n[X_{k,i}]_{-T}^T$, for which we now rely on the parameter $T$. Details of this are given in \Cref{trunc-lemma-infinite-variance}. For tight upper bounds, we choose the smallest possible $T$, subject to the squared bias term $T^{2-2v}$ not exceeding the rate stated in the lower bound of \Cref{finite-moments-LB}. This results in the choice given in~\eqref{infinite-var-ub-t-choice}. Note that this is in contrast to the previous settings, where the main role of the truncation is to upper bound the sensitivity for the purposes of privatisation. For the refinement step, we use Bernstein's inequality \citep[e.g.~Proposition 2.14 in][]{Wainwright2019high-dimensional} to bound the probability of truncating by the parameter $M$. For this to be a small probability event given an inlier observation, we require, after taking into account the bias, $M^2\gtrsim r^2\gtrsim (T/n)^2 \vee T^{2-v}/n \vee T^{2-2v}$. The choice of 
\[M\asymp r\asymp_{\log} T^{1-v} \vee \frac{\varepsilon^{-1/v}\wedge K^{1/v}}{n^{1-1/v}}\]
is one that satisfies this and leads to the optimal rate up to poly-logarithmic factors.

Notably, we have shown that \Cref{alg:homogeneous-Gaussian-mean-estimation} is an FDP procedure that is optimally robust against both batch contamination and heavy-tailed distributions with bounded $v$th moment for arbitrary $v\in (1,\infty)$. This builds on a range of recent work \citep[e.g.][]{prasad2019unified, diakonikolas2020outlier, lugosi2021robust, abdalla2024covariance} that have developed algorithms robust against both contamination and heavy-tailed distributions, unveiling a deep connection between contamination and heavy-tailed models. \cite{cherapanamjeri2025heavy} showed more generally that optimal procedures robust against adversarial contamination are also optimal for heavy-tailed distributions. To the best of our knowledge, our method is the first for an arbitrary number of moments $v\in (1,\infty)$ under batch contamination, suggesting that some of these connections may well extend to the batch contamination model.

\section{Gaussian mean estimation with heterogeneous servers}\label{sec-mean-est-hetero}
In this section, we consider the Gaussian mean estimation problem, with the class of distributions of interest $\mathcal{P}$ in~\eqref{Gaussian-distribution} and with heterogeneous batch sample sizes and privacy budgets.
For $k \in [K]$, server~$k$ holds a data set $D_k= \{X_{k,i}\}_{i=1}^{n_k}$ of size $n_k$ with
\begin{equation}\label{heterogeneous-model}
    X_k = (X_{k,1},\dots,X_{k,n_k}) \overset{\mathrm{ind.}}{\sim} P_{\varepsilon,k} =  (1 - \varepsilon) \mathcal{N}^{\otimes n_k}(\theta, 1) + \varepsilon G_k,
\end{equation}
with $\varepsilon \in [0, 1/2)$ and $G_k$ being an arbitrary distribution on $\mathbb{R}^{n_k}$.  Each server $k \in [K]$ has privacy budget $\alpha_k > 0$, jointly satisfying $\{\alpha_k\}_{k = 1}^K$-FDP as defined in \Cref{FDP-def}. Note that we do not require the batch sizes $n_k$ to be private.

To understand the fundamental limit of this problem, we partition the servers into groups of similar batch sample sizes and privacy constraints, such that all batches belonging to the same group are nearly homogeneous. Specifically, denote $n_{\max} = \max_{k\in [K]} n_k $, $n_{\min}= \min_{k\in [K]} n_k$, $\alpha_{\max} = \max_{k\in [K]} \alpha_k$ and $\alpha_{\min}= \min_{k\in [K]} \alpha_k$. For each $l_1, l_2 \in \mathbb{N}$, let 
\begin{equation}\label{server-group-def}\mathcal{I}_{l_1,l_2} = \{k \in [K] : n_k/n_{\min} \in [2^{l_1-1},2^{l_1}), \alpha_k/\alpha_{\min} \in [2^{l_2-1},2^{l_2}) \}.\end{equation}
This partitions the servers into at most $\{1+\lfloor\log_2(n_{\max}/n_{\min})\rfloor\}\{1+\lfloor\log_2(\alpha_{\max}/\alpha_{\min})\rfloor\}$ groups, such that
\[\bigcup_{l_1=1}^{1+\lfloor\log_2(n_{\max}/n_{\min})\rfloor} \bigcup_{l_2=1}^{1+\lfloor\log_2(\alpha_{\max}/\alpha_{\min})\rfloor}\mathcal{I}_{l_1,l_2} = [K].\]
 
For ease of comparison with our homogeneous results in \Cref{sec-Gaussian-mean}, let $\tilde{n}_{l_1}$ and $\tilde{\alpha}_{l_2}$ denote the representative batch sample size and privacy budget of the group $\mathcal{I}_{l_1,l_2}$ respectively, and let $K_{l_1,l_2}$ denote the number of batches belonging to the group $\mathcal{I}_{l_1,l_2}$, i.e.
\[\tilde{n}_{l_1} = 2^{l_1}n_{\min}, \quad \tilde{\alpha}_{l_2} = 2^{l_2} \alpha_{\min}\quad \text{and} \quad K_{l_1,l_2} = |\mathcal{I}_{l_1,l_2}|.\]
We further define the representative rate of group $\mathcal{I}_{l_1, l_2}$ by
\[R(\varepsilon,l_1,l_2) = \frac{1}{\tilde{n}_{l_1}K_{l_1,l_2}}\vee\frac{1}{\tilde{n}_{l_1}^2\tilde{\alpha}_{l_2}^2K_{l_1,l_2}}\vee \frac{\varepsilon^2}{\tilde{n}_{l_1}}.\]
Recalling the minimax rate 
\begin{equation} \label{eq-words-rate}
    \frac{1}{\mbox{per-server sample size}} \times \left\{\frac{1}{\mbox{\#.~server}} \vee \frac{1}{\mbox{per-server sample size} \times \mbox{\#.~server} \times \alpha^2} \vee \varepsilon^2\right\},
\end{equation}
characterised by \Cref{thm-Gaussian-mean-minimax-rate} up to poly-logarithmic factors, note that $R(\varepsilon,l_1,l_2)$ denotes the rate in the homogeneous special case, where the data set contains $K_{l_1,l_2}$ batches, all having the same sample size $\tilde{n}_{l_1}$ and privacy budget~$\tilde{\alpha}_{l_2}$. Using this batch partition framework, we establish the minimax rate of the Gaussian mean estimation problem with heterogeneous batch sample sizes and privacy budgets in the following theorem.

\begin{theorem}\label{thm-hetero-rate}
    Consider the contamination model with heterogeneous server sizes detailed in~\eqref{heterogeneous-model}. Suppose that there exist $l_1^*,l_2^*\in \mathbb{N}$ that satisfy 
    \begin{equation}\label{eq:optimal-group-def}
    R(\varepsilon,l_1^*,l_2^*) \asymp \min_{l_1,l_2\in \mathbb{N}: \mathcal{I}_{l_1,l_2} \neq \emptyset} R(\varepsilon,l_1,l_2), \quad \mathcal{I}_{l_1^*,l_2^*}\neq \emptyset, 
    \end{equation}
    and the condition that
    \begin{equation}\label{non-adaptive-thm-condition}
    \varepsilon+\sqrt{\frac{\log\{K_{l_1^*,l_2^*}(\tilde{n}_{l_1^*}\wedge \tilde{n}_{l_1^*}^2\tilde{\alpha}_{l_2^*}^2)\}}{K_{l_1^*,l_2^*}}}+ \frac{1}{\sqrt{K_{l_1^*,l_2^*}(\tilde{n}_{l_1^*}\wedge\tilde{n}_{l_1^*}^2\tilde{\alpha}_{l_2^*}^2)}}<c,
    \end{equation}
    where $c>0$ is an absolute constant.
    Then, the minimax risk defined in \eqref{general-minimax-risk} satisfies
    \begin{equation*}
        \mathcal{R}
        \asymp_{\log} \min_{\substack{l_1,l_2\in \mathbb{N}: \mathcal{I}_{l_1,l_2} \neq \emptyset}} \left\{R(\varepsilon,l_1,l_2)\right\}.
    \end{equation*}
    
\end{theorem}

\Cref{thm-hetero-rate} summarises the lower bound of \Cref{hetero-LB} in \Cref{sec-hetero-LB} and the upper bound of \Cref{non-adaptive-thm} in \Cref{sec-hetero-UB}. To interpret this result, consider the (nearly) homogeneous special case where only one group of batches $\mathcal{I}_{l_1,l_2}$ is non-empty, such that all batches have sample size of order $\tilde{n}_{l_1}$ and privacy parameter of order $\tilde{\alpha}_{l_2}$. \Cref{thm-hetero-rate} recovers the homogeneous rate of $R(\varepsilon,l_1,l_2)$ established in \Cref{thm-Gaussian-mean-minimax-rate}. Therefore, the minimax rate established in \Cref{thm-hetero-rate} can be viewed as the minimum homogeneous rate over all dyadic groups of servers. This indicates that the fundamental difficulty of the heterogeneous problem is governed, up to poly-logarithmic factors, by the single most informative dyadic group, i.e.~the minimiser of $R(\varepsilon,l_1,l_2)$. The condition \eqref{non-adaptive-thm-condition}, required by our proposed estimation procedure in \Cref{sec-hetero-UB}, thus echoes that of \eqref{UB-condition-thm1} in \Cref{thm-Gaussian-mean-minimax-rate}, except that it is now imposed on the most informative group $(l_1^*,l_2^*)$ only.

In the absence of contamination ($\varepsilon=0$), our result recovers the FDP rate $\{\sum_{k=1}^K (n_k \wedge n_k^2\alpha_k^2)\}^{-1}$ of the Gaussian mean estimation problem \citep[see e.g.][]{Li2024federated} up to poly-logarithmic factors. To see this, observe that
\begin{multline*}\left\{\log_2\left(\frac{2n_{\max}}{n_{\min}}\right)\log_2\left(\frac{2\alpha_{\max}}{\alpha_{\min}}\right)\right\}^{-1}\min_{l_1,l_2\in \mathbb{N}: \mathcal{I}_{l_1,l_2} \neq \emptyset} \left\{R(0,l_1,l_2)\right\}\\
\leq \frac{1}{\sum_{k=1}^K (n_k\wedge n_k^2\alpha_k^2)}\leq 16\min_{l_1,l_2\in \mathbb{N}: \mathcal{I}_{l_1,l_2} \neq \emptyset} \left\{R(0,l_1,l_2)\right\}.\end{multline*}
In this context, one can interpret the rate by viewing the effective sample size of the $k$th server to be $n_k \wedge n_k^2\alpha_k^2$, with an additive contribution from each server resulting in a total effective sample size of $\sum_{k=1}^K (n_k \wedge n_k^2\alpha_k^2)$. This additive phenomenon was found to hold for various other statistical tasks under the FDP constraint, such as nonparametric regression \citep[Corollary 3 in][]{cai2024optimal}, functional mean estimation \citep[e.g.~Theorem 4 in][]{xue2024optimal} and Cox regression \citep[Theorem 1 in][]{Hung2025optimal}.  In these works, weighted aggregation across servers was used to handle heterogeneity. In contrast, our proposed procedure uses only samples from the most informative dyadic group, as motivated by our batch partition framework. Our method benefits from reduced computational and data communication costs, while retaining nearly the same performance, see \Cref{sec-hetero-UB} below. 

Returning to the full problem where batch contamination is present, we discuss a key feature distinctive from the FDP setting without contamination. With the exception of the case $\varepsilon=0$, the rate $R(\varepsilon,l_1,l_2)$ does not converge to zero as the group size $|\mathcal{I}_{l_1,l_2}|=K_{l_1,l_2}\rightarrow \infty$. As a consequence, there exist scenarios where an arbitrarily small proportion of the data is sufficient to attain the minimax rate up to poly-logarithmic factors, as we will see in the example below. This serves as a key motivation for the focus on the single most informative group by our batch partition framework, especially in the regime where the number of batches is relatively large. To the best of our knowledge, \Cref{thm-hetero-rate} is the first result on batch contamination under full batch size heterogeneity. 

Another notable consequence of \Cref{thm-hetero-rate} is that heterogeneity can introduce interactions between the costs of privacy and contamination. This is in contrast to our homogeneous results (Theorems \ref{thm-Gaussian-mean-minimax-rate} and \ref{thm-heavy-tail-minimax-rate}), and the findings of \cite{Liu2021robust} under CDP and \cite{Li2023on} under LDP, where the costs disentangle.  The key distinction is that the group that is most informative under privacy constraints alone need not be the group that is most informative under contamination.  Consequently, the joint minimax risk need not be characterised by the maximum of the risks arising from privacy and contamination considered separately. 

We illustrate both of these features using the following example.

\medskip
\noindent \textbf{Example.} For simplicity, we consider the special case with only two groups of batches, each with the following batch sizes, privacy budgets and group sizes.
\begin{center}
    \setlength{\tabcolsep}{4pt}
    \begin{tabular}{c|c|c|c|c|c}
    Batch size & Group size & Privacy budget & Parametric rate &FDP rate& Contamination rate \\ \hline
    $N_1$ & $K_1$ & $N_1^{-1}$ & $1/(N_1K_1)$ & $1/K_1$& $\varepsilon^2/N_1$\\
    $1$ & $K_2$ & $N_1^{-1}$& $1/K_2$ & $N_1^2/ K_2$& $\varepsilon^2$
\end{tabular}
\end{center}
We fix the privacy budgets of both groups to $N_1^{-1}$ and the batch sample size of the second group to 1.
Suppose that $N_1> 1$, $K_2 > N_1^2 K_1 $ and $1/K_1< \varepsilon^2< N_1/K_1$.  According to \Cref{thm-hetero-rate}, the first group minimises $R(\varepsilon,l_1,l_2)$ with the rate $1/K_1$, even if the size of the second group $K_2\rightarrow \infty$. Therefore, the availability of the second group of batches does not help decrease the fundamental difficulty of the problem, regardless of the number of batches it contains. Thus, we conclude that using only the first group of batches remains minimax-optimal (up to poly-logarithmic factors), even though it represents an arbitrarily small proportion of the sample as $K_2 \rightarrow \infty$. 

To illustrate the interactivity between the costs of privacy and contamination, first consider the problem in the absence of contamination ($\varepsilon=0$). Here, the optimal rate of $N_1^2/K_2$ is attained using the second group of batches. Instead, in the absence of privacy constraints, the optimal rate of $\varepsilon^2/N_1$ is attained using the first group of batches. However, in the presence of both contamination and privacy constraints, the optimal rate of $1/K_1$ found above is strictly larger than $N_1^2/K_2 \vee \varepsilon^2/N_1$, the maximum of the costs of privacy and contamination alone. 

For the rest of this section, we present in detail the lower bound (\Cref{hetero-LB}) in \Cref{sec-hetero-LB}. In \Cref{sec-hetero-UB}, we establish the upper bound (\Cref{non-adaptive-thm}) using \Cref{alg:hetero-estimation}, which confirms that the most informative group of servers alone is sufficient for near-minimax optimal estimation. Since the identity of this group depends on the contamination level $\varepsilon$, in \Cref{sec-adaptation}, we also develop a method that adapts to $\varepsilon$ (\Cref{alg:adaptive-mean-estimation}) while achieving the optimal rate up to poly-logarithmic factors (see \Cref{hetero-UB-thm}).  

\subsection{Lower bound}\label{sec-hetero-LB}
Using the batch partition framework detailed in \eqref{server-group-def} above, we present the following lower bound on the minimax risk.
\begin{proposition}\label{hetero-LB}
    Under the contamination model with heterogeneous server sizes detailed in~\eqref{heterogeneous-model}, we have that the minimax risk defined in \eqref{general-minimax-risk} satisfies 
    \begin{equation*}
        \mathcal{R}
        \gtrsim 1\wedge \Lambda^{-1}\min_{\substack{l_1,l_2\in \mathbb{N}: \\\mathcal{I}_{l_1,l_2} \neq \emptyset}} \left\{R(\varepsilon,l_1,l_2)\right\},
    \end{equation*}
    where $\Lambda= \{1+\log_2(n_{\max}/n_{\min})\}\{1+\log_2(\alpha_{\max}/\alpha_{\min})\}$ and $\mathcal{I}_{l_1,l_2}$ is defined in \eqref{server-group-def}. 
\end{proposition}

The proof of \Cref{hetero-LB} is deferred to \Cref{appendix-hetero-LB-proof}. Before discussing the proof, we note that the quantity $\Lambda$ corresponds to the maximum total number of groups formed according to the batch partition framework. It is comprised of the product of the two logarithmic factors depending only on the heterogeneity of the batch sample sizes $n_{\max}/n_{\min}$ and that of the privacy budgets $\alpha_{\max}/\alpha_{\min}$. In the homogeneous special case where $n_{\max}=n_{\min}$ and $\alpha_{\max} = \alpha_{\min}$, we have $\Lambda = 1$ and the logarithmic factor vanishes. 

The proof uses Le Cam's two-point lemma \citep[e.g.~Lemma 1 in][]{Yu1997assouad}, and decomposes the divergence between two distributions of the privatised output $Z_{1:K}$ into contributions from each group of batches. Since the homogeneous rate associated with each group is at least that of $\min_{\substack{l_1,l_2\in \mathbb{N}: \mathcal{I}_{l_1,l_2} \neq \emptyset}} R(\varepsilon,l_1,l_2)$, each of their contributions to the KL divergence can be upper bounded according to one of the three cases -- whether the group's dominant cost is that of privacy, contamination or the classical sampling variance. Completing the bound under each of the three cases corresponds to one of the proof strategies discussed under \Cref{LB-Gaussian}.

\subsection{Upper bound}\label{sec-hetero-UB}

Motivated by \Cref{hetero-LB}, derived under the batch partition framework, we now construct an estimator that focuses on making use of one particular group of batches. Given that the chosen group of batches to be used is labelled by $(\tilde{l}_1,\tilde{l}_2)$, \Cref{alg:hetero-estimation} presented below details the construction of the estimator. For simplicity, \Cref{alg:hetero-estimation} takes the entire data set as input, but only the $(\tilde{l}_1,\tilde{l}_2)$ group is used to produce the estimator. Further, we consider only the case where $|\mathcal{I}_{\tilde{l}_1,\tilde{l}_2}|$ is even for notational simplicity. When $|\mathcal{I}_{\tilde{l}_1,\tilde{l}_2}|$ is odd, we discard one batch in the group $\mathcal{I}_{\tilde{l}_1,\tilde{l}_2}$ randomly before proceeding with \Cref{alg:hetero-estimation}. This has no effect on the stated rates in \Cref{non-adaptive-thm} and \Cref{hetero-UB-thm} below.

\begin{algorithm}
    \caption{\textsc{HeterogeneousEstimation}($\{\{X_{k,i}\}_{i=1}^{n_k}\}_{k=1}^K$,$\{\alpha_k\}_{k=1}^K$, $(\tilde{l}_1,\tilde{l}_2)$, $\delta$)}
    \label{alg:hetero-estimation}
    \begin{algorithmic}[1]
        \INPUT {Data $\{\{X_{k,i}\}_{i=1}^{n_k}\}_{k=1}^K$, privacy parameters $\{\alpha_k\}_{k=1}^K$, batch group label $(\tilde{l}_1,\tilde{l}_2)$, tuning parameter $\delta$.}
    \State Randomly split the set $\mathcal{I}_{\tilde{l}_1,\tilde{l}_2}$ defined in \eqref{server-group-def} into two sets, $\mathcal{B}$ and $\mathcal{B}'$, of equal size.
    \State Set $T \leftarrow T_{\tilde{l}_1,\tilde{l}_2}$ as defined in \eqref{eq:hetero-UB-thm-T-def}, $r \leftarrow r_{\tilde{l}_1,\tilde{l}_2}$ as defined in \eqref{eq:hetero-UB-thm-r-def}, $M \leftarrow 3r_{\tilde{l}_1,\tilde{l}_2}$.
    \State Set $\mathfrak{A}^*_{\tilde{l}_1,\tilde{l}_2} \leftarrow$ {\sc MeanLocalisation}($\{\{X_{k,i}\}_{i=1}^{n_k}\}_{k\in \mathcal{B}}$, $\{\alpha_k\}_{k\in \mathcal{B}}$, $r$, $T$) \Comment{\Cref{alg:mean-localisation}}
    
    \State Set $\hat{\theta} \leftarrow$ {\sc MeanRefinement}($\{\{X_{k,i}\}_{i=1}^{n_k}\}_{k\in \mathcal{B}'}$, $\mathfrak{A}^*_{\tilde{l}_1,\tilde{l}_2}$, $\{\alpha_k\}_{k\in \mathcal{B}'}$, $M$, $T$)\Comment{\Cref{alg:mean-refinement}}

        \State \Return $\hat{\theta}$
    \end{algorithmic}
\end{algorithm}

Extended to handle small amounts of heterogeneity, \Cref{alg:hetero-estimation} applies \Cref{alg:homogeneous-Gaussian-mean-estimation} to the group of batches labelled by $(\tilde{l}_1,\tilde{l}_2)$. An additional input parameter $\delta$ is included to allow for tuning the success probability of the procedure, which is required for the adaptive procedure, \Cref{alg:adaptive-mean-estimation}, later. The tuning parameters $T$, $r$ and $M$ analogous to the homogeneous setting are now chosen according to the inputs $(\tilde{l}_1,\tilde{l}_2)$ and $\delta$. \Cref{non-adaptive-thm} below provides the theoretical guarantees of \Cref{alg:hetero-estimation}.

\begin{proposition}\label{non-adaptive-thm}
    Under the conditions of \Cref{thm-hetero-rate},
    the output estimator $\hat{\theta}$ of \Cref{alg:hetero-estimation} satisfies
    \begin{equation}\label{non-adaptive-thm-high-prob-guarantee}
        \mathbb{P}\left(|\hat{\theta}-\theta|> L(\varepsilon,l_1^*,l_2^*,\delta)\right)\leq 2\delta
    \end{equation}
    with the batch group label $(\tilde{l}_1,\tilde{l}_2)=(l_1^*,l_2^*)$ and a tuning parameter $\delta \geq e^{-C^2K_{l_1^*,l_2^*}}$ as input, where the function $L(\cdot,\cdot,\cdot,\cdot)$ is defined in \eqref{hetero-L-def}. In particular, taking $\delta \asymp \min_{\substack{l_1,l_2\in \mathbb{N}: \mathcal{I}_{l_1,l_2} \neq \emptyset}} \left\{R(\varepsilon,l_1,l_2)\right\}$, we have that
    \begin{equation}\label{non-adaptive-risk-bound}\E \{|[\hat{\theta}]_{-1}^1-\theta|^2\} \lesssim \log^2\left(\sum_{k=1}^K n_k\right)\min_{\substack{l_1,l_2\in \mathbb{N}: \\\mathcal{I}_{l_1,l_2} \neq \emptyset}} \left\{R(\varepsilon,l_1,l_2)\right\}.\end{equation}
\end{proposition}

The proof of \Cref{non-adaptive-thm} is deferred to \Cref{appendix-hetero-UB-proof}. Guarantees on the localisation step are proved by adapting the proof techniques used in \Cref{UB-thm}. Note that, in addition to the expectation bound \eqref{non-adaptive-risk-bound}, we provide a high probability guarantee \eqref{non-adaptive-thm-high-prob-guarantee} for $\hat{\theta}$ to be used later in the proof of the theoretical guarantees of our adaptive algorithm (\Cref{hetero-UB-thm}). To prove \eqref{non-adaptive-thm-high-prob-guarantee}, we derive concentration bounds for $\hat{\theta}$ in addition to bounding its bias.

Choosing an optimal group $(l_1^*,l_2^*)$ as the input of \Cref{alg:hetero-estimation} is instrumental to attaining the upper bound \eqref{non-adaptive-risk-bound} on the minimax risk, which matches the lower bound in \Cref{hetero-LB} up to poly-logarithmic factors. When the level of contamination $\varepsilon$ is large, groups with large server sample sizes are favoured.  When $\varepsilon$ is small, the selected group is instead determined primarily by the FDP effective sample size. Note that although the computation of $(l_1^*,l_2^*)$ as defined in \eqref{eq:optimal-group-def} does require the knowledge of the contamination level $\varepsilon$, it does not require any private quantities. Therefore, the privacy guarantees of \Cref{alg:hetero-estimation} follow from those of \Cref{alg:homogeneous-Gaussian-mean-estimation}.

In the high probability guarantee, the quantity $L(\varepsilon,l_1,l_2,\delta)$ is defined by
\begin{equation}\label{hetero-L-def}
    L(\varepsilon, l_1,l_2, \delta) = C'\left\{\varepsilon M_{l_1,l_2}\vee\sqrt{\frac{\log(4/\delta)}{|\mathcal{I}_{l_1,l_2}|}}\left(M_{l_1,l_2}\vee\frac{T_{l_1,l_2}}{\tilde{n}_{l_1}\tilde{\alpha}_{l_2}}\right)\right\},
\end{equation}where $C'>0$ is a sufficiently large absolute constant introduced for later convenience, and the tuning parameters $\{T_{l_1,l_2}, r_{l_1,l_2},M_{l_1,l_2}\}_{l_1,l_2 \in \mathbb{N}: \mathcal{I}_{l_1,l_2}\neq \emptyset}$ are generalisations of those used in the homogeneous setting. Specifically, they are defined by
    \begin{equation}\label{eq:hetero-UB-thm-T-def}
        T_{l_1,l_2}= C_1\left( 1 \vee \sqrt{\log\{K_{l_1,l_2}(\tilde{n}_{l_1}\wedge\tilde{n}_{l_1}^2 \tilde{\alpha}_{l_2}^2)\}}\right),
    \end{equation}    \begin{equation}\label{eq:hetero-UB-thm-r-def}
        r_{l_1,l_2} = C_2 T_{l_1,l_2} \left(\frac{1}{\tilde{n}_{l_1}}\vee \frac{T_{l_1,l_2}^2+\log(4/\delta)}{\tilde{n}_{l_1}^2 \tilde{\alpha}_{l_2}^2|\mathcal{I}_{l_1,l_2}|} \right)^{1/2} , \quad M_{l_1,l_2}=3r_{l_1,l_2},
    \end{equation}
    where $C_1,C_2>0$ are sufficiently large absolute constants. To interpret $L(\varepsilon,l_1,l_2,\delta)$, note that with the tuning parameters defined in \eqref{eq:hetero-UB-thm-T-def} and \eqref{eq:hetero-UB-thm-r-def}, we have that 
\begin{equation}\label{approx-L-def}
L(\varepsilon,l_1,l_2,\delta)=\tilde{L}_{\varepsilon,l_1,l_2,\delta}\sqrt{R(\varepsilon,l_1,l_2)},
\end{equation}
where \begin{equation}\label{eq:L-tilde-range}\tilde{L}\equiv \tilde{L}_{\varepsilon,l_1,l_2,\delta}\in \left[1,C\left\{1\vee\log\left(\sum_{k=1}^K n_k\right)\vee \log(4/\delta)\right\}\right]\end{equation} and $C>1$ is an absolute constant. Thus, up to logarithmic factors, $L(\varepsilon,l_1^*,l_2^*,\delta)$ is of the same order as  $\min_{\substack{l_1,l_2\in \mathbb{N}: \mathcal{I}_{l_1,l_2} \neq \emptyset}} \sqrt{R(\varepsilon,l_1,l_2)}$, the minimax rate as stated in \Cref{thm-hetero-rate}.

This method, however, relies on a known contamination level $\varepsilon$. As $\varepsilon$ is often unknown in practice, we propose an adaptive procedure in \Cref{sec-adaptation} below.

\subsection{Adaptation to contamination level}\label{sec-adaptation}

We proceed to present \Cref{alg:adaptive-mean-estimation}, a Lepski-type method \citep[e.g.][]{Lepskii1991on, Lepskii1992asymptotically}, where $\varepsilon_t$ acts as the tuning parameter to be calibrated, chosen over a dyadic grid
\[
    \varepsilon_t = 2^{-t}, \quad t = 0, 1, \ldots, \lceil\log_2(K)/2\rceil.
\]
For each candidate value $\varepsilon_t$, \Cref{alg:adaptive-mean-estimation} selects the group indexed by $(l_1^t, l_2^t)$ that minimises $L(\varepsilon_t, l_1, l_2, \delta)$ over all non-empty dyadic groups.  It then applies \Cref{alg:hetero-estimation} to this group.  This produces an estimator $\hat{\theta}_{l_1^t, l_2^t}$ and a corresponding confidence interval $I_t$ centred at $\hat{\theta}_{l_1^t, l_2^t}$, with radius $L(\varepsilon_t, l_1^t, l_2^t, \delta)$. Further, noting that $L(\varepsilon,l_1,l_2,\delta)$ is constant in $\varepsilon$ for $\varepsilon < K^{-1/2}$ (since $|\mathcal{I}_{l_1,l_2}|\leq K$ and $\sqrt{\log(4/\delta)}\geq 1$, see definition \eqref{hetero-L-def}), it is only necessary to iterate over $\lceil\log_2(K)/2\rceil$ iterations.

Since the function $L(\varepsilon,l_1,l_2,\delta)$ is increasing in $\varepsilon$ and $\varepsilon_t$ decreases with $t$, the intervals $I_t$ become narrower as $t$ increases, reflecting the fact that smaller contamination levels allow for more accurate estimation. It is shown in the proof that, with high probability, all intervals corresponding to candidate values $\varepsilon_t \geq \varepsilon$ contain the true mean $\theta$. \Cref{alg:adaptive-mean-estimation}, therefore, selects $I_{\hat{t}}$, the smallest interval before the intersection with all larger intervals become empty.    It chooses the most accurate estimator whose confidence interval remains compatible with those obtained under assumptions of larger contamination levels.

\begin{algorithm}
\caption{Adaptive Gaussian mean estimation with heterogeneous servers}\label{alg:adaptive-mean-estimation}
\begin{algorithmic}[1]
    \INPUT {Data $\{\{X_{k,i}\}_{i=1}^{n_k}\}_{k=1}^K$, privacy parameters $\{\alpha_k\}_{k=1}^K$}, tuning parameter $\delta$.
    \State Set $S \leftarrow \emptyset$, $I_0 \leftarrow [-1,1]$, $t \leftarrow 0$.
    \While{$\bigcap_{t'= 0}^t I_{t'}\neq \emptyset$ and $t< \lceil\log_2(K)/2\rceil$}
    \State Set $t\leftarrow t+1$, $\varepsilon_t \leftarrow 2^{-t}$
    \State Set
    \begin{equation}\label{adaptive-alg-minimiser-condition}(l_1^t, l_2^t) \leftarrow \argmin_{l_1,l_2\in \mathbb{N}: \mathcal{I}_{l_1,l_2} \neq \emptyset} L(\varepsilon_t,l_1,l_2,\delta),\end{equation}
    where the function $L(\cdot,\cdot,\cdot,\cdot)$ is defined in \eqref{hetero-L-def}.
    \If{$(l_1^t, l_2^t)\notin S$}
    \State Set $S \leftarrow S\cup\{(l_1^t,l_2^t)\}$
    \State Set $\hat{\theta}_{l_1^t,l_2^t} \leftarrow$\textsc{HeterogeneousEstimation}($\{\{X_{k,i}\}_{i=1}^{n_k}\}_{k=1}^K$,$\{\alpha_k\}_{k=1}^K$, $(l_1^t,l_2^t)$, $\delta$)
    \State Set 
    \begin{equation}\label{eq:interval-def}I_t \leftarrow \left[\hat{\theta}_{l_1^t,l_2^t} - L(\varepsilon_t, l_1^t,l_2^t, \delta) , \hat{\theta}_{l_1^t,l_2^t} + L(\varepsilon_t, l_1^t,l_2^t, \delta) \right].\end{equation}
    \Else
    \State Set
    \[I_t \leftarrow \left[\hat{\theta}_{l_1^t,l_2^t} - L(\varepsilon_t, l_1^t,l_2^t,\delta) , \hat{\theta}_{l_1^t,l_2^t} + L(\varepsilon_t, l_1^t,l_2^t, \delta)\right].\]
    \EndIf
    \EndWhile
    \If{$\bigcap_{t'= 0}^t I_{t'}= \emptyset$}
    \State Set $\hat{t} \leftarrow t-1$
    \Else 
    \State Set $\hat{t} \leftarrow t$
    \EndIf
    \OUTPUT{$I_{\hat{t}}$}
\end{algorithmic}
\end{algorithm}

A notable feature of \Cref{alg:adaptive-mean-estimation} is that, in contrast to typical Lepski-type methods, the selection of $\varepsilon_t$ also determines the subset of data used in constructing the estimator. Occasionally, the same dyadic group is selected for multiple values of $\varepsilon_t$. In this case, the estimator for that group is computed only once and then reused.  Each server hence contributes to the procedure at most once. Further, the selection of the dyadic group depends only on $\delta$, $\varepsilon_t$ and public quantities $n_k$ and $\alpha_k$. Thus, the FDP guarantee is preserved without privacy budget composition across iterations, while reusing an already computed group estimator also reduces computation.

The theoretical guarantees of \Cref{alg:adaptive-mean-estimation} are provided in \Cref{hetero-UB-thm} below. 
\begin{theorem}\label{hetero-UB-thm}
    Consider the contamination model with heterogeneous server sizes detailed in~\eqref{heterogeneous-model}. Define the estimator $\hat{\theta}=[\hat{m}]_{-1}^1$, where $\hat{m}$ is the midpoint of the output interval $I_{\hat{t}}$ of \Cref{alg:adaptive-mean-estimation}. Recall the definition of $(l_1^t,l_2^t)$ in \eqref{adaptive-alg-minimiser-condition}. There exists an absolute constant $c>0$ such that when the condition
    \begin{equation}\label{adaptive-UB-condition}
        \max_{t\in [\lceil\log_2(K)/2\rceil]}\left\{ \sqrt{\frac{\log(\tilde{n}_{l_1^t}\wedge \tilde{n}_{l_1^t}^2\tilde{\alpha}_{l_2^t}^2K_{l_1^t,l_2^t})}{K_{l_1^t,l_2^t}}}+ \frac{1}{\sqrt{(\tilde{n}_{l_1^t}\wedge \tilde{n}_{l_1^t}^2\tilde{\alpha}_{l_2^t}^2)K_{l_1^t,l_2^t}}}\right\}+\varepsilon<c\end{equation}
    holds with \begin{equation}\label{adaptive-UB-condition-delta}
        \delta \geq \max_{t\in[\lceil\log_2(K)/2\rceil] } \exp\{-c^2K_{l_1^t,l_2^t}\},
    \end{equation} we have that
    \begin{equation}\label{hetero-UB-thm-highprob-bound}
        \mathbb{P}\left(|\hat{\theta}-\theta|> \min_{l_1,l_2\in \mathbb{N}:\mathcal{I}_{l_1,l_2}\neq \emptyset}6L(\varepsilon,l_1,l_2,\delta)\right)\leq 2\delta\log_2(K).
    \end{equation}

    In particular, taking 
    \begin{multline}\label{eq:delta-choice-hetero-UB-proof}\delta \asymp \log^2\left(\sum_{k=1}^Kn_k\right)\{\log_2(K)\}^{-1}\min_{l_1,l_2\in \mathbb{N}:\mathcal{I}_{l_1,l_2}\neq \emptyset}\left\{\frac{1}{\tilde{n}_{l_1}^2\tilde{\alpha}_{l_2}^2 K_{l_1,l_2}}\vee \frac{1}{\tilde{n}_{l_1}K_{l_1,l_2}}\right\}\\
    \vee\left\{\max_{t\in[\lceil\log_2(K)/2\rceil] } \exp\{-C^2 K_{l_1^t,l_2^t}\}\right\},\end{multline}
    we have that \begin{equation*}\E[|\hat{\theta}-\theta|^2]\lesssim \log^2\left(\sum_{k=1}^K n_k\right)\left[\min_{\substack{l_1,l_2\in \mathbb{N}:\\\mathcal{I}_{l_1,l_2}\neq \emptyset}}\left\{R(\varepsilon,l_1,l_2)\right\} \vee\max_{t\in[\lceil\log_2(K)/2\rceil] } \exp\{-C^2K_{l_1^t,l_2^t}\}\right].\end{equation*}
\end{theorem}
The proof of \Cref{hetero-UB-thm} is deferred to \Cref{appendix-hetero-UB-proof}. The probability bound~\eqref{hetero-UB-thm-highprob-bound} shows that, up to poly-logarithmic factors, the adaptive estimator performs as well as the oracle estimator that knows $\varepsilon$ and therefore knows which group of servers to use (see \Cref{non-adaptive-thm}).  Choosing $\delta$ appropriately controls the contribution from the complement of this high-probability event, leading to the stated expectation bound.  As a consequence, \Cref{hetero-UB-thm} matches the lower bound in \Cref{hetero-LB} up to poly-logarithmic factors, provided that the exponentially decaying terms are dominated by the rate.

The condition \eqref{adaptive-UB-condition} is the heterogeneous analogue of the condition \eqref{UB-condition-thm1} imposed in the homogeneous case, where each dyadic group $\mathcal{I}_{l_1,l_2}$ is viewed as a homogeneous collection of servers.  The condition \eqref{adaptive-UB-condition} therefore requires that the homogeneous upper bound condition hold uniformly over the groups used by the procedure.  Groups that contain too few servers, or whose effective sample size is too small, cannot reliably support the localisation and/or refinement steps. In applications where there is greater heterogeneity of the batch sizes and privacy budgets, one may instead use a coarser geometric partition to ensure that each group contains a sufficient number of servers. Specifically, consider instead
\[
    \mathcal{I}_{l_1,l_2}=\{k\in[K]:n_k/n_{\min}\in[b^{l_1-1},b^{l_1}),\alpha_k/\alpha_{\min}\in[b^{l_2-1},b^{l_2})\},
\]
where $b>2$ is an absolute constant. With the method's input tuning parameters adapted, the rate changes by at most constant factors.

\section{Conclusion}
This paper establishes a foundation for studying robustness and privacy jointly in federated learning by combining batch contamination with FDP.  Through a minimax analysis of mean estimation, we quantify the effects of heavy-tailedness across the full range of moments $v > 1$ and identify how heterogeneity in batch sizes and privacy budgets change the interaction between privacy and contamination.  Our results distinguish the separation of these costs in homogeneous settings from their interactions under heterogeneity.  The batch partitioning framework provides a way to characterise this interaction.  The most informative group of batches determines the minimax rate up to poly-logarithmic factors, and, in some regimes, an arbitrarily small proportion of the data suffices for near-optimal estimation.

These findings provide a basis for developing robust and private methods for other statistical tasks, including multivariate mean estimation, covariance matrix estimation, regression and hypothesis testing. 
Such extensions offer opportunities to establish broader connections between 
robustness and privacy, combining those identified under 
CDP and LDP \citep[e.g.][]{Kamath20private,Liu2021robust,acharya2021robust,Li2023on}. 
In particular, investigating which features of the heavy-tailedness phase transition and the disentanglement of privacy and contamination costs persist beyond univariate mean estimation could reveal common principles across statistical problems.  Under heterogeneity, the batch partitioning framework also motivates a fundamental question about when aggregation across groups is essential.

Another promising direction is to incorporate distributional heterogeneity alongside varying batch sizes and privacy budgets.
Our partitioning framework offers a starting point for connecting these sources of heterogeneity with existing approaches to mean shifts \citep[e.g.][]{aliakbarpour2026high} and transfer and multi-task learning \citep[e.g.][]{tian2026contaminated}.  Developing these connections could yield procedures that account jointly for differences in source distributions, sample sizes and privacy requirements, while remaining robust to contaminated batches.  These extensions would also advance the non-private theory of batch contamination and broaden the scope of rigorous statistical guarantees for federated learning.

\section*{Acknowledgements}
Berrett was supported by European Research Council Starting Grant 101163546.
Yu was partially supported by the EPSRC programme grant DASS (EP/Z531327/1), UKRI-JST grant Statistical Safeguarding and Leverhulme Trust (Philip Leverhulme Prize).

\bibliographystyle{apalike}
\bibliography{ref}
\appendix
\section*{Appendices}
\section{Proofs of results in \Cref{sec-mean-est}}\label{appendix-gaussian-mean-est}
\subsection{Proof of lower bounds}\label{sec:homogeneous-LB-proofs}
\begin{proof}[Proof of Propositions \ref{LB-Gaussian} and \ref{finite-moments-LB}] We prove the lower bounds for both the Gaussian setting (\Cref{LB-Gaussian}) and the heavy-tailed setting (\Cref{finite-moments-LB}) using the following proof strategy. In \textbf{Step~0}, we decompose the minimax risk into the maximum of three terms: the non-private risk without contamination $\mathcal{R}_\infty(0)$, the non-private risk in the presence of batch contamination $\mathcal{R}_\infty(\varepsilon)$, and the private risk without contamination $\mathcal{R}_\alpha(0)$. These are proved separately in \textbf{Steps 1}, \textbf{2} and \textbf{3}, respectively, and we combine the results in \textbf{Step 4}. We begin by introducing notation used throughout the proof.

\medskip
\noindent\textbf{Preliminary notation. }For each $k \in [K]$, let $M_{0,k}$, $M_{1,k}$, $G_{0,k}$ and $G_{1,k}$ denote distributions on $\mathbb{R}^{n}$, and let $P_0$ and $P_1$ be distributions on $\mathbb{R}$, such that 
\begin{equation}\label{homogeneous-LB-proof-distributions}M_{0,k} = (1-\varepsilon)P_0^{\otimes n}+\varepsilon G_{0,k} \quad \mathrm{and} \quad M_{1,k} = (1-\varepsilon)P_1^{\otimes n}+\varepsilon G_{1,k}.\end{equation}
We let $\theta_0$ and $\theta_1$ denote the means of $P_0$ and $P_1$ respectively, such that
\[\E_{X\sim P_0}[X] = \theta_0, \quad \text{and} \quad \E_{X\sim P_1}[X]=\theta_1.\]
Let $Q_k M_{0,k}$ and $Q_k M_{1,k}$ denote the distributions of the privatised data $Z_k$ conditioned on $Z_{1:k-1}$, where $X_{1:K}$ are generated from $M_{0,1}\times \dots \times M_{0,K}$ and $M_{1,1}\times \dots \times M_{1,K}$ respectively. Recall that the privacy mechanisms $\{Q_k\}_{k=1}^K$ are such that the privatised data $Z_k$ satisfy a common privacy constraint of $\alpha$--FDP for each $k\in [K]$. Specifically, for any measurable set $S$, we have that
\[
    Q_kM_{0,k}(S |Z_{1:k-1}=z_{1:k-1}) =Q_kM_{0,k}(S |z_{1:k-1}) =\int Q_k(S|x_k, z_{1:k-1})\,\dd M_{0,k}(x_k),
\]
noting that $X_k$ is independent of $Z_{1:k-1}$. For each $k\in [K]$, let $Q(M_{0,1}\times \dots\times M_{0,k})$ and $Q(M_{1,1}\times \dots\times M_{1,k})$ denote the distribution of $Z_{1:k}$ when $X_{1:k}$ is generated from $M_{0,1}\times \dots \times M_{0,k}$ and $M_{1,1}\times \dots \times M_{1,k}$ respectively. 

\medskip
\noindent\textbf{Step 0: Decomposition. }Recall the minimax risk $\mathcal{R}_\alpha(\varepsilon)$ defined in \eqref{Gaussian-minimax-risk}. Observing the nested nature of our problem, for any $\alpha\geq 0$ and $\varepsilon \in [0,1]$ we have \[\mathcal{Q}_\alpha \subseteq \mathcal{Q}_{\infty} \quad \text{and} \quad \mathcal{P}_0(\mathcal{P})\subseteq\mathcal{P}_\varepsilon(\mathcal{P}),\] which implies that \[\mathcal{R}_\alpha(\varepsilon)\geq \mathcal{R}_\infty(\varepsilon) \quad \text{and} \quad \mathcal{R}_\alpha(\varepsilon)\geq \mathcal{R}_\alpha(0)\] respectively. Thus, we can lower bound the minimax risk using the decomposition
    \begin{equation}\label{eq:LB-decomposition}
    \mathcal{R}_\alpha(\varepsilon) \geq \mathcal{R}_\infty(0)\vee \mathcal{R}_{\infty}(\varepsilon)\vee \mathcal{R}_\alpha(0) .
\end{equation}We proceed by lower bounding the minimax risk under these three rates separately.

\medskip
\noindent\textbf{Step 1: The non-contaminated, non-private rate $\mathcal{R}_\infty(0)$. }To obtain a lower bound on~$\mathcal{R}_\infty(0)$, we set $G_{0,k} = P_0^{\otimes n}$ and $G_{1,k}= P_1^{\otimes n}$, such that we reduce to the non-contaminated regime ($\varepsilon=0$) with $M_{0,k}=P_0^{\otimes n}$ and $M_{1,k} = P_1^{\otimes n}$. Using Le Cam's two-point lemma \citep[e.g.~Lemma 1 in][]{Yu1997assouad}, we have
\begin{align}
    \mathcal{R}_\infty(0) & \geq \frac{|\theta_0-\theta_1|^2}{8}\left\{1- \TV\left(Q\prod_{k=1}^K M_{0,k} , Q\prod_{k=1}^K M_{1,k}\right)\right\} \nonumber\\
    & \geq \frac{|\theta_0-\theta_1|^2}{8}\left\{1- \TV\left(P_0^{\otimes nK},P_1^{\otimes nK}\right)\right\} \label{parametric-lecam}
\end{align}
for any privacy mechanism $Q$, where the last inequality follows from the data processing inequality \citep[e.g.~Theorem 7.4 in][]{polyanskiy2025information}.

\textbf{Case 1: The Gaussian setting.} Choose $P_0 = \mathcal{N}(\theta_0,1)$ and $P_1 =  \mathcal{N}(\theta_1,1)$, and let $\theta_0=-1$ and $\theta_1 = 4/(3\sqrt{nK})-1$. It follows from Pinsker's inequality \citep[e.g.~Lemma 2.5 in][]{Tsybakov2009introduction} that
\[\TV\left(P_0^{\otimes nK},P_1^{\otimes nK}\right)\leq \sqrt{\KL(P_0^{\otimes nK},P_1^{\otimes nK})/2} = |\theta_0-\theta_1|\sqrt{nK}/2=2/3.\]
Hence, with \eqref{parametric-lecam}, we have
\begin{equation}\label{Gaussian-parametric-LB}\mathcal{R}_\infty(0) \geq \frac{|\theta_0-\theta_1|^2}{8}\left\{1-\TV\left(P_0^{\otimes nK},P_1^{\otimes nK}\right)\right\}\geq \frac{2}{27nK}.\end{equation}

\textbf{Case 2: The heavy-tailed setting.} For the heavy-tailed setting, we show that, for $nK\geq 1$,
\begin{equation}
    \label{eq:heavy-tailed-parametric-LB} \mathcal{R}_\infty(0)\geq A_v'\left(\frac{1}{nK}\vee \frac{1}{(nK)^{2-2/v}}\right),
\end{equation}where $A_v'>0$ is a constant depending only on $v$. To prove this, we consider two different constructions for the pair of distributions $(P_0,P_1)$; we extend the Gaussian construction to prove the rate of $1/(nK)$ observed in the Gaussian setting, and we construct two-point distributions to prove the rate $1/(nK)^{2-2/v}$ specific to the heavy-tailed setting. 

\textbf{Construction 1. } Extending the Gaussian construction, note that the $p$th central moment of the $\mathcal{N}(\theta,\sigma^2)$ distribution satisfies
\[\frac{1}{\sigma^p}\E[(X-\theta)^p] = \frac{ p!}{2^{p/2}(p/2)!}=: C_p\geq 1\]
if $p$ is even and positive. Hence, we consider a rescaling of the normal distribution with  \[P_0=\mathcal{N}(\theta_0,C_p^{-2/p}) \quad \text{and} \quad P_1=\mathcal{N}(\theta_1,C_p^{-2/p}),\]
and let $\theta_0=-1$, $\theta_1=4C_p^{-1/p}/(3\sqrt{nK})-1$ and $p$ to be the smallest even number that is at least~$v$.
Thus, by Jensen's inequality, we have that
\[\E_{X\sim P_j}[|X-\theta_j|^v]\leq \E_{X\sim P_j}[|X-\theta_j|^p]^{v/p}= 1\]
for $j \in \{0, 1\}$, verifying that $P_0,P_1\in \mathcal{P}_v$. Continuing from \eqref{parametric-lecam}, we have that
\begin{align}
    \mathcal{R}_\infty(0)&\geq \frac{|\theta_0-\theta_1|^2}{8}\left\{1- \TV\left(P_0^{\otimes nK},P_1^{\otimes nK}\right)\right\}\nonumber\\
    &\stackrel{(a)}{\geq} \frac{|\theta_0-\theta_1|^2}{8}\left(1- \sqrt{\KL\left(P_0^{\otimes nK},P_1^{\otimes nK}\right)/2}\right)\nonumber\\
    &=\frac{|\theta_0-\theta_1|^2}{8}\left(1-C_p^{1/p}|\theta_0-\theta_1|\sqrt{nK}/2\right)\nonumber\\
    &\stackrel{(b)}{=}\frac{2C_p^{-2/p}}{27nK}\nonumber\\
    &\stackrel{(c)}{=} \frac{A_v}{nK}\nonumber,
\end{align} where $A_v>0$ is a constant depending only on $v$. Note that (a) follows from Pinsker's inequality \citep[e.g. Lemma 2.5 in][]{Tsybakov2009introduction}, (b) follows from our choices of $\theta_0$ and $\theta_1$, and (c) follows from our choice of $p$.

\textbf{Construction 2. }
Consider a construction using two-point distributions similar to the proof of Theorem 3.1 in \cite{Devroye2016sub-gaussian}. Let $p=1/(nK)$ and $c=1/(2p)^{1/v}$, and define the distributions~$P_0$ and $P_1$ as
\[P_0(\{0\}) = P_1(\{0\}) = 1-p, \quad P_0(\{c\}) = P_1(\{-c\}) = p.\]
Since $\theta_0= cp=-\theta_1$, the distributions~$P_0$ and~$P_1$ are contained in the class of distributions~$\mathcal{P}_v$ defined in~\eqref{finite-moments-class} with \[\E_{X\sim P_0}[|X-\theta_0|^v] = \E_{X\sim P_1}[|X-\theta_1|^v] = c^v p(1-p) (p^{v-1} + (1-p)^{v-1})\leq 2pc^v= 1.\]
Observe that
\begin{align*}\TV(P_0^{\otimes nK},P_1^{\otimes nK})&=\sup_{A}\{P_0^{\otimes nK}(A) - P_1^{\otimes nK}(A)\}\\
&=P_0^{\otimes nK}(\{0,c\}^{nK}) - P_1^{\otimes nK}(\{0,c\}^{nK})\\
&=1-(1-p)^{nK}\\
&= 1-\left(1-\frac{1}{nK}\right)^{nK},\end{align*}
where the supremum is taken over all measurable sets $A$. Plugging this into \eqref{parametric-lecam}, we have
\begin{align}\label{eq:heavy-tailed-parametric-LB-2}
    \mathcal{R}_\infty(0)&\geq \frac{|\theta_0-\theta_1|^2}{8}\left(1-\frac{1}{nK}\right)^{nK}\stackrel{(a)}{=}\frac{1}{2^{1+2/v}(nK)^{2-2/v }}\left(1-\frac{1}{nK}\right)^{nK}\stackrel{(b)}{\geq} \frac{2^{-(3+2/v)}}{(nK)^{2-2/v}}
\end{align}
provided that $nK \geq 2$, where (a) follows from our choice of $|\theta_0-\theta_1|= 2cp=(2/nK)^{1-1/v}$, and (b) follows from the observation that the function $f(x) = (1-1/x)^x$ is increasing for $x>1$. 

\medskip
\noindent \textbf{Step 2: The contaminated, non-private rate $\mathcal{R}_{\infty}(\varepsilon)$.}  First note that, by choosing $\theta(P_0) = 0$ and $\theta(P_1)=1$ in \Cref{contamination-LB-using-Bernoulli},
    \begin{equation}\label{universal-batch-contamination-LB}\mathcal{R}_\infty(\varepsilon)\geq \frac{\varepsilon^{\lceil K/2 \rceil}}{4}\end{equation}
    holds under both settings. 
    
    Moreover, provided that 
    \begin{equation}\label{eq:contamination-LB-TV-bound}
        \tau := \TV(P_0^{\otimes n}, P_1^{\otimes n})\leq \varepsilon/(1-\varepsilon),
    \end{equation} 
    let $p_0$ and $p_1$ denote the densities of $P_0^{\otimes n}$ and $P_1^{\otimes n}$ respectively and define $G_{0,k}$ and $G_{1,k}$ by their densities \[g_0 = (1-c)p_0+ c(p_1-p_0)_+/\tau\quad \text{and} \quad g_1 = (1-c)p_1+c(p_0-p_1)_+/\tau \]respectively, where $c = \tau/\{(1+\tau)\varepsilon\}$. It then follows that $M_{0,k}=M_{1,k}$ \citep[see e.g.~the Proof of Theorem 5.1 in][]{Chen2016a}. Thus, Le Cam's two-point lemma \citep[e.g.~Lemma 1 in][]{Yu1997assouad} leads to that
    \begin{equation}\label{contamination-lecam}\mathcal{R}_\infty(\varepsilon) \geq \frac{|\theta_0-\theta_1|^2}{8}\left\{1- \TV\left(Q\prod_{k=1}^K M_{0,k} , Q\prod_{k=1}^K M_{1,k}\right)\right\} = \frac{|\theta_0-\theta_1|^2}{8},\end{equation}
    provided that \eqref{eq:contamination-LB-TV-bound} holds. In the following, we construct inlier distributions $P_0$ and $P_1$ satisfying \eqref{eq:contamination-LB-TV-bound} in both the Gaussian and heavy-tailed settings.
    
    \textbf{Case 1: The Gaussian setting.} Choose $P_0 = \mathcal{N}(\theta_0,1)$ and $P_1 =  \mathcal{N}(\theta_1,1)$, and let $\theta_0 = -1$ and $\theta_1 = 2\varepsilon/\sqrt{n}-1$. It follows from Pinsker's inequality \citep[e.g.~Lemma 2.5 in][]{Tsybakov2009introduction} that \eqref{eq:contamination-LB-TV-bound} is satisfied with
\[
    \mathrm{TV}(P_0^{\otimes n},P_1^{\otimes n}) \leq \sqrt{\mathrm{KL}(P_0^{\otimes n},P_1^{\otimes n})/2} = |\theta_0 - \theta_1|\sqrt{n}/2 = \varepsilon.
\]
Using \eqref{contamination-lecam}, we have that \begin{equation}\label{Gaussian-contamination-LB}\mathcal{R}_\infty(\varepsilon)\geq\frac{|\theta_0-\theta_1|^2}{8} = \frac{\varepsilon^2}{2n}.\end{equation}  

\textbf{Case 2: The heavy-tailed setting.} For the heavy-tailed setting, we show that
\begin{equation}\label{eq:heavy-tailed-contamination-LB}\mathcal{R}_\infty(\varepsilon)\geq B_v' \left(\frac{\varepsilon^2}{n}\vee \frac{\varepsilon^{2-2/v}}{n^{2-2/v}}\right),\end{equation}
where $B_v'>0$ is a constant depending only on $v$.
Similarly to \textbf{Step 1}, we prove the two rates separately by considering two different constructions for ($P_0,P_1$). 

\textbf{Construction 1.} We use the same Gaussian construction as in \textbf{Step 1} with $\theta_0=0$ and $\theta_1=2C_p^{-1/p}\varepsilon/\sqrt{n}$ instead, and define the distributions $P_0$ and $P_1$ as  \[P_0=\mathcal{N}(\theta_0,C_p^{-2/p}) \quad \text{and} \quad P_1=\mathcal{N}(\theta_1,C_p^{-2/p}),\]
where $p$ is the smallest even number that is at least $v$.
It follows from Pinsker's inequality \citep[e.g. Lemma 2.5 in][]{Tsybakov2009introduction} that \eqref{eq:contamination-LB-TV-bound} is satisfied with
\[\mathrm{TV}(P_0^{\otimes n}, P_1^{\otimes n}) \leq \sqrt{\mathrm{KL}(P_0^{\otimes n}, P_1^{\otimes n})/2} = C_p^{1/p}|\theta_0 - \theta_1|\sqrt{n}/2=\varepsilon.\]
Using \eqref{contamination-lecam}, we have that
\begin{equation}\label{eq:finite-moments-contamination-LB1}\mathcal{R}_\infty(\varepsilon)\geq\frac{ |\theta_0-\theta_1|^2}{8} = C_p^{-2/p}\frac{\varepsilon^2}{2n}= B_v\frac{\varepsilon^2}{n},\end{equation}
where $B_v$ is a positive constant depending only on $v$.

\textbf{Construction 2.} Consider a construction using two-point distributions similar to the proof of Proposition 4 in \cite{Li2023on}. Define the distributions $P_0$ and $P_1$ by 
    \[P_0(\{0\}) = 1,\]
    \[\begin{cases}
        P_1(\{0\}) = 1-\varepsilon n^{-1}\\
        P_1(\{(2\varepsilon n^{-1})^{-1/v}\})  = \varepsilon n^{-1}.
    \end{cases}\]
    Since $\theta_0=0$ and $\theta_1=2^{-1/v}(\varepsilon n^{-1})^{1-1/v}$, observe that $P_0$ and $P_1$ are contained in the class of distributions $\mathcal{P}_v$ defined in \eqref{finite-moments-class} with $\mathbb{E}_{X\sim P_0}[|X-\theta_0|^v] = 0$ and
    \[\mathbb{E}_{X\sim P_1}[|X-\theta_1|^v] = \frac{1}{2}\left((1-\varepsilon n^{-1})(\varepsilon n^{-1})^{v-1}+(1-\varepsilon n^{-1})^v\right)\leq 1.\]
    Further, observe that \eqref{eq:contamination-LB-TV-bound} is satisfied with
    \begin{align*}
        \TV(P_0^{\otimes n},P_1^{\otimes n}) 
        =\sup_{A}\{P_0^{\otimes n}(A) - P_1^{\otimes n}(A)\}= P_0^{\otimes n}(\{0\}^n) - P_1^{\otimes n}(\{0\}^n)=1-\left(1-\frac{\varepsilon}{n}\right)^n \leq \varepsilon,
    \end{align*}
    where the supremum is taken over all measurable sets $A$.
Using \eqref{contamination-lecam}, we have that
\begin{equation}\label{eq:finite-moments-contamination-LB2}\mathcal{R}_\infty(\varepsilon)\geq \frac{|\theta_0-\theta_1|^2}{8}= 2^{-(3+2/v)}\frac{\varepsilon^{2-2/v}}{n^{2-2/v}}.\end{equation}

\medskip
\noindent \textbf{Step 3: The private, non-contaminated rate $\mathcal{R}_\alpha(0)$. } We set $G_{0,k} = P_0^{\otimes n}$ and $G_{1,k}= P_1^{\otimes n}$, reducing to the non-contaminated regime ($\varepsilon=0$). By considering the chain rule for KL divergence, we have
\begin{align}  
    &\KL\left( Q\prod_{k=1}^K M_{0,k} ,Q \prod_{k=1}^K M_{1,k}\right)\nonumber\\
    &= \sum_{k=1}^K \int \KL(Q_kM_{0,k}(\cdot \mid z_{1:{k-1}}), Q_kM_{1,k}(\cdot \mid z_{1:{k-1}}))\dd Q(M_{0,1}\times \dots \times M_{0,k-1})(z_{1:k-1}).\label{eq:KL-chain-rule}
\end{align}
Further, let $q_0$ and $q_1$ be the densities, with respect to a measure $\mu$, of $Q_kM_{0,k}(\cdot \mid z_{1:{k-1}})$ and $Q_kM_{1,k}(\cdot \mid z_{1:{k-1}})$ respectively. Using Lemma 6.1 in \cite{karwa2017finite}, we have that
\begin{equation}\label{eq:density-ratio-bound1}
    e^{-\alpha'}\leq q_1/q_0\leq e^{\alpha'},
\end{equation}
where $\alpha' = 6\alpha n \TV(P_0,P_1)$. 
Then, for any $z_{1:k-1}$, we have that
\begin{align}\KL(Q_kM_{0,k}(\cdot \mid z_{1:{k-1}}), Q_kM_{1,k}(\cdot \mid z_{1:{k-1}}))&\leq \int \log\left(\frac{q_0}{q_1}\right)q_0 \dd \mu+ \int \log\left(\frac{q_1}{q_0}\right)q_1 \dd \mu\nonumber\\
&= \int \log\left(\frac{q_1}{q_0}\right)\left(\frac{q_1}{q_0}-1\right)q_0\dd \mu\nonumber\\
&\stackrel{(a)}{\leq} \int \{\alpha' (e^{\alpha'}-1)\vee (-\alpha')(e^{-\alpha'}-1)\}q_0\dd \mu\nonumber \\
&\leq \int \alpha' (e^{\alpha '}-1)q_0 \dd \mu\nonumber\\
&= \alpha' (e^{\alpha '}-1),\label{eq:bound-from-KV-lemma}
\end{align}
where (a) follows from \eqref{eq:density-ratio-bound1}. Hence, by Le Cam's two-point lemma, we have that 
\begin{align}\label{private-lecam}\mathcal{R}_\alpha(0) \geq \frac{|\theta_0-\theta_1|^2}{8}\left(1-\sqrt{\frac{1}{2}\KL\left( Q\prod_{k=1}^K M_{0,k} ,Q \prod_{k=1}^K M_{1,k}\right)}\right)\geq \frac{|\theta_0-\theta_1|^2}{8}\left(1-\sqrt{\frac{K\alpha'(e^{\alpha'}-1)}{2}}\right),\end{align}
where the final inequality follows from \eqref{eq:KL-chain-rule} and \eqref{eq:bound-from-KV-lemma}. Using \eqref{private-lecam}, we prove the result for each setting separately by considering specific distributions $P_0$ and $P_1$.

\textbf{Case 1: The Gaussian setting.}
Let $P_0 = \mathcal{N}(\theta_0,1)$ and $P_1 = \mathcal{N}(\theta_1,1)$, and let $\theta_0=0$ and $\theta_1= 1 \wedge 2/(9\alpha n\sqrt{K}) $. Since
\[\TV(P_0,P_1)\leq |\theta_0-\theta_1|/2 = 1/2\wedge 1/(9\alpha n\sqrt{K}),\] we have that \[\alpha' = 6\alpha n \TV(P_0,P_1)\leq \frac{2}{3\sqrt{K}}\leq 1.\]
Therefore, applying the inequality $e^{\alpha'}-1\leq 2\alpha'$ that holds for $\alpha'\leq 1$ to \eqref{private-lecam}, we have that
\begin{equation}\label{Gaussian-private-LB}\mathcal{R}_\alpha(0) \geq \frac{|\theta_0-\theta_1|^2}{8}\left(1-\sqrt{\frac{K\alpha'(e^{\alpha'}-1)}{2}}\right)\geq \left(\frac{1}{8}\wedge\frac{1}{162n^2K\alpha^2}\right)\left(1-\sqrt{K\alpha'^2}\right) \geq \frac{1}{24}\wedge\frac{1}{486n^2K\alpha^2}.\end{equation}

\textbf{Case 2: The heavy-tailed setting. }Consider the distributions $P_0$ and $P_1$ defined as $P_0(\{0\}) = 1$ and
    \[
        \begin{cases}
            P_1(\{0\}) = 1-\{1\wedge C(n\alpha \sqrt{K})^{-1}\}\\
            P_1(\{(2\{1\wedge C(n\alpha \sqrt{K})^{-1}\})^{-1/v}\}) = \{1\wedge C(n\alpha \sqrt{K})^{-1}\},
        \end{cases}
    \]
    where $C<1/6$ is an absolute constant. Since $\theta_0=0$ and $\theta_1 = 2^{-1/v}\{1\wedge C(n\alpha\sqrt{K})^{-1}\}^{1-1/v}$, the distributions~$P_0$ and~$P_1$ are contained in the class of distributions $\mathcal{P}_v$ defined in~\eqref{finite-moments-class} with~$\mathbb{E}_{X\sim P_0}[|X-\theta_0|^v] = 0$ and
    \[\E_{X\sim P_1}[|X-\theta_1|^v]=\frac{1}{2}((1- \{1\wedge C(n\alpha \sqrt{K})^{-1}\})\{1\wedge C(n\alpha \sqrt{K})^{-1}\}^{v-1}+(1-\{1\wedge C(n\alpha \sqrt{K})^{-1}\})^v)\leq 1.\]
    Since $\TV(P_0,P_1) = 1\wedge C(n\alpha \sqrt{K})^{-1}$, we have that \[\alpha' = 6\alpha n \TV(P_0,P_1)\leq 6CK^{-1/2}\leq 1\] for any $C\leq 1/6$. Therefore, applying the inequality $e^{\alpha'}-1\leq 2\alpha'$ that holds for $\alpha'\leq 1$ to \eqref{private-lecam}, we have that
    \begin{align}\mathcal{R}_\alpha(0)&\geq \frac{|\theta_0-\theta_1|^2}{8}\left(1-\sqrt{\frac{K\alpha'(e^{\alpha'}-1)}{2}}\right)\nonumber\\
    &\geq \left(1\wedge \frac{C^{2-2/v}}{(n^2K\alpha^2)^{1-1/v}}\right)\frac{1-\sqrt{K\alpha'^2}}{8\cdot 2^{2/v}}\nonumber\\
    &\geq \frac{(1-6C)}{8\cdot2^{2/v}}\left(1\wedge \frac{C^{2-2/v}}{(n^2K\alpha^2)^{1-1/v}}\right).\label{eq:heavy-tailed-private-LB}\end{align}

    \medskip
    \noindent\textbf{Step 4: Conclusion. }
    
    \textbf{Case 1: The Gaussian setting. }
    Plugging in \eqref{Gaussian-parametric-LB}, \eqref{universal-batch-contamination-LB}, \eqref{Gaussian-contamination-LB} and \eqref{Gaussian-private-LB} into \eqref{eq:LB-decomposition}, we have that
    \[\mathcal{R}_\alpha(\varepsilon)\geq C\left(\frac{1}{nK}\vee\frac{1}{1\vee n^2 K \alpha^2}\vee \frac{\varepsilon^2}{n}\vee \varepsilon^{\lceil K/2 \rceil}\right)\]
    for some absolute constant $C>0$.
    
    \textbf{Case 2: The heavy-tailed setting. }
    Plugging in \eqref{eq:heavy-tailed-parametric-LB}, \eqref{universal-batch-contamination-LB}, \eqref{eq:heavy-tailed-contamination-LB}  and \eqref{eq:heavy-tailed-private-LB} into \eqref{eq:LB-decomposition}, we have that
    \[\mathcal{R}_\alpha(\varepsilon)\geq C_v\left(\frac{1}{nK}\vee \frac{1}{(nK)^{2-2/v}}\vee\frac{1}{1\vee (n^2K\alpha^2)^{1-1/v}}\vee\frac{\varepsilon^{2-2/v}}{n^{2-2/v}}\vee\frac{\varepsilon^2}{n}\vee \varepsilon^{\lceil K/2 \rceil}\right)\]
    for some constant $C_v>0$ depending only on $v$.
    \end{proof}
\begin{proposition}\label{contamination-LB-using-Bernoulli}
    Let $\mathcal{P}$ denote a class of distributions on $\mathbb{R}$. Under the batch contamination model defined in \eqref{general-data-distribution} and \eqref{general-contamination-model}, the minimax risk satisfies
    \[\inf_{\hat{\theta}}\sup_{\{P_{\varepsilon,k}\}_{k = 1}^K\in\mathcal{P}_\varepsilon(\mathcal{P})}\mathbb{E}_{\{P_{\varepsilon,k}\}_{k=1}^K,Q}[(\hat{\theta}-\theta(P))^2] \geq \frac{|\theta(P_0)-\theta(P_1)|^2}{4}\varepsilon^{\lceil K/2 \rceil}\]
    for any fixed privacy mechanism $Q$ and any distributions $P_0, P_1 \in \mathcal{P}$.
\end{proposition}
\begin{proof}
For brevity, we denote $\theta_0=\theta(P_0)$ and $\theta_1 = \theta(P_1)$. Observe that generating 
\[X_k\sim M_{j,k}:=(1-\varepsilon)P_j^{\otimes n_k}+\varepsilon G_{j,k},\] for $j \in \{0,1\}$, is equivalent to the following:
\begin{enumerate}
    \item Independently generate $E_k\sim \mathrm{Bernoulli}(\varepsilon)$, $Y_k \sim P_j^{\otimes n_k}$, and $W_k\sim G_{j,k}$.
    \item Set $X_k = (1-E_k)Y_k+E_kW_k$.
\end{enumerate}
Here, the random variable $E_k=1$ if the batch $k$ is contaminated; otherwise $E_k=0$. Let $E = \sum_{k=1}^KE_k$ be the total number of contaminated batches. Choose \[G_{0,k} = P_1^{\otimes n_k}\quad \text{and} \quad G_{1,k} = P_0^{\otimes n_k}\] for all $k\in[K]$. Then, we can see that \begin{align}&\E_{\{M_{0,k}\}_{k=1}^K, Q}[|\hat{\theta}-\theta_0|^2|E = e]\nonumber\\
=& \sum_{\substack{(e_1,\dots,e_K)\in \{0,1\}^K,\\ \sum_k{e_k} = e}}\E_{\{M_{0,k}\}_{k=1}^K, Q}[|\hat{\theta}-\theta_0|^2|(E_1,\dots,E_K) = (e_1,\dots,e_K)]\left[{K\choose e}\right]^{-1}\nonumber\\
\stackrel{(a)}{=}& \sum_{\substack{(e_1,\dots,e_K)\in \{0,1\}^K,\\ \sum_k{e_k} = e}}\E_{\{M_{1,k}\}_{k=1}^K, Q}[|\hat{\theta}-\theta_0|^2|(E_1,\dots,E_K) = (1-e_1,\dots,1-e_K)]\left[{K\choose K-e}\right]^{-1}\nonumber\\
\label{contamination-LB-symmetry-argument}=& \E_{\{M_{1,k}\}_{k=1}^K,Q}[|\hat{\theta}-\theta_0|^2|E = K-e]\end{align}
for any fixed privacy mechanism $Q$ and estimator $\hat{\theta}$, where (a) follows from the fact that we have constructed $M_{0,k}$ and $M_{1,k}$ such that the inlier distribution under $M_{0,k}$ is equal to the outlier distribution under $M_{1,k}$, and the inlier distribution under $M_{1,k}$ is equal to the outlier distribution under $M_{0,k}$. Using \eqref{contamination-LB-symmetry-argument}, we have that
\begin{align*}
    &\sup_{\{P_{\varepsilon,k}\}_{k = 1}^K\in\mathcal{P}_\varepsilon(\mathcal{P})}\mathbb{E}_{\{P_{\varepsilon,k}\}_{k=1}^K,Q}[|\hat{\theta}-\theta(P)|^2]\\
    &\geq \frac{1}{2}\left(\E_{\{M_{0,k}\}_{k=1}^K, Q}[|\hat{\theta}-\theta_0|^2]+\E_{\{M_{1,k}\}_{k=1}^K, Q}[|\hat{\theta}-\theta_1|^2]\right)\\
    &= \frac{1}{2}\left(\sum_{e=0}^K \E_{\{M_{0,k}\}_{k=1}^K, Q}[|\hat{\theta}-\theta_0|^2|E = e]\mathbb{P}(E=e) +\sum_{e=0}^K\E_{\{M_{1,k}\}_{k=1}^K, Q}[|\hat{\theta}-\theta_1|^2|E = e]\mathbb{P}(E=e)\right)\\
    &\stackrel{(a)}{=} \frac{1}{2}\left(\sum_{e=0}^K \E_{\{M_{1,k}\}_{k=1}^K,Q}[|\hat{\theta}-\theta_0|^2|E = K-e]\mathbb{P}(E=e) +\sum_{e=0}^K\E_{\{M_{1,k}\}_{k=1}^K, Q}[|\hat{\theta}-\theta_1|^2|E = e]\mathbb{P}(E=e)\right)\\
    &= \frac{1}{2}\left(\sum_{e=0}^K \E_{\{M_{1,k}\}_{k=1}^K,Q}[|\hat{\theta}-\theta_0|^2|E = e]\mathbb{P}(E=K-e) +\sum_{e=0}^K\E_{\{M_{1,k}\}_{k=1}^K, Q}[|\hat{\theta}-\theta_1|^2|E = e]\mathbb{P}(E=e)\right)\\
    &\geq \frac{1}{2}\left(\sum_{e=0}^K \E_{\{M_{1,k}\}_{k=1}^K,Q}[|\hat{\theta}-\theta_0|^2+|\hat{\theta}-\theta_1|^2|E = e](\mathbb{P}(E=K-e)\wedge \mathbb{P}(E=e)) \right)\\
    &\geq \frac{1}{4}\left(\sum_{e=0}^K |\theta_0-\theta_1|^2(\mathbb{P}(E=K-e)\wedge \mathbb{P}(E=e)) \right)\\
    &\geq \frac{|\theta_0-\theta_1|^2}{4}\mathbb{P}(E = \lceil K/2 \rceil)\\
    &\stackrel{(b)}{=} \frac{|\theta_0-\theta_1|^2}{4} \binom{K}{\lceil K/2 \rceil} \varepsilon^{\lceil K/2 \rceil}(1-\varepsilon)^{K-\lceil K/2 \rceil}\\
    &\geq \frac{|\theta_0-\theta_1|^2}{4} 2^{K-\lceil K/2 \rceil} \varepsilon^{\lceil K/2 \rceil}(1-\varepsilon)^{K-\lceil K/2 \rceil}\\
    &\geq \frac{|\theta_0-\theta_1|^2}{4} \varepsilon^{\lceil K/2 \rceil},
\end{align*}
    where (a) follows from \eqref{contamination-LB-symmetry-argument} and (b) follows from the observation that $E\sim\mathrm{Binomial}(K,\varepsilon)$. 
\end{proof}

\subsection{Proof of upper bounds} \label{appendix:homogeneous-UB-proofs}
Before proving the upper bounds in Propositions~\ref{UB-thm}, \ref{finite-moment-UB} and~\ref{infinite-variance-UB}, and the required lemmas, we define the following notation for the specifics of \Cref{alg:homogeneous-Gaussian-mean-estimation} that we will use throughout the proofs in the rest of this subsection. For notational simplicity, we instead let the total number of servers be $2K$ rather than $K$ as specified in the main text.  Note that this change will only result in a constant change in the rates.

Regarding \Cref{alg:mean-localisation}, for each $k\in [K]$ and $j \in [\lceil2r^{-1}\rceil]$, let the location of the $j$th bin be defined by \begin{equation}\label{a_j-def}a_j = -1 +(j-1/2)r,\end{equation} and denote the smoothed votes by
    \begin{equation}\label{Z_kj-def}Z_{k,j} = \left(2r - \left|\frac{1}{n}\sum_{i=1}^n [X_{k,i}]_{-T}^T - a_j\right|\right)_+ + \frac{16 T}{\alpha n} W_{k,j},\end{equation}
    where $\{W_{k,j}\}_{k\in[K],j\in [\lceil2r^{-1}\rceil]}$ are i.i.d.~standard Laplace random variables.  We further let
    \[
        j^* = \min\argmax_{j\in [\lceil 2r^{-1}\rceil ]}\frac{1}{K}\sum_{k=1}^K Z_{k,j}
    \]
    and 
    \begin{equation}\label{J-def}
        \mathfrak{A}^*= a_{j^*}= -1 +(j^* -1/2)r. 
    \end{equation}
    
Regarding \Cref{alg:mean-refinement}, for each $k \in \{K+1,\dots, 2K\}$, let
    \begin{equation} \label{eq-def-Z_k}
        Z_k = \mathfrak{A}^*+ \left[\frac{1}{n}\sum_{i = 1}^n [X_{k,i}]_{-T}^T -\mathfrak{A}^*\right]_{-M}^M + \frac{2T}{\alpha n}L_k ,
    \end{equation}
    where $\{L_k\}_{k = K+1}^{2K}$ are i.i.d.~standard Laplace random variables, therefore
    the estimator $\hat{\theta}$ from the output of \Cref{alg:homogeneous-Gaussian-mean-estimation} can be written as
    \[\hat{\theta} = \frac{1}{K}\sum_{k=K+1}^{2K} Z_k.\]
Further, for $l \in [K]$, define 
\begin{equation}\label{R_k-def}R_l = \frac{1}{n}\sum_{i = 1}^n [X_{K+l,i}]_{-T}^T -\mathfrak{A}^*\end{equation}
and
\begin{equation}\label{S_k-def}S_l = (\mathfrak{A}^*-\theta) + \left[R_l\right]_{-M}^M,\end{equation}
such that we have
\begin{equation}\label{hat-theta-minus-theta}\hat{\theta}-\theta = \frac{1}{K}\sum_{l=1}^{K}S_l+\frac{2T}{nK\alpha}\sum_{k=K+1}^{2K} L_{k} .\end{equation}

\begin{proof}[Proof of Propositions \ref{UB-thm}, \ref{finite-moment-UB} and \ref{infinite-variance-UB}]

In this proof, we refer to the settings considered in Propositions~\ref{UB-thm}, \ref{finite-moment-UB}, and~\ref{infinite-variance-UB} by the Gaussian, finite variance and infinite variance settings, respectively. The following proof strategy is used to prove all three propositions. In \textbf{Step 1}, we decompose the mean squared error into bias and variance terms. Since the behaviour of these terms depends on the key quantity $S_l$, we derive a rough upper bound \eqref{S1-bound} on $S_l$ in \textbf{Step 2} which will be used to control the contributions from contaminated batches. In \textbf{Step 3}, we obtain \eqref{E-indicator-bound}, a control of the initial estimator $\mathfrak{A}^*$. With~\eqref{S1-bound} and~\eqref{E-indicator-bound} in hand, we proceed to \textbf{Step 4} and \textbf{5}, where we upper bound the bias and variance terms respectively. In \textbf{Step 6}, we complete the proof by combining the upper bounds on the bias and variance terms.

\medskip
\textbf{Step 1: The bias-variance decomposition. }Recalling in \eqref{J-def} the initial estimator $\mathfrak{A}^*$ output by \Cref{alg:mean-localisation}, we have that
\begin{align}
    & \mathbb{E}\{(\hat{\theta} - \theta)^2\} = \mathbb{E}[\mathbb{E}\{(\hat{\theta} - \theta)^2 | \mathfrak{A}^*\}] = \mathbb{E} \left[ \left\{\mathbb{E}(\hat{\theta}| \mathfrak{A}^*) - \theta\right\}^2 \right] + \mathbb{E} \left\{\mathrm{Var}(\hat{\theta} | \mathfrak{A}^*)\right\}\nonumber\\
    &\stackrel{(a)}{=}  \E \left[\left|\frac{1}{K}\sum_{l=1}^K \E [S_l|\mathfrak{A}^*]\right|^2\right]+\E\left( \mathrm{Var}\left(\left.\frac{1}{K}\sum_{l=1}^KS_l+\frac{2T}{nK\alpha}\sum_{k=K+1}^{2K} L_{k}\right|\mathfrak{A}^*\right)\right)\nonumber\\
    &\label{bias-variance-decomposition}\stackrel{(b)}{=}\E \left[\left|\frac{1}{K}\sum_{l=1}^K \E [S_l|\mathfrak{A}^*]\right|^2\right]+\E\left[\frac{1}{K^2}\sum_{l=1}^K \mathrm{Var}\left(S_l|\mathfrak{A}^*\right)\right]+ \frac{8T^2}{n^2K\alpha^2}\\
    &= (I) + (II) + \frac{8T^2}{n^2K\alpha^2}\nonumber,
\end{align}
where (a) follows from \eqref{hat-theta-minus-theta} and the independence of $\{L_k\}_{k=K+1}^{2K}$ from $\mathfrak{A}^*$, and (b) follows from independence of the Laplace mechanism from the $S_l$'s and that the $S_l$'s are conditionally independent given $\mathfrak{A}^*$. 
In the rest of the proof, we upper bound the bias term (I) and the variance term (II) separately. To do this, we will need the following facts about the quantities $S_l$'s and $\mathfrak{A}^*$. 

\medskip
\noindent\textbf{Step 2: An upper bound on $S_l$. }Noting that our choices for the tuning parameters in \eqref{UB-thm-tuning-params},~\eqref{infinite-var-m-choice} and~\eqref{finite-var-ub-m-choice} in each of the three settings respectively all imply that $r\lesssim M$, we have 
    \begin{align} \label{S1-bound}
        |S_l|& = \left|\left[R_l\right]_{-M}^M + (\mathfrak{A}^*-\theta)\right| \lesssim M+r\mathbf{1}\{|\mathfrak{A}^*-\theta|\leq 5r/2\}+\mathbf{{1}}\{|\mathfrak{A}^*-\theta|>5r/2\}\nonumber\\
        &\lesssim M+\mathbf{{1}}\{|\mathfrak{A}^*-\theta|>5r/2\},
    \end{align}
    where the first inequality follows from the assumption that $\theta\in [-1,1]$ in our models \eqref{Gaussian-distribution} and \eqref{finite-moments-class}, and the fact that our initial estimator $\mathfrak{A}^*$ defined in \eqref{J-def} satisfies $\mathfrak{A}^*\in [-1,1+r/2]$.  We will use \eqref{S1-bound} to bound $S_l$ when the corresponding $(K+l)$th batch is contaminated. 
    
    \medskip
    \noindent\textbf{Step 3: Upper bounding $|\mathfrak{A}^*-\theta|$ with high probability. }For each of the Gaussian, infinite variance and finite variance settings that we consider, Lemmas \ref{high-prob-lemma}, \ref{high-prob-lemma-infinite-variance} and \ref{high-prob-lemma-finite-moments} show that
    \begin{equation}\label{E-indicator-bound}
        \mathbb{P}(|\mathfrak{A}^*-\theta|>5r/2)\leq \delta
    \end{equation}
    holds for each setting respectively, provided that the condition \eqref{eq:high-prob-lemma-common-assumption} is satisfied by each of the corresponding choices of $T$, $r$ and $\delta$. These lemmas are proved separately in \Cref{sec-UB-lemma-proofs}. We now show that, with our specific choices of $T$, $r$ and $
    \delta$, the condition \eqref{UB-condition-thm1} is sufficient to ensure that \eqref{eq:high-prob-lemma-common-assumption} holds in each setting. 

    \textbf{Case 1: The Gaussian setting. } For the rest of the proof, we choose  
    \begin{equation}\label{delta-def-gaussian}\delta = \frac{12}{C_2\{K(n\wedge n^2\alpha^2)\}}\end{equation} for the Gaussian setting. Plugging \eqref{delta-def-gaussian}
    and the choices of $T$ and $r$ in \eqref{UB-thm-tuning-params} into the LHS of  \eqref{eq:high-prob-lemma-common-assumption}, we have that
    \begin{align*}
        &\varepsilon+\sqrt{\frac{\log(12/(r\delta))}{K}}+\frac{T\sqrt{\log(12/(r\delta))}}{rn\alpha\sqrt{K}}+\frac{1}{\sqrt{K}(\sqrt{n}\wedge n\alpha)}\\
        &= \varepsilon+ \sqrt{\log\left(\frac{K\{n\wedge n^2\alpha^2\}}{T\{n^{-1/2}\vee T(n^2K\alpha^2)^{-1/2}\}}\right)}\left\{\frac{1}{\sqrt{K}}+\frac{T}{rn\alpha \sqrt{K}}\right\}+\frac{1}{\sqrt{K}(\sqrt{n}\wedge n\alpha)}\\
        &\leq \varepsilon+\sqrt{\frac{3\log(K\{n\wedge n^2\alpha^2\})}{2}}\left\{\frac{1}{\sqrt{K}}+\frac{1}{C_2T}\right\}+\frac{1}{\sqrt{K}(\sqrt{n}\wedge n\alpha)}\\
        &\leq \varepsilon+\sqrt{\frac{3\log(K\{n\wedge n^2\alpha^2\})}{2K}} +\frac{1}{\sqrt{K}(\sqrt{n}\wedge n\alpha)}+\frac{\sqrt{3}}{C_1 C_2 \sqrt{2}}\\
        &<C'
    \end{align*}
    for a fixed constant $C'>0$, provided that the constants $C_1$ and $C_2$ in the choices of $T$ and $r$ in \eqref{UB-thm-tuning-params} are sufficiently large, and the constant $C$ in condition \eqref{UB-condition-thm1} is sufficiently small.

    \textbf{Case 2: The infinite variance setting. }Plugging the choices of $T$, $r$ and $\delta$ in \eqref{infinite-var-ub-t-choice}, \eqref{infinite-var-ub-r-choice} and~\eqref{infinite-var-m-choice} into the LHS of~\eqref{eq:high-prob-lemma-common-assumption}, we have that 
    \begin{align*}
        &\varepsilon+\sqrt{\frac{\log(12/(r\delta))}{K}}+\frac{T\sqrt{\log(12/(r\delta))}}{rn\alpha\sqrt{K}}+\frac{1}{\sqrt{K}(\sqrt{n}\wedge n\alpha)}\\
        &\leq \varepsilon+\sqrt{\log\left(\frac{(nK)^{2-2/v}\wedge (n^2K\alpha^2)^{1-1/v}}{T^{1-v}}\right)}\left\{\frac{1}{\sqrt{K}}+\frac{T^v}{B_1 n\alpha\sqrt{K}\sqrt{\log(nK\wedge n\alpha\sqrt{K})}}\right\}\\
        &\hspace{8cm}+\frac{1}{\sqrt{K}(\sqrt{n}\wedge n\alpha)}\\
        &\leq \varepsilon+\sqrt{\left(3-\frac{3}{v}\right)\log(nK\wedge n\alpha\sqrt{K})}\left\{\frac{1}{\sqrt{K}}+\frac{1}{B_1 \sqrt{\log(nK\wedge n\alpha\sqrt{K})}}\right\}\\
        &\hspace{8cm}+\frac{1}{\sqrt{K}(\sqrt{n}\wedge n\alpha)}\\
        &\leq \varepsilon+\sqrt{\left(3-\frac{3}{v}\right)\frac{\log(nK\wedge n\alpha\sqrt{K})}{K}}+\frac{1}{B_1}\sqrt{\left(3-\frac{3}{v}\right)}+\frac{1}{\sqrt{K}(\sqrt{n}\wedge n\alpha)}\\
        &\leq C'
    \end{align*}
    for a fixed constant $C'>0$, provided that the constant $B_1$ in the choice of $r$ is sufficiently large and the constant $C$ in condition \eqref{UB-condition-thm1} is sufficiently small.

    \textbf{Case 3: The finite variance setting. }For the rest of the proof, we choose 
    \begin{equation}\label{delta-choice-finite-variance}
        \delta = \frac{12}{A_2}\left(\frac{1}{nK}\vee \frac{1}{(n^2K\alpha^2)^{1-1/v}}\right)
    \end{equation}
    for the finite variance setting. Plugging \eqref{delta-choice-finite-variance} and the choices of $T$ and $r$ given in \eqref{finite-var-ub-t-choice} into the LHS of \eqref{eq:high-prob-lemma-common-assumption}, we have that
    \begin{align*}
    &\varepsilon+\sqrt{\frac{\log(12/(r\delta))}{K}}+\frac{T\sqrt{\log(12/(r\delta))}}{rn\alpha\sqrt{K}}+\frac{1}{\sqrt{K}(\sqrt{n}\wedge n\alpha)}\\
    &\leq \varepsilon+\sqrt{\log\left(\frac{nK\wedge (n^2K\alpha^2)^{1-1/v}}{n^{-1/2}\vee T^{1-v}}\right)}\left\{\frac{1}{\sqrt{K}}+ \frac{T^v}{A_2 n\alpha\sqrt{K}\sqrt{\log(K(n\wedge n^2\alpha^2))}}\right\}+\frac{1}{\sqrt{K}(\sqrt{n}\wedge n\alpha)}\\
    &\leq \varepsilon+\sqrt{\frac{3\log(K\{n\wedge n^2\alpha^2\})}{2}}\left\{\frac{1}{\sqrt{K}}+ \frac{A_1^v}{A_2 \sqrt{\log(K(n\wedge n^2\alpha^2))}}\right\}+\frac{1}{\sqrt{K}(\sqrt{n}\wedge n\alpha)}\\
    &\leq \varepsilon +\sqrt{\frac{3\log(K\{n\wedge n^2\alpha^2\})}{2K}}+\frac{A_1^v\sqrt{3}}{A_2\sqrt{2}}+\frac{1}{\sqrt{K}(\sqrt{n}\wedge n\alpha)}\\
    &\leq C'
    \end{align*}
    for a fixed constant $C'>0$, provided that the constant $A_2$ in the choice of $r$ is sufficiently large and the constant $C$ in condition \eqref{UB-condition-thm1} is sufficiently small. Note that the second inequality holds provided that $A_1\geq 1$.

\medskip
\noindent \textbf{Step 4: The bias term (I).} Having verified that \eqref{E-indicator-bound} holds in each setting, we focus on the high probability event where our initial estimator $\mathfrak{A}^*$ is accurate in the following steps. Note that
\begin{align}
    \E\left[\left|\frac{1}{K}\sum_{l=1}^K \E [S_l|\mathfrak{A}^*]\right|^2\right]
&\leq \E\left[\max_{l\in [K]}|\E [S_l|\mathfrak{A}^*]|^2\right]\nonumber\\
&\lesssim \E\left[ \left\{\varepsilon (M+\mathbf{{1}}\{|\mathfrak{A}^*-\theta|>5r/2\})+(1-\varepsilon)\max_{l\in[K]}|\E_{X_{K+l}\sim P^{\otimes n}}[S_l|\mathfrak{A}^*]|\right\}^2\right]\nonumber\\
&\label{bias-decomposition}\lesssim \varepsilon^2M^2+\varepsilon^2\delta+\E\left[\max_{l\in[K]}|\E_{X_{K+l}\sim P^{\otimes n}}[S_l|\mathfrak{A}^*]|^2\right],\end{align}
    where the second inequality uses \eqref{S1-bound} and the fact that $X_{K+l}$ follows a mixture distribution, and the final inequality uses \eqref{E-indicator-bound}. Upper bounds on $|\E_{X_{K+l}\sim P^{\otimes n}}[S_l|\mathfrak{A}^*]|$ and choices of tuning parameters differ in the three problem settings. We proceed by considering each setting separately.
    
    \textbf{Case 1: The Gaussian setting.} Recall that $M=3r$ with $r$ defined in \eqref{UB-thm-tuning-params}. Applying \eqref{UB-thm-tuning-params} and \Cref{E-S-given-J-lemma} to \eqref{bias-decomposition}, 
    \begin{align}
        \E \left[\left|\frac{1}{K}\sum_{l=1}^K \E [S_l|\mathfrak{A}^*]\right|^2\right]&\lesssim T^2\left(\frac{\varepsilon^2}{n}\vee \frac{\varepsilon^2T^2}{n^2K\alpha^2}\right)+\varepsilon^2\delta+\frac{1}{nK}\vee \frac{1}{n^2K\alpha^2}+\mathbb{P}(|\mathfrak{A}^*-\theta|>5r/2)\nonumber\\
        \label{bias-bound-Gaussian}&\lesssim T^2\left(\frac{\varepsilon^2}{n}\vee \frac{\varepsilon^2T^2}{n^2K\alpha^2}\right)\vee \frac{1}{nK}\vee \frac{1}{n^2K\alpha^2},
    \end{align}
    where the final inequality follows from \eqref{E-indicator-bound} and the definition of $\delta$ in \eqref{delta-def-gaussian}.

    \textbf{Case 2: The infinite variance setting.} Applying \Cref{E-S-given-J-lemma-infinite-variance} to \eqref{bias-decomposition}, we have that
    \begin{align}
        \E \left[\left|\frac{1}{K}\sum_{l=1}^K \E [S_l|\mathfrak{A}^*]\right|^2\right]&\lesssim \varepsilon^2 M^2+\varepsilon^2 \delta + \left(T^{2-2v}+\delta^2\log^2(K/\delta)+\mathbb{P}(|\mathfrak{A}^*-\theta|>5r/2)\right)\nonumber\\
        &\lesssim \varepsilon^2 \left(T^{2-2v}\{1\vee\log(nK\wedge n\alpha\sqrt{K})\}+\frac{\varepsilon^{-2/v}\wedge K^{2/v}}{n^{2-2/v}}\log^2(K/\delta)\right)\nonumber\\
        & \hspace{4cm}+T^{2-2v}+\delta^2\log^2(K/\delta)+\delta\nonumber\\
        \label{bias-bound-infinite-variance}&\lesssim \{1\vee\log^2(nK)\}\left(\frac{\varepsilon^{2-2/v}}{n^{2-2/v}}\vee\frac{1}{(nK)^{2-2/v}}\vee\frac{1}{(n^2K\alpha^2)^{1-1/v}}\right),
    \end{align}
    where the second inequality holds by plugging in the choice of $M$ \eqref{infinite-var-m-choice} into \eqref{E-indicator-bound}, and the final inequality follows from the choice of $\delta$ in \eqref{infinite-var-m-choice} and the consequent fact that  $\log(K/\delta)\lesssim \log(nK)$. 
    
    \textbf{Case 3: The finite variance setting.} Applying \Cref{E-S-given-J-lemma-finite-moments} to \eqref{bias-decomposition}, we have that
    \begin{align}
        \E \left[\left|\frac{1}{K}\sum_{l=1}^K \E [S_l|\mathfrak{A}^*]\right|^2\right]&\lesssim\varepsilon^2 M^2+\varepsilon^2 \delta+\left(\frac{e^{-2C'_1 nM^2}}{n^2M^2}+(nM)^{2-2v}+T^{2-2v}+\mathbb{P}(|\mathfrak{A}^*-\theta|>5r/2)\right)\nonumber\\
        &\stackrel{(a)}{\lesssim} \varepsilon^2 M^2 + \frac{e^{-2C'_1 nM^2}}{n^2M^2}+(nM)^{2-2v}+T^{2-2v}+\delta\nonumber\\
        &\stackrel{(b)}{\lesssim}\left(\frac{\varepsilon^2}{n} +T^{2-2v}\right)\{1\vee\log(K(n\wedge n^2\alpha^2)) \}+\frac{1}{nK}+\frac{\varepsilon^{2-2/v}}{n^{2-2/v}}+ \frac{1}{(n^2 K\alpha^2)^{1-1/v}}\nonumber\\
        \label{bias-bound-finite-variance}&\stackrel{(c)}{\lesssim} \left(\frac{\varepsilon^2}{n} +\frac{1}{(n^2 K\alpha^2)^{1-1/v}}\right)\{1\vee\log(K(n\wedge n^2\alpha^2))\} +\frac{1}{nK}+\frac{\varepsilon^{2-2/v}}{n^{2-2/v}},
    \end{align}
    where (a) follows from \eqref{E-indicator-bound}, (b) follows from the definitions of $M$ and $\delta$ in \eqref{finite-var-ub-m-choice} and \eqref{delta-choice-finite-variance} respectively, noting that the constant $A_3$ in \eqref{finite-var-ub-m-choice} is chosen to be sufficiently large, and (c) follows from the definition of $T$ in \eqref{finite-var-ub-t-choice}.

\medskip  
    \noindent\textbf{Step 5: The variance term (II).} Observing that $X_{K+l}$ follows a mixture distribution, we have, from the law of total variance, that 
    \begin{align*}\mathrm{Var}(S_l|\mathfrak{A}^*) &=  (1-\varepsilon)\mathrm{Var}_{X_{K+l}\sim P^{\otimes n}}(S_l|\mathfrak{A}^*)+ \varepsilon\mathrm{Var}_{X_{K+l}\sim G_{K+l}}(S_l|\mathfrak{A}^*)\\
    &\hspace{4cm}+\varepsilon(1-\varepsilon)\left(\E_{X_{K+l}\sim P^{\otimes n}}(S_l|\mathfrak{A}^*)-\E_{X_{K+l}\sim G_{K+l}}(S_l|\mathfrak{A}^*)\right)^2\\
    &\leq (1-\varepsilon)\mathrm{Var}_{X_{K+l}\sim P^{\otimes n}}(S_l|\mathfrak{A}^*)+\varepsilon M^2+\varepsilon(1-\varepsilon)(2M)^2\\
    &\lesssim \varepsilon M^2+ \mathrm{Var}_{X_{K+l}\sim P^{\otimes n}}\left(\left.\left[\frac{1}{n}\sum_{i = 1}^n [X_{k,i}]_{-T}^T -\mathfrak{A}^*\right]_{-M}^M\right|\mathfrak{A}^*\right)\\
    &\lesssim \varepsilon M^2+\mathrm{Var}_{X_{K+l}\sim P^{\otimes n}}\left(\frac{1}{n}\sum_{i = 1}^n [X_{K+l,i}]_{-T}^T\right)\\
    &\lesssim \varepsilon M^2+ \frac{1}{n}\mathrm{Var}_{X_{K+l}\sim P^{\otimes n}}([X_{K+l,1}]_{-T}^T),\end{align*}
    where the second and third inequalities follow from the definition of $S_l$ in \eqref{S_k-def}, and the penultimate inequality follows from the fact that $\mathrm{Var}([X]_{-M}^M)\leq \mathrm{Var}(X)$ for any random variable~$X$. As a result, we have 
    \begin{equation}\label{variance-decomposition}
         \E\left[\frac{1}{K^2}\sum_{l=1}^K \mathrm{Var}\left(S_l|\mathfrak{A}^*\right)\right]\lesssim \frac{\varepsilon M^2}{K}+\frac{1}{nK} \mathrm{Var}_{X_{K+1}\sim P^{\otimes n}}([X_{K+1,1}]_{-T}^T) .
    \end{equation}

        \textbf{Case 1: The Gaussian setting.} 
        Observing that $P= \mathcal{N}(\theta,1)$ for some $\theta\in[-1,1]$, we have that
        \[\mathrm{Var}_{X_{K+1}\sim P^{\otimes n}}([X_{K+1,1}]_{-T}^T)\leq \mathrm{Var}_{X_{K+1}\sim P^{\otimes n}}(X_{K+1,1})= 1.\]
        Thus, from \eqref{variance-decomposition}, we have that
        \begin{align}\E\left[\frac{1}{K^2}\sum_{l=1}^K \mathrm{Var}\left(S_l|\mathfrak{A}^*\right)\right]&\lesssim \frac{\varepsilon M^2}{K}+\frac{1}{nK} \mathrm{Var}_{X_{K+1}\sim P^{\otimes n}}([X_{K+1,1}]_{-T}^T) \nonumber\\&\lesssim \varepsilon T^2\left(\frac{1}{nK}\vee \frac{T^2/K}{n^2K\alpha^2}\right)+\frac{1}{nK}  \nonumber\\
        \label{variance-bound-Gaussian}&\lesssim T^2\left(\frac{\varepsilon^2}{n}\vee \frac{1}{n^2K\alpha^2}\right)+\frac{1}{nK}   ,\end{align}
        where in the second inequality we have used that $M=3r$ with $r$ defined in \eqref{UB-thm-tuning-params}, and the final inequality follows from the fact that $T^2/K\lesssim 1$ by \eqref{UB-condition-thm1}.
        
        \textbf{Case 2: The infinite variance setting. }
        Continuing from \eqref{variance-decomposition}, we have that
        \begin{align}
            \E\left[\frac{1}{K^2}\sum_{l=1}^K \mathrm{Var}\left(S_l|\mathfrak{A}^*\right)\right]&\lesssim \frac{\varepsilon M^2}{K}+\frac{1}{nK} \mathrm{Var}_{X_{K+1}\sim P^{\otimes n}}([X_{K+1,1}]_{-T}^T) \nonumber\\
            &\stackrel{(a)}{\lesssim} \frac{\varepsilon M^2}{K}+ \frac{T^{2-v}}{nK}\nonumber\\
            &\stackrel{(b)}{\lesssim} \frac{\varepsilon T^{2-2v}}{K}\{1\vee\log(nK\wedge n\alpha\sqrt{K})\}+\frac{\varepsilon}{K}\cdot\frac{\varepsilon^{-2/v}\wedge K^{2/v}}{n^{2-2/v}}\log^2(K/\delta)+ \frac{T^{2-v}}{nK}\nonumber\\
            &\stackrel{(c)}{\lesssim} T^{2-2v}+\{1\vee\log^2(nK)\}{\left(\frac{\varepsilon^{2-2/v}}{n^{2-2/v}}\vee\frac{1}{(nK)^{2-2/v}}\right)}+\frac{T^{2-v}}{nK}\nonumber\\
            &\label{variance-bound-infinite-variance}\stackrel{(d)}{\lesssim} \{1\vee\log^2(nK)\}\left(\frac{\varepsilon^{2-2/v}}{n^{2-2/v}}\vee\frac{1}{(nK)^{2-2/v}}\right)+ \frac{1}{(n^2K\alpha^2)^{1-1/v}},
        \end{align}
        where (a) follows from \eqref{infinite-var-trunc-lemma-var}, (b) follows from the definition of $M$ in \eqref{infinite-var-m-choice}, (c) follows from assumption \eqref{UB-condition-thm1}, the fact that $\varepsilon/K \lesssim \varepsilon^2 \vee K^{-2}$ and noting that $\log(K/\delta)\lesssim 1\vee\log(nK)$ due to the definition of $\delta$ in \eqref{infinite-var-m-choice}, and (d) follows from the definition of $T$ in \eqref{infinite-var-ub-t-choice}.
        
        \textbf{Case 3: The finite variance setting.} 
        From \eqref{variance-decomposition}, we have the bound 
        \begin{align}
            \E\left[\frac{1}{K^2}\sum_{l=1}^K \mathrm{Var}\left(S_l|\mathfrak{A}^*\right)\right]&\lesssim \frac{\varepsilon M^2}{K}+\frac{1}{nK} \mathrm{Var}_{X_{K+1}\sim P^{\otimes n}}([X_{K+1,1}]_{-T}^T) \nonumber\\
            &\lesssim \left(\frac{\varepsilon }{nK}+\frac{\varepsilon T^{2-2v}}{K}\right)\{1\vee\log(K(n\wedge n^2\alpha^2))\}+ {\frac{\varepsilon^{1-2/v}}{n^{2-2/v}K}}+\frac{1}{nK}\nonumber\\
            \label{variance-bound-finite-variance}&\lesssim \frac{\varepsilon^2}{n}\{1\vee\log(K(n\wedge n^2\alpha^2))\}+\frac{1}{nK}+\frac{1}{(n^2K\alpha^2)^{1-1/v}},
        \end{align}
        where the second inequality follows from \eqref{finite-var-trunc-lemma-var} and the definition of $M$ in \eqref{finite-var-ub-m-choice}, and the final inequality uses the definition of $T$ in \eqref{finite-var-ub-t-choice}, recalls the assumption that $v\geq2$, and notes that $\log(K(n\wedge n^2\alpha^2))\lesssim K$ by assumption \eqref{UB-condition-thm1}.

        \medskip
        \noindent \textbf{Step 6: The mean squared error. } We combine the bounds on the bias and variance to obtain the final result.
        
        \textbf{Case 1: The Gaussian setting. } 
        Plugging \eqref{bias-bound-Gaussian} and \eqref{variance-bound-Gaussian} into \eqref{bias-variance-decomposition}, we have
        \[\E \{(\hat{\theta} - \theta)^2\}\lesssim \frac{1}{nK}+\frac{\varepsilon^2 T^2}{n}+\frac{T^2}{n^2 K\alpha^2}(1+\varepsilon^2T^2).\]

        \textbf{Case 2: The infinite variance setting. }Plugging \eqref{bias-bound-infinite-variance} and \eqref{variance-bound-infinite-variance} into \eqref{bias-variance-decomposition}, and noting that $T\lesssim (n^2K\alpha^2)^{1/(2v)}$ due to \eqref{infinite-var-ub-t-choice}, we have
        \[\E \{(\hat{\theta} - \theta)^2\}\lesssim\{1\vee\log^2(nK)\} \left( \frac{1}{(nK)^{2-2/v}}\vee\frac{1}{(n^2K\alpha^2)^{1-1/v}}\vee\frac{\varepsilon^{2-2/v}}{n^{2-2/v}}\right).\]
  
        \textbf{Case 3: The finite variance setting. } Plugging \eqref{bias-bound-finite-variance} and \eqref{variance-bound-finite-variance} into \eqref{bias-variance-decomposition}, and noting that $T \asymp (n^2K\alpha^2)^{1/(2v)}$ as defined in \eqref{finite-var-ub-t-choice}, we have
        \[\E \{(\hat{\theta} - \theta)^2\}\lesssim \frac{1}{nK}+ \frac{\varepsilon^{2-2/v}}{n^{2-2/v}}+\left(\frac{\varepsilon^2}{n} +\frac{1}{(n^2 K\alpha^2)^{1-1/v}}\right)\{1\vee\log(K(n\wedge n^2\alpha^2))\}.\]

\end{proof}

\subsubsection{Proofs of auxiliary results}\label{sec-UB-lemma-proofs}
In Lemmas~\ref{high-prob-lemma}, \ref{high-prob-lemma-infinite-variance}, and \ref{high-prob-lemma-finite-moments} stated below, we control the initial estimator $\mathfrak{A}^*$ with high probability in the Gaussian, infinite variance, and finite variance settings respectively. All three lemmas make the common assumption that
\begin{equation}\label{eq:high-prob-lemma-common-assumption}\varepsilon+\sqrt{\frac{\log(12/(r\delta))}{K}}+\frac{T\sqrt{\log(12/(r\delta))}}{rn\alpha\sqrt{K}}+\frac{1}{\sqrt{K}(\sqrt{n}\wedge n\alpha)}<C,\end{equation}
where $C>0$ is an absolute constant.

\begin{lemma}\label{high-prob-lemma}
    Consider the Gaussian mean estimation problem detailed in \eqref{Gaussian-distribution} and \eqref{Gaussian-model}. Given $\delta>0$ and tuning parameters $r$ and $T$, suppose that the condition 
    \eqref{eq:high-prob-lemma-common-assumption}
    holds for a sufficiently small absolute constant $C>0$. Then, the output $\mathfrak{A}^*$ of \Cref{alg:mean-localisation} satisfies
    \begin{equation*}\mathbb{P}(|\mathfrak{A}^*-\theta|\leq 5r/2)\geq 1-\delta\end{equation*}
    with input tuning parameters
    \begin{equation}\label{homo-Gaussian-T-condition}
        T =  C_1 \left\{1\vee \sqrt{\log(K(n \wedge n^2 \alpha^2))}\right\} \quad \text{and} \quad r = C_2 T\left(\frac{1}{n}\vee \frac{T^2}{n^2K\alpha^2}\right)^{1/2},
    \end{equation}
    where $C_1, C_2>0$ are sufficiently large absolute constants.

\end{lemma}
\begin{lemma}
    \label{high-prob-lemma-infinite-variance} Consider the heavy-tailed mean estimation problem detailed in \eqref{finite-moments-class} and \eqref{finite-moments-model} with $v\in(1,2]$.
    Given $\delta>0$ and tuning parameters $r$ and $T$, suppose that the condition 
    \eqref{eq:high-prob-lemma-common-assumption}
    holds for a sufficiently small absolute constant $C>0$. Then, the output $\mathfrak{A}^*$ of \Cref{alg:mean-localisation} satisfies
    \begin{equation*}\mathbb{P}(|\mathfrak{A}^*-\theta|\leq 5r/2)\geq 1-\delta\end{equation*}
    with input tuning parameters 
    \begin{equation}\label{T-def-infinite-variance}T = (n^2K\alpha^2)^{1/(2v)}\wedge (nK)^{1/v}\wedge (n/\varepsilon)^{1/v}\end{equation} and \begin{equation}\label{r-def-infinite-variance}r \geq  B_1\left(T^{1-v}\left\{1\vee\sqrt{\log(nK\wedge n\alpha\sqrt{K})}\right\}\vee \frac{\varepsilon^{-1/v}\wedge K^{1/v}}{n^{1-1/v}}\right),\end{equation} where $B_1>0$ is a sufficiently large constant depending only on $v$.
    
\end{lemma}

\begin{lemma}
    \label{high-prob-lemma-finite-moments} Consider the heavy-tailed mean estimation problem detailed in \eqref{finite-moments-class} and \eqref{finite-moments-model} with fixed $v\geq 2$. Given $\delta>0$ and tuning parameters $r$ and $T$, suppose that the condition 
    \eqref{eq:high-prob-lemma-common-assumption}
    holds for a sufficiently small absolute constant $C>0$. Then, the output $\mathfrak{A}^*$ of \Cref{alg:mean-localisation} satisfies
    \begin{equation*}\mathbb{P}(|\mathfrak{A}^*-\theta|\leq 5r/2)\geq 1-\delta\end{equation*}
    with input tuning parameters
    \begin{equation}\label{T-def-finite-variance}T = A_1 (n^2K\alpha^2)^{1/(2v)}
    \quad \text{and} \quad 
    r = A_2\left(\left[1\vee\sqrt{\log(K(n\wedge n^2\alpha^2))}\right]\left\{\frac{1}{\sqrt{n}}\vee T^{1-v}\right\}\vee \frac{\varepsilon^{-1/v}\wedge K^{1/v}}{n^{1-1/v}}\right),\end{equation}
    where $A_1,A_2>0$ are sufficiently large constants depending only on $v$.
\end{lemma}
\begin{proof}[Proof of Lemmas \ref{high-prob-lemma}, \ref{high-prob-lemma-infinite-variance}, and \ref{high-prob-lemma-finite-moments}.]
    The quality of the initial estimator $\mathfrak{A}^*$ will depend on the properties of the smoothed votes $Z_{k,j}$, defined in \eqref{Z_kj-def}. In \textbf{Step 1}, we derive a concentration bound for the smoothed votes around their expectation. In \textbf{Step 2}, we upper bound the expected smoothed votes assigned to bins that are far from $\theta$. In \textbf{Step 3}, we lower bound the expected smoothed votes assigned to bins that are close to $\theta$. Combining the concentration bound in \textbf{Step 1} with the expectation bounds in \textbf{Steps 2} and \textbf{3}, we show that the estimator $\mathfrak{A}^*$ is close to $\theta$ with high probability in \textbf{Step 4}. Since the bound $|\mathfrak{A}^*-\theta|\leq 2+r/2$ always holds, we only need to consider the case $r <1$ for the rest of the proof. 

    \medskip
    \noindent\textbf{Step 1. }Recall that \Cref{alg:mean-localisation} takes each batch $k\in [K]$ and produces privatised smoothed votes $\{Z_{k,j}\}_{j\in[\lceil 2r^{-1}\rceil ]}$ detailed in \eqref{Z_kj-def}.  In this step, we control the deviation of the aggregated smoothed votes $\sum_{k=1}^K Z_{k,j}$ from its expectation. Regarding the added privacy noise, note that $16T\sum_{k=1}^KW_{k,j}/(nK\alpha)$ is sub-exponential with parameters \[\left(32T/\sqrt{n^2K\alpha^2},32T/(nK\alpha)\right).\]
    Using a union bound to control the privacy noise separately, we have that, for any $t>0$ and $\gamma \in (0,1)$,
    \begin{align}
        &\mathbb{P}\left(\left|\frac{1}{K}\sum_{k=1}^K (Z_{k,j}-\mathbb{E}Z_{k,j})\right|\geq t\right)\nonumber\\
        &= \mathbb{P}\left(\left|\frac{1}{K}\sum_{k=1}^K \left\{\left(2r - \left|\frac{1}{n}\sum_{i=1}^n [X_{k,i}]_{-T}^T - a_j\right|\right)_+ + \frac{16 T}{\alpha n} W_{k,j}-\mathbb{E}Z_{k,j}\right\}\right|\geq t\right)\nonumber\\
        &\leq \mathbb{P}\left(\left|\frac{1}{K}\sum_{k=1}^K \left\{\left(2r-\left|\frac{1}{n}\sum_{i=1}^n[X_{k,i}]_{-T}^T-a_j\right|\right)_+-\mathbb{E}Z_{k,j}\right\}\right|\geq (1-\gamma)t\right)+ \mathbb{P}\left(\left|\frac{16T}{nK\alpha}\sum_{k=1}^K W_{k,j}\right|\geq \gamma t\right)\nonumber\\
        &\leq 2\exp\left(-\frac{K(1-\gamma)^2t^2}{2 r^2}\right)+ 2\exp\left(-\frac{K}{64}\min\left\{\frac{n^2\alpha^2\gamma^2 t^2}{32T^2},\frac{n\alpha \gamma t}{T}\right\}\right)\label{eq:smoothed-vote-concentration-bound}.
    \end{align}
    where the final inequality follows from Hoeffding's inequality \citep[e.g. Proposition 2.5 in][]{Wainwright2019high-dimensional} and Bernstein's inequality \citep[e.g. Proposition 2.9 in][]{Wainwright2019high-dimensional}. Define the event
    \begin{equation}\label{high-prob-event}\Omega_t = \left\{\max_{j\in[\lceil 2r^{-1}\rceil ]}\left|\frac{1}{K}\sum_{k=1}^K (Z_{k,j}-\mathbb{E}Z_{k,j})\right|< t\right\}.\end{equation}
    Given $\delta \in (0,1)$, choose $\gamma = (1+rn\alpha/(32T))^{-1}$ and 
    \begin{align}\label{t-choice-complicated}t &= \max\left\{\frac{2r}{(1-\gamma)}\sqrt{\frac{\log(12/(r\delta))}{2K}},\frac{32T}{n\alpha\gamma}\sqrt{\frac{2\log(12/(r\delta))}{K}},\frac{64T}{nK\alpha\gamma}\log(12/(r\delta))\right\}\\
    \label{t-choice-infinite-variance}&= 2\left(r+\frac{32T}{n\alpha}\right)\sqrt{\frac{\log(12/(r\delta))}{2K}},
    \end{align}
    where the final equality follows from the observation that assumption \eqref{eq:high-prob-lemma-common-assumption} (for $C\leq1/\sqrt{2}$) implies that we have $2\log(12/(r\delta))\leq K$. Then, applying a further union bound to \eqref{eq:smoothed-vote-concentration-bound} and plugging in \eqref{t-choice-complicated}, we have that
    \begin{equation}\label{Omega-event-prob-bound}\mathbb{P}\left(\Omega_t^c\right)\leq \frac{6}{r}\exp\left(-\frac{K(1-\gamma)^2t^2}{2r^2}\right)+ \frac{6}{r}\exp\left(-\frac{1}{64}\min\left\{\frac{n^2K\alpha^2\gamma^2 t^2}{32T^2},\frac{nK\alpha \gamma t}{T}\right\}\right)\leq \delta.\end{equation}

    \medskip
    \noindent\textbf{Step 2.} For $j$ such that $|\theta-a_j|>5r/2$, we can upper bound $\mathbb{E}Z_{k,j}$ as follows:
    \begin{align}
        \mathbb{E}Z_{k,j}&\leq 2r\mathbb{P}\left(\left|\frac{1}{n}\sum_{i=1}^n[X_{1,i}]_{-T}^T - a_j\right|\leq 2r\right)\nonumber\\
        &\leq  2r\left\{\varepsilon+(1-\varepsilon)\mathbb{P}_{X_1\sim P^{\otimes n}}\left(\left|\frac{1}{n}\sum_{i=1}^n[X_{1,i}]_{-T}^T - \theta\right|\geq |\theta - a_j| - 2r\right)\right\}\nonumber\\
        \label{bin-weight-UB}&\leq 2r\left\{\varepsilon+ \mathbb{P}_{X_1\sim P^{\otimes n}} \left(\left|\frac{1}{n}\sum_{i=1}^n [X_{1,i}]_{-T}^T -\theta\right|\geq \frac{1}{2}r\right)\right\}.
    \end{align}

We proceed by analysing the three settings separately using the respective assumptions on $P$ and choices of tuning parameters.

\textbf{Case 1: The Gaussian setting. }Recall our choices of the tuning parameters $T$ and $r$ in \eqref{homo-Gaussian-T-condition}. With sufficiently large constants $C_1$ and $C_2$ in \eqref{homo-Gaussian-T-condition}, we have
\begin{align}
    &\mathbb{P}_{X_1\sim \mathcal{N}^{\otimes n}(\theta,1)} \left(\left|\frac{1}{n}\sum_{i=1}^n [X_{1,i}]_{-T}^T -\theta\right|\geq \frac{1}{2}r\right)\nonumber\\
    &\leq \mathbb{P}_{X_1\sim \mathcal{N}^{\otimes n}(\theta,1)} \left(\left|\frac{1}{n}\sum_{i=1}^n [X_{1,i}]_{-T}^T -\E_{X_1\sim \mathcal{N}^{\otimes n}(\theta,1)} \left\{[X_{1,1}]_{-T}^T\right\}\right|\geq \frac{1}{2}r-\left|\E _{X_1\sim \mathcal{N}^{\otimes n}(\theta,1)}\left\{[X_{1,1}]_{-T}^T\right\}-\theta\right|\right)\nonumber\\
    &\stackrel{(a)}{\leq} \mathbb{P}_{X_1\sim \mathcal{N}^{\otimes n}(\theta,1)} \left(\left|\frac{1}{n}\sum_{i=1}^n [X_{1,i}]_{-T}^T -\E_{X_1\sim \mathcal{N}^{\otimes n}(\theta,1)} \left\{[X_{1,1}]_{-T}^T\right\}\right|\geq \frac{1}{2}r-\sqrt{2}e^{-T^2/16}\right)\nonumber\\
    &\stackrel{(b)}{\leq} \mathbb{P}_{X_1\sim \mathcal{N}^{\otimes n}(\theta,1)} \left(\left|\frac{1}{n}\sum_{i=1}^n [X_{1,i}]_{-T}^T -\E_{X_1\sim \mathcal{N}^{\otimes n}(\theta,1)} \left\{[X_{1,1}]_{-T}^T\right\}\right|\geq \frac{1}{2}r-\frac{\sqrt{2}}{\sqrt{K}(\sqrt{n}\wedge n\alpha)}\right)\nonumber\\
    &\stackrel{(c)}{\leq} \mathbb{P}_{X_1\sim \mathcal{N}^{\otimes n}(\theta,1)} \left(\left|\frac{1}{n}\sum_{i=1}^n [X_{1,i}]_{-T}^T -\E_{X_1\sim \mathcal{N}^{\otimes n}(\theta,1)} \left\{[X_{1,1}]_{-T}^T\right\}\right|\geq C_1' r\right)\nonumber\\
    &\stackrel{(d)}{\leq} 2e^{-C_1'^2nr^2/16}\nonumber\\
    &\stackrel{(e)}{\leq} 2e^{-C_2'T^2}\nonumber\\
    &\label{eq:Gaussian-pre-bin-weight-ub}\stackrel{(f)}{\leq} \frac{2}{\sqrt{K}(\sqrt{n}\wedge n\alpha)},
\end{align}
where $C_1'\in (0,1/2)$ and $C_2'>0$ are positive constants, (a) follows from \eqref{Gaussian-trunc-lemma-bias}, (b) and (f) follow from the definition of $T$ in \eqref{homo-Gaussian-T-condition} (provided that the constant $C_1$ is sufficiently large), (c) and (e) follow from the definition of $r$ in \eqref{homo-Gaussian-T-condition} (provided that the constant $C_2$ is sufficiently large), and (d) follows from \eqref{Gaussian-trunc-lemma-concentration}. Plugging \eqref{eq:Gaussian-pre-bin-weight-ub} into \eqref{bin-weight-UB}, we have that
\begin{equation}\label{bin-weight-ub-gaussian}
    \E Z_{k,j}\leq 2r\left(\varepsilon+\frac{2}{\sqrt{K}(\sqrt{n}\wedge n\alpha)}\right).
\end{equation}

\textbf{Case 2: The infinite variance setting. }
Recall our choices of the tuning parameters $T$ and~$r$ in \eqref{T-def-infinite-variance} and \eqref{r-def-infinite-variance} respectively. 
For the rest of the proof, we note that \[T\geq (\sqrt{K}(\sqrt{n}\wedge n\alpha))^{1/v}\wedge \varepsilon^{-1/v}\stackrel{(a)}{\geq}1/C^{1/v}\stackrel{(b)}{\geq} 2,\] where (a) follows from assumption \eqref{eq:high-prob-lemma-common-assumption} and (b) from the assumption that $C>0$ is sufficiently small. With a sufficiently large constant $B_1$ in \eqref{r-def-infinite-variance}, we have that
\begin{align}
        &\mathbb{P}_{X_1\sim P^{\otimes n}} \left(\left|\frac{1}{n}\sum_{i=1}^n [X_{1,i}]_{-T}^T -\theta\right|\geq \frac{1}{2}r\right)\nonumber\\
        &\leq \mathbb{P}_{X_1\sim P^{\otimes n}} \left(\left|\frac{1}{n}\sum_{i=1}^n [X_{1,i}]_{-T}^T -\E_{X_1\sim P^{\otimes n}}\left\{[X_{1,1}]_{-T}^T\right\}\right|\geq \frac{1}{2}r- \left|\E_{X_1\sim P^{\otimes n}}\left\{[X_{1,1}]_{-T}^T\right\}-\theta\right| \right) \nonumber\\
        &\stackrel{(a)}{\leq} \mathbb{P}_{X_1\sim P^{\otimes n}} \left(\left|\frac{1}{n}\sum_{i=1}^n [X_{1,i}]_{-T}^T -\E_{X_1\sim P^{\otimes n}}\left\{[X_{1,1}]_{-T}^T\right\}\right|\geq \frac{1}{2}r- 2T^{1-v}\right) \nonumber\\
        &\stackrel{(b)}{\leq} 2\exp\left(-\frac{n(r/2-2T^{1-v})^2}{2\left\{17T^{2-v}+2T(r/2-2T^{1-v})/3\right\}}\right)\nonumber\\
        &= 2\exp\left(-\frac{3n(r-4T^{1-v})^2}{8T(47T^{1-v}+r)}\right)\nonumber\\
        &\stackrel{(c)}{\leq} 2e^{-nr/(12T)}\nonumber\\
        &\label{eq:infinite-var-pre-bin-weight-ub}\stackrel{(d)}{\leq} 1/4,
    \end{align}
    where (a) follows from \eqref{infinite-var-trunc-lemma-bias}, (b) follows from \eqref{infinite-var-trunc-lemma-concentration} and that $r/2 > 2T^{1-v}$ with $B_1$ sufficiently large in the definition \eqref{r-def-infinite-variance} of $r$, (c) follows from the definition of~$r$ in \eqref{r-def-infinite-variance} (provided that $B_1\geq 47$), and (d) follows from the fact that the definitions of $T$ and $r$ in \eqref{T-def-infinite-variance} and \eqref{r-def-infinite-variance} respectively imply $nr \geq B_1\left\{(nK)^{1/v}\wedge (n/\varepsilon)^{1/v}\right\}\geq B_1 T\geq 12T\log8$. Plugging \eqref{eq:infinite-var-pre-bin-weight-ub} into \eqref{bin-weight-UB}, we have that
    \begin{equation}\label{bin-weight-ub-infinite-variance}
        \E Z_{k,j}\leq 2r(\varepsilon+1/4).
    \end{equation}

\textbf{Case 3: The finite variance setting. }Recall our choices of the tuning parameters $T$ and $r$ in \eqref{T-def-finite-variance}. 
With sufficiently large constants $A_1$ and $A_2$ in \eqref{T-def-finite-variance}, we have
\begin{align}
    &\mathbb{P}_{X_1\sim P^{\otimes n}} \left(\left|\frac{1}{n}\sum_{i=1}^n [X_{1,i}]_{-T}^T -\theta\right|\geq \frac{1}{2}r\right)\nonumber\\
    &\leq \mathbb{P}_{X_1\sim P^{\otimes n}} \left(\left|\frac{1}{n}\sum_{i=1}^n [X_{1,i}]_{-T}^T -\E_{X_1\sim P^{\otimes n}}\left\{[X_{1,1}]_{-T}^T\right\}\right|\geq \frac{1}{2}r- \left|\E_{X_1\sim P^{\otimes n}}\left\{[X_{1,1}]_{-T}^T\right\}-\theta\right| \right) \nonumber\\
    &\stackrel{(a)}{\leq} \mathbb{P}_{X_1\sim P^{\otimes n}} \left(\left|\frac{1}{n}\sum_{i=1}^n [X_{1,i}]_{-T}^T -\E\left\{[X_{1,i}]_{-T}^T\right\}\right|\geq \frac{1}{2}r-\left(\frac{T}{2}\right)^{1-v}\right)\nonumber\\
    &\stackrel{(b)}{\leq}\mathbb{P}_{X_1\sim P^{\otimes n}} \left(\left|\frac{1}{n}\sum_{i=1}^n [X_{1,i}]_{-T}^T -\E\left\{[X_{1,i}]_{-T}^T\right\}\right|\geq A_1'r\right)\nonumber\\
    &\stackrel{(c)}{\leq} 2e^{-A_2' nr^2}+A_3'n^{1-v}r^{-v}\nonumber\\
    \label{eq:finite-var-pre-bin-weight-ub}&\stackrel{(d)}{\leq} \frac{2}{\sqrt{K}(\sqrt{n}\wedge n\alpha)}+\varepsilon+\frac{1}{4K},
\end{align}
where $A_1'\in (0,1/2)$, $A_2'$, and $A_3'$ are positive constants depending only on $v$, (a) follows from \eqref{finite-var-trunc-lemma-bias}, (b) follows from assumption \eqref{eq:high-prob-lemma-common-assumption} and the definition of $r$ in \eqref{T-def-finite-variance} (provided that the constant $A_2$ is sufficiently large), (c) follows from \eqref{finite-var-trunc-lemma-concentration}, and (d) follows from the definition of $r$ in \eqref{T-def-finite-variance} (provided that the constant $A_2$ is sufficiently large). Plugging \eqref{eq:finite-var-pre-bin-weight-ub} into \eqref{bin-weight-UB}, we have that
\begin{equation}\label{bin-weight-ub-finite-variance}
    \E Z_{k,j} \leq 2r\left(2\varepsilon+\frac{2}{\sqrt{K}(\sqrt{n}\wedge n\alpha)}+\frac{1}{4}\right).
\end{equation}

\textbf{The general case. }From \eqref{bin-weight-ub-gaussian}, \eqref{bin-weight-ub-infinite-variance} and \eqref{bin-weight-ub-finite-variance}, the bound
\begin{equation}\label{bin-weight-ub-general}
    \E Z_{k,j}\leq r\left(\frac{1}{2}+4\varepsilon+\frac{4}{\sqrt{K}(\sqrt{n}\wedge n\alpha)}\right)
\end{equation}applies to all three settings.

\medskip
\noindent\textbf{Step 3.} For $j$ such that $|\theta - a_j|\leq r/2$, we can lower bound $\mathbb{E}Z_{k,j}$ as follows:
\begin{align}
    \E Z_{k,j} &\geq (1-\varepsilon)\E_{X_1\sim P^{\otimes n}}\left[2r - \left|\frac{1}{n}\sum_{i=1}^n[X_{1,i}]_{-T}^T - a_j\right|\right]\nonumber+\varepsilon\cdot 0\\
        &\geq (1-\varepsilon) \left(2r-\E_{X_1\sim P^{\otimes n}}\left[\left|\frac{1}{n}\sum_{i=1}^n[X_{1,i}]_{-T}^T-\theta\right|+|\theta-a_j|\right]\right)\nonumber\\
        &\geq (1-\varepsilon) \left(\frac{3}{2}r-\left\{\E_{X_1\sim P^{\otimes n}}\left[\left|\frac{1}{n}\sum_{i=1}^n[X_{1,i}]_{-T}^T-\theta\right|^2\right]\right\}^{1/2}\right)\nonumber\\
        \label{bin-weight-lb}&\geq (1-\varepsilon)\left(\frac{3}{2}r-\left\{\frac{1}{n}\mathrm{Var}_{X_1\sim P^{\otimes n}}\left([X_{1,1}]_{-T}^T\right)+|\E_{X_1\sim P^{\otimes n}}\left\{[X_{1,1}]_{-T}^T\right\}-\theta|^2 \right\}^{1/2}\right)
\end{align}
We again proceed by analysing the three setting separately using the respective assumptions on $P$.

\textbf{Case 1: The Gaussian setting. } Plugging \eqref{Gaussian-trunc-lemma-bias} and \eqref{Gaussian-trunc-lemma-var} into \eqref{bin-weight-lb} and using similar arguments to \textbf{Step 2}, we have that
\begin{align}
    \E Z_{k,j} &\geq (1-\varepsilon)\left(\frac{3}{2}r-\frac{1}{\sqrt{n}}-\sqrt{2}e^{-T^2/16}\right)\nonumber\\
    &\stackrel{(a)}{\geq} (1-\varepsilon)\left(\frac{3}{2}r-\frac{1}{\sqrt{n}}-\frac{\sqrt{2}}{\sqrt{K}(\sqrt{n}\wedge n\alpha)}\right)\nonumber\\
    \label{bin-weight-lb-Gaussian}&\stackrel{(b)}{\geq} (1-\varepsilon)r,
\end{align}
where (a) and (b) follows from the definition of $T$ and $r$ in \eqref{homo-Gaussian-T-condition} respectively, provided that the constants $C_1$ and $C_2$ in their respective definitions are sufficiently large.

\textbf{Case 2: The infinite variance setting. } Plugging \eqref{infinite-var-trunc-lemma-bias} and \eqref{infinite-var-trunc-lemma-var} into \eqref{bin-weight-lb} and using similar arguments to \textbf{Step 2}, we have that
    \begin{align}
    \E Z_{k,j} &\geq (1-\varepsilon)\left(\frac{3}{2}r-\frac{\sqrt{17}T^{1-v/2}}{\sqrt{n}}-2T^{1-v}\right)\nonumber\\
    &\stackrel{(a)}{\geq} (1-\varepsilon)\left(\frac{3}{2}r-\frac{\sqrt{17}(\varepsilon^{-1}\wedge K)^{1/v-1/2}}{n^{1-1/v}}-2T^{1-v}\right)\nonumber\\
    \label{bin-weight-lb-infinite-var}&\stackrel{(b)}{\geq} (1-\varepsilon) r,
\end{align}
where (a) and (b) follows from the definitions of $T$ and $r$ in \eqref{T-def-infinite-variance} and~\eqref{r-def-infinite-variance} respectively, provided that the constant $B_1$ in \eqref{r-def-infinite-variance} is sufficient large.

\textbf{Case 3: The finite variance setting. } Plugging \eqref{finite-var-trunc-lemma-bias} and \eqref{finite-var-trunc-lemma-var} into \eqref{bin-weight-lb} and using similar arguments to \textbf{Step 2}, we have that 
\begin{align}
    \E Z_{k,j} &\geq (1-\varepsilon)\left(\frac{3}{2}r-\frac{1}{\sqrt{n}}-\left(\frac{T}{2}\right)^{1-v}\right)\nonumber\\
    \label{bin-weight-lb-finite-var}&\geq (1-\varepsilon) r,
\end{align}
where the last inequality follows from assumption \eqref{eq:high-prob-lemma-common-assumption} and the definition of $r$ in \eqref{T-def-finite-variance}, provided that the constant $A_2$ in \eqref{T-def-finite-variance} is sufficiently large.

    \textbf{The general case. } From \eqref{bin-weight-lb-Gaussian}, \eqref{bin-weight-lb-infinite-var} and \eqref{bin-weight-lb-finite-var}, the bound
    \begin{equation}\label{bin-weight-lb-general}\E Z_{k,j} \geq (1-\varepsilon)r\end{equation}
    applies to all three settings.

\medskip    
\noindent\textbf{Step 4.} On the event $\Omega_t$ defined in \eqref{high-prob-event}, observe that we have
\[\frac{1}{K}\sum_{k=1}^K Z_{k,j}-\frac{1}{K}\sum_{k=1}^K Z_{k,l}> (\E Z_{k,j}-t)-(\E Z_{k,l}+t)\geq r\left(\frac{1}{2}-5\varepsilon-\frac{4}{\sqrt{K}(\sqrt{n}\wedge n\alpha)}\right)-2t\]
for $j$ satisfying $|\theta-a_j|\leq r/2$ and $l$ satisfying $|\theta-a_l|> 5r/2$, where the final inequality follows from \eqref{bin-weight-ub-general} and \eqref{bin-weight-lb-general}.
Hence, plugging in the choice of $t$ in \eqref{t-choice-infinite-variance}, we have $|\mathfrak{A}^*-\theta |\leq 5r/2$ on the event $\Omega_t$ provided that
    \[2t= 4\left(r+\frac{32T}{n\alpha}\right)\sqrt{\frac{\log(12/(r\delta))}{2K}}<r\left(\frac{1}{2}-5\varepsilon-\frac{4}{\sqrt{K}(\sqrt{n}\wedge n\alpha)}\right).\]
    Rearranging this, we conclude by \eqref{Omega-event-prob-bound} that there exists a constant $C>0$ such that the condition
    \[\sqrt{\frac{\log(12/(r\delta))}{K}}+\frac{T\sqrt{\log(12/(r\delta))}}{rn\alpha\sqrt{K}}+\varepsilon+\frac{1}{\sqrt{K}(\sqrt{n}\wedge n\alpha)}<C\]
is sufficient to guarantee that $\mathbb{P}(|\mathfrak{A}^*-\theta|\leq 5r/2)\geq 1-\delta$. 
\end{proof}

We state some basic properties of the truncated mean under the Gaussian, infinite variance, and finite variance settings in Lemmas \ref{trunc-lemma}, \ref{trunc-lemma-infinite-variance}, and \ref{trunc-lemma-finite-moments} respectively.

\begin{lemma}\label{trunc-lemma}
    For $i\in [n]$, let $X_{i}\overset{\mathrm{i.i.d.}}{\sim} \mathcal{N}(\theta,1)$ for some $\theta \in [-1,1]$. Suppose that the truncation parameter $T$ satisfies $T \geq 2$. Then, the truncated mean $n^{-1}\sum_{i=1}^n [X_i]_{-T}^T$ satisfies the following properties:
    \begin{enumerate}
    
    \item (Bias.) \begin{equation}\label{Gaussian-trunc-lemma-bias}|\E\left\{[X_{1}]_{-T}^T\right\}-\theta|\leq \sqrt{2}e^{-T^2/16}.\end{equation}
    \item (Variance.) \begin{equation}\label{Gaussian-trunc-lemma-var}\mathrm{Var}\left([X_{1}]_{-T}^T\right)\leq 1.\end{equation}
    \item (Concentration bound.) For any $t\geq 0$,  \begin{equation}\label{Gaussian-trunc-lemma-concentration}\mathbb{P}\left(\left|\frac{1}{n}\sum_{i=1}^n[X_{i}]_{-T}^T - \E \left\{[X_{1}]_{-T}^T\right\}\right|\geq t\right)\leq 2e^{-nt^2/16}.\end{equation}
    \end{enumerate}

\end{lemma}
\begin{proof}

\begin{enumerate}
\item Since $T\geq 2 \geq 2|\theta|$, we have that \begin{align*}
        |\E\{ [X_1]_{-T}^T\} - \theta| &= |\E\{ [X_1]_{-T}^T- X_1\}|\\
        &\leq \E[ \mathbf{1}\{|X_1|>T\} |\sign(X_{1})T - X_{1}| ]\\
        &\leq \E[ \mathbf{1}\{|X_{1}|>T\} |\theta - X_{1}| ]\\
        &\leq (\mathbb{P}(|X_{1}|>T))^{1/2}(\E [|X_{1}-\theta|^2])^{1/2}\\
        &\leq (\mathbb{P}(|X_{1}-\theta|>T-|\theta|))^{1/2}\\
        &\leq (\mathbb{P}(|X_{1}-\theta|>T/2))^{1/2}\\
        &\leq \sqrt{2}e^{-T^2/16},
    \end{align*}
    where the third inequality follows from the Cauchy--Schwarz inequality and the last inequality follows from Hoeffding's inequality.
    \item This follows from the general fact that \[\mathrm{Var}([X]_{-T}^T)\leq \mathrm{Var}(X)\]
    for any random variable $X$. Let $X'$ be an independent copy of $X$, $Y = [X]_{-T}^T$, and $Y'= [X']_{-T}^T$. Then, we have that
    \[2\mathrm{Var}([X]_{-T}^T)= \E[(Y-Y')^2]\leq \E[(X-X')^2]= 2\mathrm{Var}(X).\]
    \item We first show that $[X_{1}]_{-T}^T - \E\left\{[X_{1}]_{-T}^T\right\} $ is sub-Gaussian with variance parameter 8. Note that for $T\geq2$,
\begin{align*}
\E \left| [X_{1}]_{-T}^T - \E\left\{[X_{1}]_{-T}^T\right\}\right|^{2q}&\leq 2^{2q-1}\left(\E\left| [X_{1}]_{-T}^T - \theta\right|^{2q} +\left| \theta - \E\left\{[X_{1}]_{-T}^T\right\}\right|^{2q}\right)\\
&\leq 2^{2q-1}\left(\E\left| X_{1} - \theta\right|^{2q} +\left| \theta - \E\left\{[X_{1}]_{-T}^T\right\}\right|^{2q}\right)\\
&\leq 2^{2q-1}\left(\frac{(2q)!}{2^q q!} +2^q e^{-qT^2/8}\right)\\
&\leq 2^{2q} \frac{(2q)!}{2^q q!}
\end{align*}
for every integer $q \geq 2$, where the penultimate inequality follows from \eqref{Gaussian-trunc-lemma-bias}. Further, the $q=1$ case follows from \eqref{Gaussian-trunc-lemma-var}. Sub-Gaussianity then follows from standard properties \citep[e.g.~Theorem 2.6 in][]{Wainwright2019high-dimensional}. The result follows from Hoeffding's inequality \citep[e.g.~Proposition 2.5 in][]{Wainwright2019high-dimensional}. 

    \end{enumerate}
\end{proof}

\begin{lemma}\label{trunc-lemma-infinite-variance}
    For fixed $v \in (1,2]$ and $i\in [n]$, let $X_{i}\overset{\mathrm{i.i.d.}}{\sim} P\in \mathcal{P}_v$, where $\mathcal{P}_v$ is defined by \eqref{finite-moments-class}, and suppose $\theta = \E X_{i}\in [-1,1]$. Suppose that the truncation parameter $T$ satisfies $T\geq2$. Then, the truncated mean $n^{-1}\sum_{i=1}^n [X_i]_{-T}^T$ satisfies the following properties:
    \begin{enumerate}
    \item (Bias.)\begin{equation}\label{infinite-var-trunc-lemma-bias}\left|\E\left\{[X_{1}]_{-T}^T\right\}-\theta\right|\leq 1\wedge(T-|\theta|)^{1-v}\leq 1\wedge2^{v-1}T^{1-v}\end{equation}
    \item (Variance.) \begin{equation}\label{infinite-var-trunc-lemma-var}\mathrm{Var}([X_{1}]_{-T}^T)\leq\E\left[([X_{1}]_{-T}^T-\theta)^2 \right]\leq 4(2^{v}+2^{-v})T^{2-v}\leq 17T^{2-v}\end{equation}
       
        \item (Concentration bound.) For any $t\geq 0$, \begin{equation}\label{infinite-var-trunc-lemma-concentration}\mathbb{P}\left(\left|\frac{1}{n}\sum_{i = 1}^n [X_{i}]_{-T}^T - \E \left\{[X_1]_{-T}^T\right\}\right|\geq t\right)\leq 2\exp\left(-\frac{nt^2}{2\{17T^{2-v}+2Tt/3\}}\right)\end{equation}

    \end{enumerate}
\end{lemma}
\begin{proof}
    \begin{enumerate}
    \item Since $T\geq 2>|\theta|$, we have that
    \begin{align*}
        |\E\{ [X_{K+1, 1}]_{-T}^T\} - \theta| &\leq \E[ \mathbf{1}\{|X_{K+1,1}|>T\} |\sign(X_{K+1,1})T - X_{K+1,1}| ]\\
        &\leq \E[ \mathbf{1}\{|X_{K+1,1}|>T\} |\theta - X_{K+1,1}| ]\\
        &\leq (\mathbb{P}(|X_{K+1,1}|>T))^{1-1/v}(\E [|X_{K+1,1}-\theta|^v])^{1/v}\\
        &\leq (\mathbb{P}(|X_{K+1,1}-\theta|>T-|\theta|))^{1-1/v}\\
        &\leq 1 \wedge (T-|\theta|)^{1-v},
    \end{align*}
    where the third inequality follows from H\"older's inequality.
    
    \item Since $T\geq2|\theta|$ and $v\in (1,2]$, 
    \begin{align*}&\E\left[([X_{1}]_{-T}^T-\theta)^2 \right] \\
    &\leq \E \left[ 4T^2\mathbf{1}\{|X_{1}|>T\}+|X_{1}-\theta|^v|X_{1}-\theta|^{2-v}\mathbf{1}\{|X_{1}|\leq T\}\right]\\
    &\leq 4T^2(T-|\theta|)^{-v}+(2T)^{2-v}\\
    &\leq 4T^2(T/2)^{-v}+(2T)^{2-v},
    \end{align*}
    where the second inequality follows from Markov's inequality.
       
    \item We use Bernstein's inequality, specialised to bounded random variables \citep[e.g.~Proposition 2.14 in][]{Wainwright2019high-dimensional}, to obtain
    \begin{align*}
        \mathbb{P}\left(\left|\frac{1}{n}\sum_{i = 1}^n [X_{i}]_{-T}^T - \E \{[X_{i}]_{-T}^T\}\right|\geq t\right)&\leq 2\exp\left(-\frac{nt^2}{2\{\mathrm{Var}([X_{i}]_{-T}^T)+2Tt/3\}}\right).
    \end{align*}
    The results follows from \eqref{infinite-var-trunc-lemma-var}.
    
    \end{enumerate}
\end{proof}
\begin{lemma}\label{trunc-lemma-finite-moments}
    For fixed $v\geq 2$ and $i\in [n]$, let  $X_{i}\overset{\mathrm{i.i.d.}}{\sim} P\in \mathcal{P}_v$, where $\mathcal{P}_v$ is defined by \eqref{finite-moments-class}, and suppose $\theta=\E X_{i}\in [-1,1]$. Suppose that the truncation parameter $T$ satisfies $T\geq 2$. Then, the truncated mean $n^{-1}\sum_{i=1}^n [X_i]_{-T}^T$ satisfies the following properties:
    \begin{enumerate}
    \item (Bias.) \begin{equation}\label{finite-var-trunc-lemma-bias}\left|\E\left\{[X_{1}]_{-T}^T\right\}-\theta\right|\leq 1\wedge(T-|\theta|)^{1-v}\leq 1\wedge2^{v-1}T^{1-v}.\end{equation}
    \item (Variance.) \begin{equation}\label{finite-var-trunc-lemma-var}\mathrm{Var}([X_{1}]_{-T}^T)\leq1.\end{equation}

    \item (Concentration bound.) For any $t\geq 0$, \begin{equation}\label{finite-var-trunc-lemma-concentration}\mathbb{P}\left(\left|\frac{1}{n}\sum_{i=1}^n\left([X_{K+1,i}]_{-T}^T-\E\left\{[X_{K+1,i}]_{-T}^T\right\}\right)\right|\geq x\right)
    \leq 2e^{-C_1nx^2}+C_2 n^{1-v}x^{-v},\end{equation}
    where $C_1$ and $C_2$ are positive constants depending only on $v$.

    \end{enumerate}
\end{lemma}

\begin{proof}[Proof of \Cref{trunc-lemma-finite-moments}]
\begin{enumerate}
\item Same as \Cref{trunc-lemma-infinite-variance}.
\item Following from the proof of \eqref{Gaussian-trunc-lemma-var} in \Cref{trunc-lemma}, we have that \begin{align*}\mathrm{Var}\left([X_{K+1,1}]_{-T}^T\right)&\leq \mathrm{Var}\left(X_{K+1,1}\right)\\ &\leq (\E |X_{K+1,1}-\theta|^v)^{2/v}\\&\leq 1,\end{align*}
where the second inequality follows from Jensen's inequality.
    \item By a version of the Fuk--Nagaev inequality \citep[e.g.~Equation (1) in][]{jirak2025concentration}, we have, for positive constants $C_1$, $C_2$, and $C_2'$ that depend only on $v$,
    \begin{align*}
    &\mathbb{P}\left(\left|\frac{1}{n}\sum_{i=1}^n\left([X_{K+1,i}]_{-T}^T-\E\left\{[X_{K+1,i}]_{-T}^T\right\}\right)\right|\geq x\right)\\\leq &1 \wedge \left\{2e^{-C_1nx^2}+C_2(n^{1-v}x^{-v})^2+ \mathbb{P}\left(\max_{i\in [n]}\left|[X_{K+1,i}]_{-T}^T-\E\left\{[X_{K+1,i}]_{-T}^T\right\}\right|>\frac{nx}{4}\right)\right\}\\
    \leq&1 \wedge \left\{2e^{-C_1nx^2}+C_2(n^{1-v}x^{-v})^2+8^v n^{1-v}x^{-v}\right\}\\
    \leq &2e^{-C_1nx^2}+C_2' n^{1-v}x^{-v},
    \end{align*}
    where the second inequality follows from applying Markov's inequality with a union bound, and the final inequality follows from the observation that the bound is vacuous unless \[(n^{1-v}x^{-v})^2\leq n^{1-v}x^{-v}\leq 1.\] 
    \end{enumerate}
\end{proof}

Recall that the control of the quantity $S_l$, defined in \eqref{S_k-def}, is critical in the proof of Propositions~\ref{UB-thm},~\ref{finite-moment-UB} and~\ref{infinite-variance-UB}. In Lemmas \ref{E-S-given-J-lemma}, \ref{E-S-given-J-lemma-infinite-variance} and \ref{E-S-given-J-lemma-finite-moments} stated below, we control, for clean batches, its conditional expectation given the initial estimator $\mathfrak{A}^*$ under the Gaussian, infinite variance, and finite variance settings respectively.
\begin{lemma}\label{E-S-given-J-lemma}
    Let $X_{K+1}\sim \mathcal{N}^{\otimes n}(\theta,1)$ for some $\theta \in [-1,1]$. Assume that $n^2K\alpha^2\geq 1$ and suppose that we choose tuning parameters \begin{equation}\label{eq:gaussian-E-S-given-J-lemma-tuning-params}
        T =  C_1 \left\{1\vee \sqrt{\log(K(n \wedge n^2 \alpha^2))}\right\}, \quad r = C_2 T\left(\frac{1}{n}\vee \frac{T^2}{n^2K\alpha^2}\right)^{1/2}, \quad \text{and} \quad M=3r,
    \end{equation}
    where $C_1, C_2>0$ are sufficiently large absolute constants.  Then, for any $\mathfrak{A}^*\in [-1,1+r/2]$,
    \[|\mathbb{E}[S_1|\mathfrak{A}^*]|\lesssim \frac{1}{\sqrt{K}(\sqrt{n}\wedge n\alpha)}+\mathbf{1}\{|\mathfrak{A}^*-\theta|>5r/2\}.\]
\end{lemma}
\begin{proof}
    By considering the three possible cases of the truncation, 
    \begin{multline}\label{Gaussian-three-cases-truncation}
        |\mathbb{E}[S_1|\mathfrak{A}^*]|\leq \mathbb{E}\left[\left|\mathfrak{A}^*-\theta+M\right|\mathbf{1}\{R_1>M\}
        +\left|\mathfrak{A}^*-\theta-M\right|\mathbf{1}\{R_1<-M\}\right|\mathfrak{A}^*]\\+\left.\left|\E\left[\frac{1}{n}\sum_{i=1}^n([X_{K+1,i}]_{-T}^T-\theta) \mathbf{1}\{|R_1|\leq M\}\right|\mathfrak{A}^*\right]\right|.
    \end{multline}
    Note that
    \begin{align*}
        &\left.\left|\E\left[\frac{1}{n}\sum_{i=1}^n([X_{K+1,i}]_{-T}^T-\theta) \mathbf{1}\{|R_1|\leq M\}\right|\mathfrak{A}^*\right]\right|\\
        &=\left.\left|\E\left[\frac{1}{n}\sum_{i=1}^n([X_{K+1,i}]_{-T}^T-\theta) (1-\mathbf{1}\{|R_1|> M\})\right|\mathfrak{A}^*\right]\right|\\
        &\leq \left|\E\left\{[X_{K+1,1}]_{-T}^T-\theta \right\}\right|+\left.\left|\E\left[\frac{1}{n}\sum_{i=1}^n([X_{K+1,i}]_{-T}^T-\theta) \mathbf{1}\{|R_1|> M\}\right|\mathfrak{A}^*\right]\right|\\
        &\lesssim e^{-T^2/16}+\left(\E \left|\frac{1}{n}\sum_{i=1}^n [X_{K+1,i}]_{-T}^T -\theta\right|^2\right)^{1/2}\left\{\mathbb{P}\left(\left.\left|\frac{1}{n}\sum_{i = 1}^n [X_{K+1,i}]_{-T}^T - \mathfrak{A}^*\right|>M\right|\mathfrak{A}^*\right)\right\}^{1/2}\\
        &\lesssim e^{-T^2/16}+\frac{1}{n^{1/2}}\left\{\mathbb{P}\left(\left.\left|\frac{1}{n}\sum_{i = 1}^n [X_{K+1,i}]_{-T}^T - \mathfrak{A}^*\right|>M\right|\mathfrak{A}^*\right)\right\}^{1/2}\\
        &=e^{-T^2/16}+\frac{1}{n^{1/2}}\left\{\mathbb{P}\left(\left.\left|R_1\right|>M\right|\mathfrak{A}^*\right)\right\}^{1/2},
    \end{align*}
    where the penultimate inequality follows from \eqref{Gaussian-trunc-lemma-bias} and the Cauchy--Schwarz inequality, the final inequality follows from \eqref{Gaussian-trunc-lemma-bias} and \eqref{Gaussian-trunc-lemma-var} using the bias-variance decomposition, and the final equality follows from the definition of $R_l$ in \eqref{R_k-def}. Hence, noting that $r\vee n^{-1/2}\lesssim M$ from \eqref{eq:gaussian-E-S-given-J-lemma-tuning-params} and $|\mathfrak{A}^*-\theta|\leq 2+r/2$, \eqref{Gaussian-three-cases-truncation} becomes 
    \begin{multline}\label{Gaussian-three-cases-truncation-2}|\E[S_1|\mathfrak{A}^*]|\\\lesssim e^{-T^2/16}+((M+1)\mathbf{1}\{|\mathfrak{A}^*-\theta|>5r/2\}+M\mathbf{1}\{|\mathfrak{A}^*-\theta|\leq5r/2\})\left\{\mathbb{P}\left(\left.\left|R_1\right|>M\right|\mathfrak{A}^*\right)\right\}^{1/2}.\end{multline}
Applying the definition of $R_l$ in \eqref{R_k-def}, note that, provided that the absolute constants $C_1$ and $C_2$ are sufficiently large, 
    \begin{align}
        &\mathbb{P}\left(\left.\left|R_1\right|>M\right|\mathfrak{A}^*\right)\nonumber\\
        &=\mathbb{P}\left(\left.\left|\frac{1}{n}\sum_{i = 1}^n [X_{K+1,i}]_{-T}^T - \mathfrak{A}^*\right|>M\right|\mathfrak{A}^*\right)\nonumber\\
        &\leq \mathbb{P}\left(\left.\left|\frac{1}{n}\sum_{i = 1}^n [X_{K+1,i}]_{-T}^T - \E\left\{[X_{K+1,1}]_{-T}^T\right\}\right|>M-|\mathfrak{A}^*-\theta|-\left|\E\left\{[X_{K+1,1}]_{-T}^T-\theta \right\}\right|\right|\mathfrak{A}^*\right)\nonumber\\
        &\stackrel{(a)}{\lesssim} \mathbf{1}\{|\mathfrak{A}^*-\theta|\leq5r/2\}\exp\left(-n(M-5r/2-\sqrt{2}e^{-T^2/16})^2/16\right)\nonumber\\
        &\hspace{1cm}+\mathbf{1}\{|\mathfrak{A}^*-\theta|>5r/2\}\left[\exp\left(-n(M-2-r/2-\sqrt{2}e^{-T^2/16})^2/16\right)\mathbf{1}\{M>3\}+\mathbf{1}\{M\leq 3\}\right]\nonumber\\
        &\stackrel{(b)}{\leq} e^{-n(M/12)^2/16}+\mathbf{1}\{|\mathfrak{A}^*-\theta|>5r/2\}\mathbf{1}\{M\leq 3\}\nonumber\\
        \label{eq:gaussian-E-S-given-J-prob-bound}&\stackrel{(c)}{\leq} e^{-C_2^2T^2/256}+\mathbf{1}\{|\mathfrak{A}^*-\theta|>5r/2\}\mathbf{1}\{M\leq 3\},
    \end{align}
    where (a) follows from \eqref{Gaussian-trunc-lemma-bias} and \eqref{Gaussian-trunc-lemma-concentration}, (b) follows from the facts that $M=3r$ and
    \[\sqrt{2}e^{-T^2/16}<3/4 \wedge r/4\] for $C_1$ sufficiently large in the definition of $T$ in \eqref{eq:gaussian-E-S-given-J-lemma-tuning-params}, and (c) follows from the definition of $M$ and $r$ in \eqref{eq:gaussian-E-S-given-J-lemma-tuning-params}. Substituting \eqref{eq:gaussian-E-S-given-J-prob-bound} into \eqref{Gaussian-three-cases-truncation-2},
    \begin{align*}|\E[S_1|\mathfrak{A}^*]|&\lesssim e^{-T^2/16}+Me^{-C_2^2T^2/512}+\{(M+1)e^{-C_2^2T^2/512}+1\}\mathbf{1}\{|\mathfrak{A}^*-\theta|>5r/2\}\\
    &\lesssim \frac{1}{\sqrt{K}(\sqrt{n}\wedge n\alpha)}+\mathbf{1}\{|\mathfrak{A}^*-\theta|>5r/2\},\end{align*}
    provided that the constants $C_1$ and $C_2$ are sufficiently large.
\end{proof}
\begin{lemma}\label{E-S-given-J-lemma-infinite-variance} For fixed $v \in (1,2]$, let $X_{K+1}\sim P^{\otimes n}$ with $P\in \mathcal{P}_v$, defined by \eqref{finite-moments-class}, and suppose $\theta = \E X_{K+1,1}\in [-1,1]$. Assume that $n^2K\alpha^2\geq 1$ and suppose that we choose tuning parameters 
\begin{equation*}
        T = \left(nK\wedge n\alpha\sqrt{K}\wedge  n/\varepsilon\right)^{1/v},
    \end{equation*} 
    \[r =  B_1\left(T^{1-v}\left\{1\vee\sqrt{\log(nK\wedge n\alpha\sqrt{K})}\right\}\vee \frac{\varepsilon^{-1/v}\wedge K^{1/v}}{n^{1-1/v}}\log(K/\delta)\right),\]
    \[M = B_2 r, \quad \text{and} \quad \delta \asymp \left(\frac{1}{(nK)^{2-2/v}} \vee  \frac{1}{(n^2 K\alpha^2)^{1-1/v}}\right),\]
    where $B_1,B_2>0$ are sufficiently large absolute constants.
Then, for any $\mathfrak{A}^*\in [-1,1+r/2]$ we have that
    \[|\E[S_1|\mathfrak{A}^*]| \lesssim T^{1-v}+\delta\log(K/\delta)+\mathbf{1}\{|\mathfrak{A}^*-\theta|>5r/2\}\]
\end{lemma}
\begin{proof}
    By considering the three possible cases of the truncation, 
    \begin{multline}\label{infinite-var-three-cases-truncation}
        |\mathbb{E}[S_1|\mathfrak{A}^*]|\leq \mathbb{E}\left[\left|\mathfrak{A}^*-\theta+M\right|\mathbf{1}\{R_1>M\}
        +\left|\mathfrak{A}^*-\theta-M\right|\mathbf{1}\{R_1<-M\}\right|\mathfrak{A}^*]\\+\left.\left|\E\left[\frac{1}{n}\sum_{i=1}^n([X_{K+1,i}]_{-T}^T-\theta) \mathbf{1}\{|R_1|\leq M\}\right|\mathfrak{A}^*\right]\right|.
    \end{multline}
    Note that
    \begin{align*}
        &\left.\left|\E\left[\frac{1}{n}\sum_{i=1}^n([X_{K+1,i}]_{-T}^T-\theta) \mathbf{1}\{|R_1|\leq M\}\right|\mathfrak{A}^*\right]\right|\\
        &=\left.\left|\E\left[\frac{1}{n}\sum_{i=1}^n([X_{K+1,i}]_{-T}^T-\theta) (1-\mathbf{1}\{|R_1|> M\})\right|\mathfrak{A}^*\right]\right|\\
        &\leq \left|\E\left\{[X_{K+1,1}]_{-T}^T-\theta \right\}\right|+\left.\left|\E\left[\frac{1}{n}\sum_{i=1}^n([X_{K+1,i}]_{-T}^T-\theta) \mathbf{1}\{|R_1|> M\}\right|\mathfrak{A}^*\right]\right|\\
        &\lesssim T^{1-v}+\left(\E \left|\frac{1}{n}\sum_{i=1}^n [X_{K+1,i}]_{-T}^T -\theta\right|^2\right)^{1/2}\left\{\mathbb{P}\left(\left.\left|\frac{1}{n}\sum_{i = 1}^n [X_{K+1,i}]_{-T}^T - \mathfrak{A}^*\right|>M\right|\mathfrak{A}^*\right)\right\}^{1/2}\\
        &\lesssim T^{1-v}+\frac{T^{1-v/2}}{n^{1/2}}\left\{\mathbb{P}\left(\left.\left|\frac{1}{n}\sum_{i = 1}^n [X_{K+1,i}]_{-T}^T - \mathfrak{A}^*\right|>M\right|\mathfrak{A}^*\right)\right\}^{1/2},
    \end{align*}
    where the penultimate inequality follows from \eqref{infinite-var-trunc-lemma-bias} and the Cauchy--Schwarz inequality, and the final inequality follows from \eqref{infinite-var-trunc-lemma-bias} and \eqref{infinite-var-trunc-lemma-var} using the bias-variance decomposition.
    Hence, noting that $r\vee T^{1-v/2}/n^{1/2}\lesssim M$ and $|\mathfrak{A}^*-\theta|\leq 2+r/2$, \eqref{infinite-var-three-cases-truncation} becomes
    \begin{multline}\label{infinite-var-three-cases-truncation-2}|\E[S_1|\mathfrak{A}^*]|\lesssim T^{1-v}+\\((M+1)\mathbf{1}\{|\mathfrak{A}^*-\theta|>5r/2\}+M\mathbf{1}\{|\mathfrak{A}^*-\theta|\leq5r/2\})\left\{\mathbb{P}\left(\left.\left|\frac{1}{n}\sum_{i = 1}^n [X_{K+1,i}]_{-T}^T - \mathfrak{A}^*\right|>M\right|\mathfrak{A}^*\right)\right\}^{1/2}.\end{multline}
    Provided that the absolute constants $B_1$ and $B_2$ are sufficiently large, we have, for some absolute constant $B_3>0$, 
    \begin{align}
        &\mathbb{P}\left(\left.\left|\frac{1}{n}\sum_{i = 1}^n [X_{K+1,i}]_{-T}^T - \mathfrak{A}^*\right|>M\right|\mathfrak{A}^*\right)\nonumber\\
        &\leq \mathbb{P}\left(\left.\left|\frac{1}{n}\sum_{i = 1}^n [X_{K+1,i}]_{-T}^T - \E\left\{[X_{K+1,1}]_{-T}^T\right\}\right|>M-|\mathfrak{A}^*-\theta|-\left|\E\left\{[X_{K+1,1}]_{-T}^T-\theta \right\}\right|\right|\mathfrak{A}^*\right)\nonumber\\
        &\stackrel{(a)}{\lesssim} \mathbf{1}\{|\mathfrak{A}^*-\theta|\leq5r/2\}\exp\left(-\frac{n(M-5r/2-2T^{1-v})^2}{2\{17T^{2-v}+2T(M-5r/2-2T^{1-v})/3\}}\right)\nonumber\\
        &\hspace{0.5cm}+\mathbf{1}\{|\mathfrak{A}^*-\theta|>5r/2\}\left(\exp\left(-\frac{n(M-2-r/2-2T^{1-v})^2}{2\{17T^{2-v}+2T(M-2-r/2-2T^{1-v})/3\}}\right)\mathbf{1}\{M> 3\}+\mathbf{1}\{M\leq 3\}\right)\nonumber\\
        &\leq \mathbf{1}\{|\mathfrak{A}^*-\theta|\leq5r/2\}\exp\left(-\frac{1}{2}\left\{\frac{n(M-5r/2-2T^{1-v})^2}{34T^{2-v}}\wedge \frac{n(M-5r/2-2T^{1-v})}{4T/3}\right\}\right) \nonumber\\
        &\hspace{0.5cm}+\mathbf{1}\{|\mathfrak{A}^*-\theta|>5r/2\}\bigg(\exp\left(-\frac{1}{2}\left\{\frac{n(M-2-r/2-2T^{1-v})^2}{34T^{2-v}}\wedge \frac{n(M-2-r/2-2T^{1-v})}{4T/3}\right\}\right)\mathbf{1}\{M> 3\}\nonumber\\
        &\hspace{8cm}+\mathbf{1}\{M\leq 3\}\bigg)\nonumber\\
        &\stackrel{(b)}{\leq}e^{-B_3 nM/T}+\mathbf{1}\{|\mathfrak{A}^*-\theta|>5r/2\}\mathbf{1}\{M\leq 3\}\nonumber\\
        \label{eq:infinite-var-E-S-given-J-prob-bound}&\stackrel{(c)}{\leq} \delta^2/K^2 + \mathbf{1}\{|\mathfrak{A}^*-\theta|>5r/2\}\mathbf{1}\{M\leq 3\}
    \end{align}
    where (a) follows from \eqref{infinite-var-trunc-lemma-bias} and \eqref{infinite-var-trunc-lemma-concentration}, (b) and (c) follows from the fact that \[M\gtrsim r \vee T^{1-v}\vee T\log(K/\delta)/n\]
    and that the absolute constants $B_1$ and $B_2$ are sufficiently large. Substituting \eqref{eq:infinite-var-E-S-given-J-prob-bound} into \eqref{infinite-var-three-cases-truncation-2},

\begin{align*}|\E[S_1|\mathfrak{A}^*]|&\lesssim T^{1-v}+M\delta/K +(M+1)\mathbf{1}\{|\mathfrak{A}^*-\theta|>5r/2\}\left(\delta^2K^{-2}+\mathbf{1}\{M\leq 3\}\right)^{1/2}\\
&\lesssim T^{1-v}+M\delta/K+\mathbf{1}\{|\mathfrak{A}^*-\theta|>5r/2\}\\
&\lesssim T^{1-v}+\delta\log(K/\delta)+\mathbf{1}\{|\mathfrak{A}^*-\theta|>5r/2\},\end{align*}
where the final line follows from the fact that $M/K\lesssim \log(K/\delta)$ from our choice of $M$.
    
\end{proof}
\begin{lemma}\label{E-S-given-J-lemma-finite-moments} For fixed $v\geq 2$, let $X_{K+1}\sim P^{\otimes n}$ with $P\in \mathcal{P}_v$, defined by \eqref{finite-moments-class}, and suppose $\theta = \E X_{K+1,1}\in [-1,1]$. Assume that $n^2K\alpha^2\geq 1$ and suppose that we choose tuning parameters
\[T = A_1 (n^2K\alpha^2)^{1/(2v)}
    , \quad 
    r = A_2\left(\left[1\vee\sqrt{\log(K(n\wedge n^2\alpha^2))}\right]\left\{\frac{1}{\sqrt{n}}\vee T^{1-v}\right\}\vee \frac{\varepsilon^{-1/v}\wedge K^{1/v}}{n^{1-1/v}}\right),\]
and $M=A_3 r$, where $A_1$, $A_2$ and $A_3$ are sufficiently large constants depending only on $v$. Then, for any $\mathfrak{A}^*\in [-1,1+r/2]$, we have that
    \begin{equation*} |\E[S_1|\mathfrak{A}^*]| \lesssim  \frac{1}{nM}e^{-C'_1nM^2}+(nM)^{1-v} +T^{1-v}+\mathbf{1}\{|\mathfrak{A}^*-\theta|>5r/2\}.\end{equation*}
\end{lemma}
\begin{proof} Decomposing the expectation using the triangle inequality, we have that
    \begin{align}
        |\mathbb{E}[S_1|\mathfrak{A}^*]|&= \left|\mathbb{E}\left[\left.\mathfrak{A}^*-\theta+ \left[R_1\right]_{-M}^M\right|\mathfrak{A}^*\right]\right|\nonumber\\
        &= \left|\mathbb{E}\left[\left.\left[R_1\right]_{-M}^M-\left(\frac{1}{n}\sum_{i=1}^n[X_{K+1,i}]_{-T}^T-\mathfrak{A}^*\right)+\frac{1}{n}\sum_{i=1}^n[X_{K+1,i}]_{-T}^T-\theta\right| \mathfrak{A}^*\right]\right|\nonumber\\
        &\leq \mathbb{E}\left[\left.\left|\left[R_1\right]_{-M}^M-\left(\frac{1}{n}\sum_{i=1}^n[X_{K+1,i}]_{-T}^T-\mathfrak{A}^*\right)\right|\right| \mathfrak{A}^*\right]+\left|\E\left[\left.\frac{1}{n}\sum_{i=1}^n[X_{K+1,i}]_{-T}^T-\theta\right|\mathfrak{A}^*\right]\right|\nonumber\\
        &= \mathbb{E}\left[\left.\left(\left|\frac{1}{n}\sum_{i=1}^n[X_{K+1,i}]_{-T}^T-\mathfrak{A}^*\right|-M\right)_+\right|\mathfrak{A}^*\right]+\left|\E\left[\frac{1}{n}\sum_{i=1}^n[X_{K+1,i}]_{-T}^T-\theta\right]\right|\nonumber\\
        &\lesssim \E\Bigg[\Bigg(\Bigg|\frac{1}{n}\sum_{i=1}^n[X_{K+1,i}]_{-T}^T-\E\left\{[X_{K+1,i}]_{-T}^T\right\} \nonumber \\
        & \hspace{4cm} +\E\left\{[X_{K+1,i}]_{-T}^T\right\}-\theta+\theta-\mathfrak{A}^*\bigg|-M\bigg)_+\bigg|\mathfrak{A}^*\bigg] 
        +T^{1-v}\nonumber\\
        &\leq \E\left[\left.\left(\left|\frac{1}{n}\sum_{i=1}^n\left([X_{K+1,i}]_{-T}^T-\E\left\{[X_{K+1,i}]_{-T}^T\right\}\right)\right|-(M-(T/2)^{1-v}-|\mathfrak{A}^*-\theta|)\right)_+\right|\mathfrak{A}^*\right]\nonumber\\ \label{bound1}&\hspace{8cm}+T^{1-v},
    \end{align}
    where the final two inequalities follows from \eqref{finite-var-trunc-lemma-bias}.
    Applying the bound 
    \begin{equation}\label{bound2}|\E[S_1| \mathfrak{A}^*]|\leq |\mathfrak{A}^*-\theta|+M.\end{equation}
    and \eqref{bound1} to the following decomposition, we have, for some absolute constant $c\in(0,1)$, 
    \begin{align}&|\E[S_1|\mathfrak{A}^*] |\nonumber\\&= |\E[S_1|\mathfrak{A}^*]|\left\{\mathbf{1}\{|\mathfrak{A}^*-\theta|\leq 5r/2\}+\mathbf{1}\{|\mathfrak{A}^*-\theta|>5r/2\}(\mathbf{1}\{M>3\}+\mathbf{1}\{M\leq 3\})\right\}\nonumber\\
    &\stackrel{(a)}{\leq} \E\left[\left(\left|\frac{1}{n}\sum_{i=1}^n\left([X_{K+1,i}]_{-T}^T-\E\left\{[X_{K+1,i}]_{-T}^T\right\}\right)\right|-(M-(T/2)^{1-v}-5r/2)\right)_+\right]\mathbf{1}\{|\mathfrak{A}^*-\theta|\leq 5r/2\}\nonumber\\
    &\hspace{1cm}+\E\left[\left(\left|\frac{1}{n}\sum_{i=1}^n\left([X_{K+1,i}]_{-T}^T-\E\left\{[X_{K+1,i}]_{-T}^T\right\}\right)\right|-(M-(T/2)^{1-v}-2-r/2)\right)_+\right]\nonumber\\
    &\hspace{2cm}\cdot\mathbf{1}\{|\mathfrak{A}^*-\theta|>5r/2\}\mathbf{1}\{M>3\}+T^{1-v}+(2+r/2+3)\mathbf{1}\{|\mathfrak{A}^*-\theta|>5r/2\}\mathbf{1}\{M\leq 3\}\nonumber \\
    &\label{bound3}\stackrel{(b)}{\lesssim} \E\left[\left(\left|\frac{1}{n}\sum_{i=1}^n\left([X_{K+1,i}]_{-T}^T-\E\left\{[X_{K+1,i}]_{-T}^T\right\}\right)\right|-cM\right)_+\right]+T^{1-v}+\mathbf{1}\{|\mathfrak{A}^*-\theta|>5r/2\},\end{align}
    where (a) follows from \eqref{bound1}, \eqref{bound2}, and that $|\mathfrak{A}^*-\theta|\leq 2+r/2$; (b) follows from $M\gtrsim T^{1-v} \vee r$, provided that the constants $A_2$ and $A_3$ are sufficiently large. We now focus on bounding the expectation. 
    For some positive constants $c$, $C_1$, and $C_1'$ depending only on $v$, we have that
    \begin{align}
        &\E\left[\left(\left|\frac{1}{n}\sum_{i=1}^n\left([X_{K+1,i}]_{-T}^T-\E\left\{[X_{K+1,i}]_{-T}^T\right\}\right)\right|-cM\right)_+\right]\nonumber\\
        \leq & \int_0^\infty \mathbb{P}\left(\left|\frac{1}{n}\sum_{i=1}^n\left([X_{K+1,i}]_{-T}^T-\E\left\{[X_{K+1,i}]_{-T}^T\right\}\right)\right|\geq cM +x\right) \dd x\nonumber\\
        \stackrel{(a)}{\lesssim} & \int_{cM}^\infty e^{-C_1nx^2}+ n^{1-v}x^{-v} \dd x\nonumber\\
        \leq & \int_{cM}^\infty \frac{x}{cM}e^{-C_1nx^2}\dd x+ (cnM)^{1-v}\nonumber\\
        = & \frac{1}{2C_1 cnM}e^{-c^2C_1nM^2}+ (cnM)^{1-v}\nonumber\\
        \label{expectation-bound}\lesssim  &\frac{1}{nM}e^{-C'_1nM^2}+(nM)^{1-v},
    \end{align}
    where (a) follows from \eqref{finite-var-trunc-lemma-concentration}. Plugging \eqref{expectation-bound} into \eqref{bound3}, we conclude that
    \[|\E[S_1|\mathfrak{A}^*]| \lesssim  \frac{1}{nM}e^{-C'_1nM^2}+(nM)^{1-v} +T^{1-v}+\mathbf{1}\{|\mathfrak{A}^*-\theta|>5r/2\},\]
    where $C_1'$ is a positive constant depending only on $v$.
\end{proof}
\section{Proofs of results in \Cref{sec-mean-est-hetero}}\label{appendix-mean-est-hetero}
\subsection{Proof of lower bound}\label{appendix-hetero-LB-proof}

\begin{proof}[Proof of \Cref{hetero-LB}]
\textbf{Preliminary notation. }Recall the notation used in \Cref{sec:homogeneous-LB-proofs}. For each $k \in [K]$, let $M_{0,k}$, $M_{1,k}$, $G_{0,k}$ and $G_{1,k}$ denote distributions on $\mathbb{R}^{n_k}$, such that 
\begin{equation}\label{hetero-LB-proof-distributions}M_{0,k} = (1-\varepsilon)\mathcal{N}^{\otimes n_k}(\theta_0,1)+\varepsilon G_{0,k} \quad \mathrm{and} \quad M_{1,k} = (1-\varepsilon)\mathcal{N}^{\otimes n_k}(\theta_1,1)+\varepsilon G_{1,k}.\end{equation}
Let $Q_k M_{0,k}$ and $Q_k M_{1,k}$ denote the distribution of the privatised data $Z_k$ conditioned on $Z_{1:k-1}$, where $X_{1:K}$ is generated from $M_{0,1}\times \dots \times M_{0,K}$ and $M_{1,1}\times \dots \times M_{1,K}$ respectively. Recall that the privacy mechanisms $\{Q_k\}_{k=1}^K$ are such that the privatised data $Z_k$ satisfy $\alpha_k$--FDP for each $k\in[K]$. Specifically,
\[Q_kM_{0,k}(S |Z_{1:k-1}=z_{1:k-1}) =Q_kM_{0,k}(S |z_{1:k-1}) =\int Q_k(S|x_k, z_{1:k-1})\dd M_{0,k}(x_k),\]
noting that $X_k$ is independent of $Z_{1:k-1}$. For each $k\in [K]$, let $Q(M_{0,1}\times \dots\times M_{0,k})$ and $Q(M_{1,1}\times \dots\times M_{1,k})$ denote the distribution of $Z_{1:k}$ when $X_{1:k}$ is generated from $M_{0,1}\times \dots \times M_{0,k}$ and $M_{1,1}\times \dots \times M_{1,k}$ respectively. 

\medskip
\noindent \textbf{Step 1: Choice of the mean parameters $\theta_0$ and $\theta_1$. } 
Let $N = 2n_{\max}/n_{\min}$ and $A = 2\alpha_{\max}/\alpha_{\min}$. Choose~$\theta_0=0$ and
\begin{equation}\label{hetero-LB-proof-theta-choice}\theta_1^2 =1\wedge \frac{1}{\log_2(N)\log_2(A)}\min_{\substack{l_1,l_2\in \mathbb{N}:\\ \mathcal{I}_{l_1,l_2} \neq \emptyset}} \left\{\frac{2}{ 2^{l_1} n_{\min}|\mathcal{I}_{l_1,l_2}|}\vee \frac{\varepsilon^2\log_2(N)\log_2(A)}{2^{l_1} n_{\min}}\vee\frac{1/18}{ 2^{2(l_1+l_2)}n_{\min}^2\alpha_{\min}^2|\mathcal{I}_{l_1,l_2}|}\right\},\end{equation}
where $\mathcal{I}_{l_1,l_2}$ is as defined in \eqref{server-group-def}.
Then, for each $l_1, l_2\in \mathbb{N}$ such that $\mathcal{I}_{l_1,l_2} \neq \emptyset$, one of the following cases must hold: 
\begin{enumerate}
    \item \begin{equation}\label{eq:hetero-LB-case1} \theta_1^2 \leq \frac{\varepsilon^2}{2^{l_1} n_{\min}};\end{equation}
    \item \begin{equation}\label{eq:hetero-LB-case2}\theta_1^2\leq \frac{2}{ 2^{l_1} n_{\min}|\mathcal{I}_{l_1,l_2}|\log_2(N)\log_2(A)};\end{equation}
    \item \begin{equation}\label{eq:hetero-LB-case3}\theta_1^2\leq  \frac{1}{18\times 2^{2(l_1+l_2)}n_{\min}^2 \alpha_{\min}^2|\mathcal{I}_{l_1,l_2}|\log_2(N)\log_2(A)}.\end{equation}
\end{enumerate}
Using the above choices of $\theta_0$ and $\theta_1$ for the inlier distributions, we continue the proof by choosing the corresponding outlier distributions $\{G_{0,k}\}_{k=1}^K$ and $\{G_{1,k}\}_{k=1}^K$, from which we derive an upper bound on the separation between the data generating mechanisms.

\medskip
\noindent\textbf{Step 2: An upper bound on the KL divergence. }In this step, we show that there exist distributions $\{G_{0,k}\}_{k=1}^K$ and $\{G_{1,k}\}_{k=1}^K$ such that 
\begin{align}
    &\KL\left( Q\prod_{k=1}^K M_{0,k} ,Q \prod_{k=1}^K M_{1,k}\right)\nonumber\\
    &\stackrel{(a)}{=} \sum_{k=1}^K \int \KL(Q_kM_{0,k}(\cdot \mid z_{1:{k-1}}), Q_kM_{1,k}(\cdot \mid z_{1:{k-1}}))\dd Q(M_{0,1}\times \dots \times M_{0,k-1})(z_{1:k-1})\nonumber\\
    &= \sum_{l_1,l_2: \mathcal{I}_{l_1,l_2} \neq \emptyset}\sum_{k\in \mathcal{I}_{l_1,l_2}} \int \KL(Q_kM_{0,k}(\cdot \mid z_{1:{k-1}}), Q_kM_{1,k}(\cdot \mid z_{1:{k-1}}))\dd Q(M_{0,1}\times \dots \times M_{0,k-1})(z_{1:k-1})\nonumber\\
    &\stackrel{(b)}{\leq} \sum_{l_1=1}^{\lfloor\log_2(N)\rfloor}\sum_{l_2=1}^{\lfloor\log_2(A)\rfloor}\frac{1}{\log_2(N)\log_2(A)}\nonumber\\
    &\label{eq:KL-UB}\leq 1.
\end{align}
Note that (a) follows from the chain rule of KL divergences \citep[e.g.~Theorem 2.15 in][]{polyanskiy2025information}, and (b) holds provided that there exist distributions $\{G_{0,k}\}_{k=1}^K$ and $\{G_{1,k}\}_{k=1}^K$ such that 
\begin{multline}\label{hetero-LB-proof-KLbound}\sum_{k\in \mathcal{I}_{l_1,l_2}}\int\KL(Q_kM_{0,k}(\cdot \mid z_{1:{k-1}}), Q_kM_{1,k}(\cdot \mid z_{1:{k-1}}))\dd Q(M_{0,1}\times \dots \times M_{0,k-1})(z_{1:k-1})\\
\leq \frac{1}{\log_2(N)\log_2(A)}\end{multline}
holds for each $l_1 \in [\lfloor\log_2(N) \rfloor]$ and $l_2 \in [\lfloor \log_2(A)\rfloor]$.
We proceed to show that \eqref{hetero-LB-proof-KLbound} holds for each of the three cases in~\eqref{eq:hetero-LB-case1}, \eqref{eq:hetero-LB-case2} and~\eqref{eq:hetero-LB-case3}.

\textbf{Case 1. }Since $n_k < 2^{l_1}n_{\min}$ for each $k\in \mathcal{I}_{l_1,l_2}$, \eqref{eq:hetero-LB-case1} implies that 
\[\tau :=\TV(\mathcal{N}^{\otimes n_k}(\theta_0,1),\mathcal{N}^{\otimes n_k}(\theta_1,1))\leq|\theta_0-\theta_1|\sqrt{n_k}/2\leq \varepsilon\leq \varepsilon/(1-\varepsilon).\] 
Let $p_0$ and $p_1$ denote the densities of $\mathcal{N}^{\otimes n_k}(\theta_0,1)$ and $\mathcal{N}^{\otimes n_k}(\theta_1,1)$ respectively. Defining $G_{0,k}$ and $G_{1,k}$ by their densities 
\[g_0 = (1-c)p_0+ c(p_1-p_0)_+/\tau\quad \text{and} \quad g_1 = (1-c)p_1+c(p_0-p_1)_+/\tau \]
respectively, where $c = \tau/\{(1+\tau)\varepsilon\}$, we have that $M_{0,k} = M_{1,k}$ \citep[see e.g.~the Proof of Theorem 5.1 in][]{Chen2016a}. Thus, for each $k\in \mathcal{I}_{l_1,l_2}$, we have that
\begin{align*}
    \KL(Q_kM_{0,k}(\cdot \mid z_{1:{k-1}}), Q_kM_{1,k}(\cdot \mid z_{1:{k-1}}))&\leq \KL(M_{0,k}(\cdot \mid z_{1:{k-1}}), M_{1,k}(\cdot \mid z_{1:{k-1}}))\\
    &= \KL(M_{0,k}, M_{1,k})\\
    &=0
\end{align*}
where the first line follows from the data processing inequality \citep[e.g.~Theorem~7.4 in][]{polyanskiy2025information}, and the second line follows from the independence of the random variables $\{X_k\}_{k=1}^K$.

\textbf{Case 2. }We set $G_{0,k} = \mathcal{N}^{\otimes n_k}(\theta_0,1)$ and $G_{1,k}= \mathcal{N}^{\otimes n_k}(\theta_1,1)$, such that we have $M_{0,k}=\mathcal{N}^{\otimes n_k}(\theta_0,1)$ and $M_{1,k} = \mathcal{N}^{\otimes n_k}(\theta_1,1)$. With this choice for each $k\in \mathcal{I}_{l_1,l_2}$, we have that
\begin{align*}&\sum_{k\in \mathcal{I}_{l_1,l_2}} \int \KL(Q_kM_{0,k}(\cdot \mid z_{1:{k-1}}), Q_kM_{1,k}(\cdot \mid z_{1:{k-1}}))\dd Q(M_{0,1}\times \dots \times M_{0,k-1})(z_{1:k-1})\\
&\stackrel{(a)}{\leq} \sum_{k\in \mathcal{I}_{l_1,l_2}} \int \KL(M_{0,k},M_{1,k})\dd Q(M_{0,1}\times \dots \times M_{0,k-1})(z_{1:k-1})\\
&= \sum_{k\in \mathcal{I}_{l_1,l_2}}\KL(M_{0,k},M_{1,k})\\
&= \sum_{k\in \mathcal{I}_{l_1,l_2}}\KL(\mathcal{N}^{\otimes n_k}(\theta_0,1),\mathcal{N}^{\otimes n_k}(\theta_1,1))\\
&\leq \sum_{k\in \mathcal{I}_{l_1,l_2}}n_k|\theta_0-\theta_1|^2/2\\
&\stackrel{(b)}{\leq} 2^{l_1-1}n_{\min} |\mathcal{I}_{l_1,l_2}| |\theta_0-\theta_1|^2\\
&\stackrel{(c)}{\leq} \frac{1}{\log_2(N)\log_2(A)},
\end{align*}
where (a) follows from the data processing inequality \citep[e.g. Theorem 7.4 in][]{polyanskiy2025information}, (b) follows from the fact that $n_k<2^{l_1}n_{\min}$ for each $k\in \mathcal{I}_{l_1,l_2}$, and (c) follows from~\eqref{eq:hetero-LB-case2}.

\textbf{Case 3. }As in \textbf{Case 2}, we set $G_{0,k}$ and $G_{1,k}$ such that $M_{0,k}=\mathcal{N}^{\otimes n_k}(\theta_0,1)$ and $M_{1,k} = \mathcal{N}^{\otimes n_k}(\theta_1,1)$. Let $q_{0,k}$ and $q_{1,k}$ be the densities, with respect to a measure $\mu_k$, of $Q_kM_{0,k}(\cdot \mid z_{1:{k-1}})$ and $Q_kM_{1,k}(\cdot \mid z_{1:{k-1}})$ respectively. Using Lemma 6.1 in \cite{karwa2017finite}, we have that
\begin{equation}\label{eq:density-ratio-bound}
    e^{-\alpha_k'}\leq q_{1,k}/q_{0,k}\leq e^{\alpha'_k},
\end{equation}
where \begin{equation}\label{eq:alpha-prime-bound}\alpha'_k = 6\alpha_k n_k \TV(\mathcal{N}(\theta_0,1),\mathcal{N}(\theta_1,1)) \leq 3\alpha_k n_k |\theta_0-\theta_1|
\leq \frac{1}{\sqrt{2|\mathcal{I}_{l_1,l_2}|\log_2(N)\log_2(A)}}\leq 1,\end{equation}
by noting \eqref{eq:hetero-LB-case3} and that $n_k\leq 2^{l_1}n_{\min}$ and $\alpha_k\leq 2^{l_2}\alpha_{\min}$ in the penultimate inequality.
Thus, for any $z_{1:k-1}$, we have that
\begin{align*}&\sum_{k\in \mathcal{I}_{l_1,l_2}} \int \KL(Q_kM_{0,k}(\cdot \mid z_{1:{k-1}}), Q_kM_{1,k}(\cdot \mid z_{1:{k-1}}))\dd Q(M_{0,1}\times \dots \times M_{0,k-1})(z_{1:k-1})\\
&\leq \sum_{k\in \mathcal{I}_{l_1,l_2}} \int \left[\int \log\left(\frac{q_{0,k}}{q_{1,k}}\right)q_{0,k} \dd \mu_k+ \int \log\left(\frac{q_{1,k}}{q_{0,k}}\right)q_{1,k} \dd \mu_k\right]\dd Q(M_{0,1}\times \dots \times M_{0,k-1})(z_{1:k-1})\\
&\leq \sum_{k\in \mathcal{I}_{l_1,l_2}} \int \left[\int \log\left(\frac{q_{1,k}}{q_{0,k}}\right)\left(\frac{q_{1,k}}{q_{0,k}}-1\right)q_{0,k}\dd \mu_k\right]\dd Q(M_{0,1}\times \dots \times M_{0,k-1})(z_{1:k-1})\\
&\stackrel{(a)}{\leq} \sum_{k\in \mathcal{I}_{l_1,l_2}} \int \left[\int \{\alpha'_k (e^{\alpha '_k}-1)\vee (-\alpha_k')(e^{-\alpha_k'}-1)\}q_{0,k} \dd \mu_k\right]\dd Q(M_{0,1}\times \dots \times M_{0,k-1})(z_{1:k-1})\\
&\leq \sum_{k\in \mathcal{I}_{l_1,l_2}} \int \left[\alpha'_k (e^{\alpha '_k}-1)\right]\dd Q(M_{0,1}\times \dots \times M_{0,k-1})(z_{1:k-1})\\
&\leq \sum_{k\in \mathcal{I}_{l_1,l_2}}\alpha'_k (e^{\alpha'_k}-1)\\
&\stackrel{(b)}{\leq} \sum_{k\in \mathcal{I}_{l_1,l_2}}2\alpha'^2_k\\
&\stackrel{(c)}{\leq} \sum_{k\in \mathcal{I}_{l_1,l_2}}\frac{1}{|\mathcal{I}_{l_1,l_2}|\log_2(N)\log_2(A)}\\
&\leq \frac{1}{\log_2(N)\log_2(A)}
\end{align*}
where (a) follows from \eqref{eq:density-ratio-bound}, (b) follows from noting that $e^x - 1 \leq 2x$ for $x\in [0,1]$, and (c) follows from \eqref{eq:alpha-prime-bound}.  

\medskip
\textbf{Step 3: Applying Le Cam's two-point lemma. }Using Le Cam's two-point lemma \citep[e.g.~Lemma 1 in][]{Yu1997assouad}, we have that
    \begin{align*}
        \mathcal{R}&\geq \frac{(\theta_0-\theta_1)^2}{8}\left(1- \TV\left(Q\prod_{k=1}^K M_{0,k} , Q\prod_{k=1}^K M_{1,k}\right)\right)\\
        &\stackrel{(a)}{\geq}\frac{(\theta_0-\theta_1)^2}{8}\left(1- \sqrt{\frac{1}{2}\KL\left( Q\prod_{k=1}^K M_{0,k} ,Q \prod_{k=1}^K M_{1,k}\right)}\right)\\
        &\stackrel{(b)}{\geq} \frac{1-1/\sqrt{2}}{8}(\theta_0-\theta_1)^2\\
        &\stackrel{(c)}{\gtrsim} 1\wedge \frac{1}{\log_2(N)\log_2(A)}\min_{\substack{l_1,l_2\in \mathbb{N}:\\ \mathcal{I}_{l_1,l_2} \neq \emptyset}} \left\{\frac{2}{ 2^{l_1} n_{\min}|\mathcal{I}_{l_1,l_2}|}\vee \frac{\varepsilon^2\log_2(N)\log_2(A)}{2^{l_1} n_{\min}}\vee\frac{1/18}{ 2^{2(l_1+l_2)}n_{\min}^2\alpha_{\min}^2|\mathcal{I}_{l_1,l_2}|}\right\}\\
        &\gtrsim 1\wedge \frac{1}{\log_2(N)\log_2(A)}\min_{\substack{l_1,l_2\in \mathbb{N}:\\ \mathcal{I}_{l_1,l_2} \neq \emptyset}} \left\{\frac{1}{ 2^{l_1} n_{\min}|\mathcal{I}_{l_1,l_2}|}\vee \frac{\varepsilon^2}{2^{l_1} n_{\min}}\vee\frac{1}{ 2^{2(l_1+l_2)}n_{\min}^2\alpha_{\min}^2|\mathcal{I}_{l_1,l_2}|}\right\},
    \end{align*}
    where (a) follows from Pinsker's inequality \citep[e.g.~Lemma 2.5 in][]{Tsybakov2009introduction}, (b) follows from \eqref{eq:KL-UB} and (c) follows from \eqref{hetero-LB-proof-theta-choice}.
\end{proof}
\subsection{Proof of upper bound}\label{appendix-hetero-UB-proof}
Before proving \Cref{hetero-UB-thm} and the required lemmas, we extend the notation defined in \Cref{appendix:homogeneous-UB-proofs} to the heterogeneous case.

In the heterogeneous setting, we require control of the outputs $\mathfrak{A}^*_{l_1,l_2}$ and $\hat{\theta}_{l_1,l_2}$ of the estimation procedure when it is given, as input, one particular group of batches $\mathcal{I}_{l_1,l_2}$ only. To be specific, given any fixed $(l_1,l_2)$, we randomly split $\mathcal{I}_{l_1,l_2}$ into two halves $\mathcal{B}$ and $\mathcal{B}'$ such that we have $|\mathcal{B}|=|\mathcal{B}'|=|\mathcal{I}_{l_1,l_2}|/2$. The first half $\mathcal{B}$ is used as input data for \Cref{alg:mean-localisation} in the construction of the initial estimator \begin{equation}\label{eq:A_l1l2-def}\mathfrak{A}_{l_1,l_2}^* \leftarrow \text{\sc MeanLocalisation}(\{\{X_{k,i}\}_{i=1}^{n_k}\}_{k\in \mathcal{B}}, \{\alpha_k\}_{k\in \mathcal{B}}, r_{l_1,l_2}, T_{l_1,l_2}),\end{equation}
and combined with the second half $\mathcal{B}'$, we have from \Cref{alg:mean-refinement} the estimator
\begin{equation}\label{theta-hat_l1l2-def}
    \hat{\theta}_{l_1,l_2} \leftarrow \text{\sc MeanRefinement}(\{\{X_{k,i}\}_{i=1}^{n_k}\}_{k\in \mathcal{B}'}, \mathfrak{A}^*_{l_1,l_2}, \{\alpha_k\}_{k\in \mathcal{B}'}, M_{l_1,l_2}, T_{l_1,l_2}).
\end{equation}
Regarding \Cref{alg:mean-localisation}, for each $k\in  \mathcal{I}_{l_1,l_2}$ and $j\in [\lceil 2r_{l_1,l_2}^{-1}\rceil]$, let the location of the $j$th bin be defined by \[a_j = -1+(j-1/2)r_{l_1,l_2},\] and denote the smoothed votes by 
\begin{equation}\label{Z_kj-def-hetero}Z_{k,j} = \left(2r_{l_1,l_2} - \left|\frac{1}{n_k}\sum_{i=1}^{n_k} [X_{k,i}]_{-T_{l_1,l_2}}^{T_{l_1,l_2}} - a_j\right|\right)_+ + \frac{16 T_{l_1,l_2}}{\alpha_k n_k} W_{k,j},\end{equation}
    where $\{W_{k,j}\}_{k\in\mathcal{I}_{l_1,l_2},j\in [\lceil 2r_{l_1,l_2}^{-1}\rceil]}$ are i.i.d.~standard Laplace random variables. The initial estimator $\mathfrak{A}^*_{l_1,l_2}$ is defined by
    \[\mathfrak{A}^*_{l_1,l_2}=a_{j^*},\]
    where \[j^* = \min\argmax_{j\in [\lceil 2r_{l_1,l_2}^{-1}\rceil]}\frac{1}{|\mathcal{B}|}\sum_{k\in \mathcal{B}} Z_{k,j}.\]
    Regarding \Cref{alg:mean-refinement}, for each $k\in \mathcal{B}'$, let 
    \begin{equation} \label{eq-def-Z_k-hetero}
        Z_k = \mathfrak{A}^*_{l_1,l_2}+ \left[\frac{1}{n_k}\sum_{i = 1}^{n_k} [X_{k,i}]_{-T_{l_1,l_2}}^{T_{l_1,l_2}} -\mathfrak{A}^*_{l_1,l_2}\right]_{-M_{l_1,l_2}}^{M_{l_1,l_2}} + \frac{2T_{l_1,l_2}}{\alpha_k n_k}L_k ,
    \end{equation}
    where $\{L_k\}_{k\in \mathcal{B}'}$ are i.i.d.~standard Laplace random variables. The estimator $\hat{\theta}_{l_1,l_2}$ can be written as
    \begin{equation}\label{theta-hat_l1l2-def2}\hat{\theta}_{l_1,l_2} = \frac{1}{|\mathcal{B}'|}\sum_{k\in \mathcal{B}'}Z_k.\end{equation}
Recall the tuning parameters $\{T_{l_1,l_2}, r_{l_1,l_2},M_{l_1,l_2}\}_{l_1,l_2 \in \mathbb{N}: \mathcal{I}_{l_1,l_2}\neq \emptyset}$ used in \Cref{alg:hetero-estimation} are defined by
    \begin{equation*}
        T_{l_1,l_2}= C_1\left( 1 \vee \sqrt{\log(|\mathcal{I}_{l_1,l_2}|(2^{l_1}n_{\min}\wedge2^{2(l_1+l_2)}n_{\min}^2\alpha_{\min}^2)}\right),
    \end{equation*}    \begin{equation*}
        r_{l_1,l_2} = C_2 T_{l_1,l_2} \left(\frac{1}{2^{l_1} n_{\min}}\vee \frac{T_{l_1,l_2}^2+\log(4/\delta)}{2^{2(l_1+l_2)}n_{\min}^2\alpha_{\min}^2|\mathcal{I}_{l_1,l_2}|} \right)^{1/2} , \quad M_{l_1,l_2}=3r_{l_1,l_2},
    \end{equation*}
    where $C_1,C_2>0$ are sufficiently large absolute constants.
The function
\begin{equation*}
    L(\varepsilon, l_1,l_2, \delta) = C'\left\{\varepsilon M_{l_1,l_2}\vee\sqrt{\frac{\log(4/\delta)}{|\mathcal{I}_{l_1,l_2}|}}\left(M_{l_1,l_2}\vee\frac{T_{l_1,l_2}}{2^{l_1+l_2}n_{\min}\alpha_{\min}}\right)\right\},
\end{equation*}
where $C'>0$ is a sufficiently large absolute constant, will be crucial in the following proofs. 

\begin{proof}[Proof of \Cref{non-adaptive-thm}]
    Note that, the condition \eqref{non-adaptive-thm-condition} in \Cref{thm-hetero-rate} implies that the condition \eqref{non-adaptive-prop-condition} in \Cref{non-adaptive-prop} holds. Hence, the high probability guarantee \eqref{non-adaptive-thm-high-prob-guarantee} follows immediately from \Cref{non-adaptive-prop} with $(\tilde{l}_1,\tilde{l}_2)=(l_1^*,l_2^*)$.

To complete the proof of \eqref{non-adaptive-risk-bound}, we first note that 
\begin{align*}
    \delta \asymp  
    R(\varepsilon,l_1^*,l_2^*)
    \geq \frac{1}{K_{l_1^*,l_2^*}(\tilde{n}_{l_1^*}\wedge \tilde{n}_{l_1^*}^2\tilde{\alpha}_{l_2^*}^2)}
    \stackrel{}{>}  e^{-C^2 K_{l_1^*,l_2^*}},
\end{align*}
where the final inequality follows from \eqref{non-adaptive-thm-condition}. Hence, the choice of $\delta$ is valid. Finally, recall the assumption that $\theta\in [-1,1]$ and, using the law of total expectation, we have that
\begin{align*}
    \E \{|[\hat{\theta}]_{-1}^1-\theta|^2\}&= \E\left[|[\hat{\theta}]_{-1}^1-\theta|^2\big| |\hat{\theta}-\theta|\leq L(\varepsilon,l_1^*,l_2^*,\delta)\right]\mathbb{P}\left(|\hat{\theta}-\theta|\leq L(\varepsilon,l_1^*,l_2^*,\delta)\right)\nonumber\\
    &\hspace{2cm}+\E\left[|[\hat{\theta}]_{-1}^1-\theta|^2\big| |\hat{\theta}-\theta|>L(\varepsilon,l_1^*,l_2^*,\delta)\right]\mathbb{P}\left(|\hat{\theta}-\theta|> L(\varepsilon,l_1^*,l_2^*,\delta)\right)\nonumber\\
    &\leq \{L(\varepsilon,l_1^*,l_2^*,\delta)\}^2 + 8\delta\\
    &\lesssim \log^2\left(\sum_{k=1}^K n_k\right)\min_{\substack{l_1,l_2\in \mathbb{N}: \\\mathcal{I}_{l_1,l_2} \neq \emptyset}} \left\{R(\varepsilon,l_1,l_2)\right\},
\end{align*}
    where the penultimate inequality follows from \eqref{non-adaptive-thm-condition}, and the final inequality follows from the observation made in \eqref{approx-L-def}.
\end{proof}

\begin{proof}[Proof of \Cref{hetero-UB-thm}] The proof proceeds in three steps. In \textbf{Step 1}, we show that, for each iteration $t$, the estimator $\hat{\theta}_{l_1^t,l_2^t}$ constructed using the group of batches labelled by the pair $(l_1^t,l_2^t)$ lies within distance $L(\varepsilon,l_1^t,l_2^t,\delta)$ of $\theta$ with high probability. In \textbf{Step 2}, we show that the adaptive estimator $\hat{\theta}$ concentrates within an interval of near-optimal length of the order $L(\varepsilon,l_1^*,l_2^*,\delta)$ around $\theta$ with high probability. In \textbf{Step 3}, we show an upper bound on the minimax risk as a corollary of the high probability bound obtained in \textbf{Step 2}.

\medskip
\noindent\textbf{Step 1.} Note that the conditions \eqref{adaptive-UB-condition} and \eqref{adaptive-UB-condition-delta} imply that the conditions \eqref{non-adaptive-prop-condition} and \eqref{non-adaptive-prop-condition-delta} hold respectively for any input batch group label $(\tilde{l}_1,\tilde{l}_2) = (l_1^t,l_2^t)$ to \Cref{alg:hetero-estimation}. Thus, applying \Cref{non-adaptive-prop}, we have that
\begin{equation}\label{hetero-group-estimator-high-prob-guarantee}
    \mathbb{P}\left(|\hat{\theta}_{l_1^t,l_2^t}-\theta|> L(\varepsilon,l_1^t,l_2^t,\delta)\right)\leq 2\delta
\end{equation}
for any fixed $t\in [\lceil\log_2(K)/2\rceil]$.

\medskip
\noindent\textbf{Step 2. }In this step, we prove \eqref{hetero-UB-thm-highprob-bound} by showing that the final estimator $\hat{\theta}$ satisfies
\begin{equation}\label{hetero-final-estimator-high-prob-guarantee}\mathbb{P}\left(|\hat{\theta}-\theta|> 6L(\varepsilon,l_1^*,l_2^*,\delta)\right)\leq 2\delta\log_2(K),\end{equation}
where $l_1^*,l_2^*\in \mathbb{N}$ satisfy
\[L(\varepsilon, l_1^*,l_2^*,\delta) =\min_{l_1,l_2\in \mathbb{N}: \mathcal{I}_{l_1,l_2} \neq \emptyset} L(\varepsilon,l_1,l_2,\delta).\]
Note that since $|\mathcal{I}_{l_1,l_2}|\leq K$, having a contamination level $\varepsilon\leq K^{-1/2}$ implies that
\[L(\varepsilon,l_1,l_2,\delta) = C'\left\{\sqrt{\frac{\log(4/\delta)}{|\mathcal{I}_{l_1,l_2}|}}\left(M_{l_1,l_2}\vee\frac{T_{l_1,l_2}}{2^{l_1+l_2}n_{\min}\alpha_{\min}}\right)\right\}= L(K^{-1/2},l_1,l_2,\delta)\]
for any fixed $l_1$, $l_2$, and $\delta$. Hence, we need only consider the case where $\varepsilon\geq K^{-1/2}$ and iterate through $t = 1,\ldots, \lceil\log_2(K)/2\rceil$ only. For the rest of \textbf{Step 2}, we assume that $\varepsilon\geq K^{-1/2}$.
Define $t^*= \lfloor\log_2(1/\varepsilon)\rfloor$, such that the true contamination level $\varepsilon\in (2^{-(t^*+1)},2^{-t^*}]$. Noting that $t^* = \lfloor \log_2(1/\varepsilon)\rfloor\leq \lceil\log_2(K)/2\rceil$, observe that we have the event
\begin{equation}\label{hetero-good-event-def}\Omega := \bigcap_{t=0}^{t^*}\{\theta \in I_t\} \subseteq \{\hat{t}\geq t^*\}\end{equation}
by the definitions of $I_t$ and $\hat{t}$ in \Cref{alg:adaptive-mean-estimation}. Recall that $\hat{\theta}=[\hat{m}]_{-1}^1$, where $\hat{m}$ is defined to be the midpoint of the interval $I_{\hat{t}}$. Since $\bigcap_{t=0}^{\hat{t}}I_t \neq\emptyset$ by the definition of $\hat{t}$ in \Cref{alg:adaptive-mean-estimation}, and that $\hat{t}\geq t^*$ on the event $\Omega$, there exists a point $x\in \bigcap_{t=0}^{\hat{t}}I_t\subseteq I_{\hat{t}}\cap I_{t^*}$. Using the triangle inequality, note that, on the event $\Omega$, we have 
\begin{align}
    |\hat{\theta}-\theta|&\leq |\hat{m}-\theta| \\
    &\leq |\hat{m}-x|+|x- \hat{\theta}_{l_1^{t^*},l_2^{t^*}}| + |\hat{\theta}_{l_1^{t^*},l_2^{t^*}}-\theta|\nonumber\\
&\stackrel{(a)}{\leq} L(\varepsilon_{\hat{t}},l_1^{\hat{t}},l_2^{\hat{t}},\delta)+L(\varepsilon_{t^*},l_1^{t^*},l_2^{t^*},\delta)+L(\varepsilon_{t^*},l_1^{t^*},l_2^{t^*},\delta)\nonumber\\
&\stackrel{(b)}{\leq} 3L(\varepsilon_{t^*},l_1^{t^*},l_2^{t^*},\delta)\nonumber\\
&\stackrel{(c)}{\leq} 3L(\varepsilon_{t^*},l_1^*,l_2^*,\delta)\nonumber\\
&\stackrel{(d)}{\leq} 6 L(\varepsilon,l_1^*,l_2^*,\delta)\nonumber,
\end{align}
where (a) follows from the definition \eqref{eq:interval-def} that $I_t$ is an interval of radius $L(\varepsilon_t,l_1^t,l_2^t,\delta)$ centred at $\hat{\theta}_{l_1^t,l_2^t}$, (b) follows from the facts that $\hat{t} \geq t^*$ and $L(\varepsilon,l_1,l_2,\delta)$ is increasing in $\varepsilon$ for any $l_1,l_2\in \mathbb{N}$ and $\delta>0$, as defined in \eqref{hetero-L-def}, (c) follows from the definition of $l_1^t,l_2^t$ in \eqref{adaptive-alg-minimiser-condition}, and (d) follows from the definition \eqref{hetero-L-def} of $L(\varepsilon,l_1,l_2,\delta)$ and the fact that $\varepsilon_{t^*}\leq 2\varepsilon$. We complete the proof of \eqref{hetero-final-estimator-high-prob-guarantee} by noting that the event $\Omega$ defined in \eqref{hetero-good-event-def} satisfies
\begin{align*}
    \mathbb{P}(\Omega^c)&=   \mathbb{P}\left(\bigcup_{t= 1}^{t^*}\{\theta \notin I_t\}\right)\\
    &\stackrel{(a)}{=} \mathbb{P}\left(\bigcup_{t=1}^{t^*} \{|\hat{\theta}_{l_1^t,l_2^t}-\theta| > L(\varepsilon_t,l_1^t,l_2^t,\delta)\}\right)\\
    &\stackrel{(b)}{\leq} \mathbb{P}\left(\bigcup_{t=1}^{t^*} \{|\hat{\theta}_{l_1^t,l_2^t}-\theta| > L(\varepsilon,l_1^t,l_2^t,\delta)\}\right)\\
    &\stackrel{(c)}{\leq}\sum_{t=1}^{t^*}\mathbb{P}\left(|\hat{\theta}_{l_1^t,l_2^t}-\theta|> L(\varepsilon,l_1^t,l_2^t,\delta)\right)\\
    &\stackrel{(d)}{\leq}  2\delta \lceil \log_2(K)/2\rceil,
\end{align*}
where (a) follows from the definition of $I_t$ in \eqref{eq:interval-def}, (b) follows from the fact that $\varepsilon\leq \varepsilon_{t^*}\leq \varepsilon_t$ for any $t\in \{1,\dots,t^*\}$ and that $L(\varepsilon,l_1,l_2,\delta)$ is increasing in $\varepsilon$ for any $l_1,l_2\in \mathbb{N}$ and $\delta>0$ as defined in \eqref{hetero-L-def}, (c) follows from the union bound, and (d) follows from \eqref{hetero-group-estimator-high-prob-guarantee} and that $t^*\leq \lceil \log_2(K)/2\rceil$. 

\medskip
\noindent\textbf{Step 3.} Using \eqref{hetero-final-estimator-high-prob-guarantee} and the assumption that $\theta\in [-1,1]$, the law of total expectation implies that
\begin{align}
    \E [|\hat{\theta}-\theta|^2]&= \E\left[|\hat{\theta}-\theta|^2\big| |\hat{\theta}-\theta|\leq6L(\varepsilon,l_1^*,l_2^*,\delta)\right]\mathbb{P}\left(|\hat{\theta}-\theta|\leq 6L(\varepsilon,l_1^*,l_2^*,\delta)\right)\nonumber\\
    &\hspace{2cm}+\E\left[|\hat{\theta}-\theta|^2\big| |\hat{\theta}-\theta|>6L(\varepsilon,l_1^*,l_2^*,\delta)\right]\mathbb{P}\left(|\hat{\theta}-\theta|> 6L(\varepsilon,l_1^*,l_2^*,\delta)\right)\nonumber\\
    &\label{eq:risk-bound-hetero-UB-proof}\leq 36 \{L(\varepsilon,l_1^*,l_2^*,\delta)\}^2+8\delta\log_2(K).
\end{align}

Noting that the choice of $\delta$ in \eqref{eq:delta-choice-hetero-UB-proof} implies that 
\begin{align}\label{eq:delta-bound-hetero-ub-proof}
    \log(4/\delta) \leq \log\left(\max_{l_1,l_2\in \mathbb{N}:\mathcal{I}_{l_1,l_2}\neq \emptyset}\{2^{l_1}n_{\min}|\mathcal{I}_{l_1,l_2}|\}\right)\lesssim \log\left(\sum_{k=1}^K n_k\right).
\end{align}
Further, we have that, for each $(l_1,l_2)$ such that $\mathcal{I}_{l_1,l_2}\neq \emptyset$,
\begin{align}\label{eq:T-bound-hetero-ub-proof}
    T_{l_1,l_2}^2 \lesssim \log(|\mathcal{I}_{l_1,l_2}|2^{l_1}n_{\min})\lesssim \log\left(\sum_{k=1}^K n_k\right),
\end{align}
and that
\begin{align}\label{eq:L-simplified}
    L(\varepsilon,l_1,l_2,\delta)\lesssim T_{l_1,l_2}\left\{\frac{\varepsilon}{\sqrt{2^{l_1}n_{\min}}}\vee \frac{\varepsilon T_{l_1,l_2}\vee \sqrt{\log(4/\delta)}}{2^{l_1+l_2}n_{\min}\alpha_{\min}\sqrt{|\mathcal{I}_{l_1,l_2}|}}\vee\frac{\sqrt{\log(4/\delta)}}{\sqrt{2^{l_1}n_{\min}|\mathcal{I}_{l_1,l_2}|}}\right\}
\end{align}
by the definition of $L(\varepsilon,l_1,l_2,\delta)$ in \eqref{hetero-L-def}.
Therefore, plugging in \eqref{eq:delta-choice-hetero-UB-proof}, \eqref{eq:delta-bound-hetero-ub-proof}, \eqref{eq:T-bound-hetero-ub-proof} and \eqref{eq:L-simplified} into \eqref{eq:risk-bound-hetero-UB-proof}, we have
\begin{align}
    \E[|\hat{\theta}-\theta|^2]&\lesssim \log^2\left(\sum_{k=1}^K n_k\right)\min_{l_1,l_2\in \mathbb{N}:\mathcal{I}_{l_1,l_2}\neq \emptyset}\left\{\frac{\varepsilon^2}{2^{l_1}n_{\min}}\vee\frac{1}{2^{2(l_1+l_2)}n_{\min}^2\alpha_{\min}^2 |\mathcal{I}_{l_1,l_2}|}\vee \frac{1}{2^{l_1}n_{\min}|\mathcal{I}_{l_1,l_2}|}\right\}\nonumber\\
    &\hspace{4cm} \vee\log_2(K)\max_{t\in[\lceil\log_2(K)/2\rceil] } \exp\{-C^2 K_{l_1^t,l_2^t}\}\nonumber\\
    &\lesssim \log^2\left(\sum_{k=1}^K n_k\right)\nonumber\\
    \label{eq:final-risk-bound}&\left[\min_{\substack{l_1,l_2\in \mathbb{N}:\\\mathcal{I}_{l_1,l_2}\neq \emptyset}}\left\{\frac{\varepsilon^2}{2^{l_1}n_{\min}}\vee\frac{1}{2^{2(l_1+l_2)}n_{\min}^2\alpha_{\min}^2 |\mathcal{I}_{l_1,l_2}|}\vee \frac{1}{2^{l_1}n_{\min}|\mathcal{I}_{l_1,l_2}|}\right\} \vee\max_{t\in[\lceil\log_2(K)/2\rceil] } e^{-C^2 K_{l_1^t,l_2^t}}\right].
\end{align}
\end{proof}
\subsubsection{Proofs of auxiliary results}
The following proposition provides the theoretical guarantees of \Cref{alg:hetero-estimation}.
\begin{proposition}\label{non-adaptive-prop}
    Consider the contamination model with heterogeneous server sizes detailed in~\eqref{Gaussian-distribution} and~\eqref{heterogeneous-model}. 
    There exists an absolute constant $C>0$ such that when the condition
    \begin{equation}\label{non-adaptive-prop-condition}
         \varepsilon+\sqrt{\frac{\log(2^{\tilde{l}_1}n_{\min}\wedge 2^{2(\tilde{l}_1+\tilde{l}_2)}n_{\min}^2\alpha_{\min}^2|\mathcal{I}_{\tilde{l}_1,\tilde{l}_2}|)}{|\mathcal{I}_{\tilde{l}_1,\tilde{l}_2}|}}+ \frac{1}{\sqrt{(2^{\tilde{l}_1}n_{\min}\wedge2^{2(\tilde{l}_1+\tilde{l}_2)}n_{\min}^2\alpha_{\min}^2)|\mathcal{I}_{\tilde{l}_1,\tilde{l}_2}|}}<C\end{equation}
    holds for the input batch group label $(\tilde{l}_1,\tilde{l}_2)$ with 
    \begin{equation}\label{non-adaptive-prop-condition-delta}
        \delta \geq e^{-C^2|\mathcal{I}_{\tilde{l}_1,\tilde{l}_2}|},
    \end{equation}
    the output estimator $\hat{\theta}$ of \Cref{alg:hetero-estimation} satisfies
    \begin{equation*}
        \mathbb{P}\left(|\hat{\theta}-\theta|> L(\varepsilon,\tilde{l}_1,\tilde{l}_2,\delta)\right)\leq 2\delta,
    \end{equation*}
    where the function $L(\cdot,\cdot,\cdot,\cdot)$ is defined in \eqref{hetero-L-def}.
\end{proposition}

\begin{proof}
    \noindent\textbf{Step 1.} In this step, we show that, given $\delta \geq e^{-C^2|\mathcal{I}_{\tilde{l}_1,\tilde{l}_2}|}$, the initial estimator $\mathfrak{A}^*_{\tilde{l}_1,\tilde{l}_2}$ defined in \Cref{alg:hetero-estimation} satisfies
\begin{equation}\label{eq:hetero-high-prob-bound}\mathbb{P}(|\mathfrak{A}^*_{\tilde{l}_1,\tilde{l}_2}-\theta|\leq 5r_{\tilde{l}_1,\tilde{l}_2}/2)\geq 1-\delta,\end{equation}
provided that the condition \eqref{non-adaptive-prop-condition} holds. To prove this, we apply \Cref{heterogeneity-step1-lemma}.
Observe that, due to the assumption $\theta\in [-1,1]$, the bound \eqref{eq:hetero-high-prob-bound} holds trivially for any $\delta\geq 0$ whenever $r_{\tilde{l}_1,\tilde{l}_2}\geq 1$, without needing to invoke \Cref{heterogeneity-step1-lemma}. Thus, we proceed by verifying the condition \eqref{step1-heterogeneous-delta-condition} stated in \Cref{heterogeneity-step1-lemma} with $l_1=\tilde{l}_1$ and $l_2=\tilde{l}_2$ assuming that $r_{\tilde{l}_1,\tilde{l}_2}\leq 1$. Substituting the definition of $r_{l_1,l_2}$ in \eqref{eq:hetero-UB-thm-r-def}, we have that
\begin{align}
    &\varepsilon+\frac{1}{\sqrt{(2^{\tilde{l}_1}n_{\min}\wedge2^{2(\tilde{l}_1+\tilde{l}_2)}n_{\min}^2\alpha_{\min}^2)|\mathcal{I}_{\tilde{l}_1,\tilde{l}_2}|}}+\sqrt{\frac{\log(12/(r_{\tilde{l}_1,\tilde{l}_2}\delta))}{|\mathcal{I}_{\tilde{l}_1,\tilde{l}_2}|}}\nonumber\\
    \stackrel{(a)}{\leq} &\varepsilon+\frac{1}{\sqrt{(2^{\tilde{l}_1}n_{\min}\wedge2^{2(\tilde{l}_1+\tilde{l}_2)}n_{\min}^2\alpha_{\min}^2)|\mathcal{I}_{\tilde{l}_1,\tilde{l}_2}|}}
    +\sqrt{\frac{\log(2^{\tilde{l}_1}n_{\min}\wedge 2^{(\tilde{l}_1+\tilde{l}_2)}n_{\min}^2\alpha_{\min}^2|\mathcal{I}_{\tilde{l}_1,\tilde{l}_2}|)+2\log(1/\delta)}{2|\mathcal{I}_{\tilde{l}_1,\tilde{l}_2}|}}\nonumber\\
    \stackrel{(b)}{<} &C+\sqrt{\frac{\log(1/\delta)}{|\mathcal{I}_{\tilde{l}_1,\tilde{l}_2}|}}\nonumber\\
    \stackrel{(c)}{<} &2C\nonumber,
\end{align}
as required by \eqref{step1-heterogeneous-delta-condition}, provided that the constant $C$ in condition \eqref{non-adaptive-prop-condition} is sufficiently small. Note that we have used the fact that $T_{\tilde{l}_1,\tilde{l}_2}\geq 1$, $\log(12/(r_{\tilde{l}_1,\tilde{l}_2}\delta))\geq 1$ and that the constant $C_2$ in the definition \eqref{eq:hetero-UB-thm-r-def} of $r_{l_1,l_2}$ is sufficiently large in (a), (b) follows from \eqref{non-adaptive-prop-condition}, and (c) follows from \eqref{non-adaptive-prop-condition-delta}. Thus, we conclude that \eqref{eq:hetero-high-prob-bound} holds.

\medskip
\noindent\textbf{Step 2. }In this step, we show that, given \[\delta \geq e^{-C^2|\mathcal{I}_{\tilde{l}_1,\tilde{l}_2}|},\] 
the estimator $\hat{\theta}$ satisfies
\begin{equation*}\mathbb{P}\left(|\hat{\theta} - \theta|> L(\varepsilon, \tilde{l}_1,\tilde{l}_2, \delta)\right) \leq  2\delta,
\end{equation*} 
where $L(\varepsilon,l_1,l_2,\delta)$ is defined in \eqref{hetero-L-def}. We first note that $\hat{\theta} = \hat{\theta}_{\tilde{l}_1,\tilde{l}_2}$ defined in \eqref{theta-hat_l1l2-def}. Further, since $|\mathcal{I}_{\tilde{l}_1,\tilde{l}_2}|>1/C^2$ due to condition \eqref{non-adaptive-prop-condition}, we have that
\[e^{-C^2|\mathcal{I}_{\tilde{l}_1,\tilde{l}_2}|}\geq 4e^{-|\mathcal{I}_{\tilde{l}_1,\tilde{l}_2}|/4}\]
for a sufficiently small constant $C>0$. Thus, applying \Cref{hetero-refinement-high-prob-bound} and using the law of total probability, we have that
\begin{align}
    &\mathbb{P}\left(|\hat{\theta}_{\tilde{l}_1,\tilde{l}_2} - \theta|> L(\varepsilon, \tilde{l}_1,\tilde{l}_2, \delta)\right)\nonumber\\ 
    \leq& \mathbb{P}\left(|\hat{\theta}_{\tilde{l}_1,\tilde{l}_2} - \theta|> L(\varepsilon, \tilde{l}_1,\tilde{l}_2, \delta)\big||\mathfrak{A}^*_{\tilde{l}_1,\tilde{l}_2}-\theta|\leq 5r_{\tilde{l}_1,\tilde{l}_2}/2\right)\mathbb{P}(|\mathfrak{A}^*_{\tilde{l}_1,\tilde{l}_2}-\theta|\leq 5r_{\tilde{l}_1,\tilde{l}_2}/2)\nonumber\\
    &\hspace{2cm}+\mathbb{P}(|\mathfrak{A}^*_{\tilde{l}_1,\tilde{l}_2}-\theta|> 5r_{\tilde{l}_1,\tilde{l}_2}/2)\nonumber\\
    \leq& 2\delta,\nonumber
\end{align}
where the final inequality follows from \Cref{hetero-refinement-high-prob-bound} and \eqref{eq:hetero-high-prob-bound}.
\end{proof}
In \Cref{heterogeneity-step1-lemma}, we control, for any fixed $(l_1,l_2)$ satisfying $\mathcal{I}_{l_1,l_2}\neq \emptyset$, the initial estimator $\mathfrak{A}^*_{l_1,l_2}$ defined in \eqref{eq:A_l1l2-def} under the assumption that the condition~\eqref{step1-heterogeneous-delta-condition} holds.
\begin{lemma} \label{heterogeneity-step1-lemma}Consider the contamination model with heterogeneous server sizes detailed in~\eqref{Gaussian-distribution} and~\eqref{heterogeneous-model}. Fix any $(l_1,l_2)$ satisfying that $\mathcal{I}_{l_1,l_2} \neq \emptyset$. Given $\delta>0$ and the tuning parameter $r_{l_1,l_2}$, suppose that the condition 
\begin{equation}\label{step1-heterogeneous-delta-condition}
\varepsilon+ \sqrt{\frac{\log(12/(r_{l_1,l_2}\delta))}{|\mathcal{I}_{l_1,l_2}|}}+\frac{1}{\sqrt{(2^{l_1}n_{\min}\wedge2^{2(l_1+l_2)}n_{\min}^2\alpha_{\min}^2)|\mathcal{I}_{l_1,l_2}|}}<C
\end{equation}
holds for a sufficiently small absolute constant $C>0$. Then, $\mathfrak{A}^*_{l_1,l_2}$ defined in \eqref{eq:A_l1l2-def} satisfies
\[\mathbb{P}(|\mathfrak{A}^*_{l_1,l_2}-\theta|\leq 5r_{l_1,l_2}/2)\geq 1-\delta\]
provided that the tuning parameters $T_{l_1,l_2}$ and $r_{l_1,l_2}$ are as defined in \eqref{eq:hetero-UB-thm-T-def} and \eqref{eq:hetero-UB-thm-r-def}, where $C_1$ and $C_2$ are sufficiently large absolute constants.
\end{lemma}
\begin{proof} We adapt the proof strategy for \Cref{high-prob-lemma} to account for the small amount of heterogeneity within the batch group labelled by $(l_1,l_2)$. In \textbf{Step 1}, we derive a concentration bound for the aggregated smoothed votes around its expectation. In \textbf{Step 2}, we upper bound the expected smoothed votes assigned to bins that are far from $\theta$, and in \textbf{Step 3}, we lower bound the expected smoothed votes assigned to bins that are close to $\theta$. In \textbf{Step 4}, we combine the concentration bound in \textbf{Step 1} with the bounds on the expected smoothed votes in \textbf{Steps 2} and \textbf{3} to obtain the result corresponding to this heterogeneous setup. Since the bound $|\mathfrak{A}^*_{l_1,l_2}-\theta|\leq 2+r_{l_1,l_2}/2$ always holds, we only need to consider the case $r_{l_1,l_2} <1$ for the rest of the proof.

    \textbf{Step 1. }In this step, we control the deviation of the aggregated smoothed votes $|\mathcal{B}|^{-1}\sum_{k\in \mathcal{B}}Z_{k,j}$ from its expectation. Note that the aggregated privacy noise $16 T_{l_1,l_2}\sum_{k\in \mathcal{B}}{W_{k,j}/(\alpha_k n_k |\mathcal{B}|)} $ is sub-exponential with parameters
    \[\left(32T_{l_1,l_2}/(2^{l_1+l_2-2}n_{\min}\alpha_{\min}|\mathcal{B}|^{1/2}),32T_{l_1,l_2}/(2^{l_1+l_2-2}n_{\min}\alpha_{\min}|\mathcal{B}|)\right)\]
    due to the observation that $n_k \geq 2^{l_1-1}n_{\min}$ and $\alpha_k \geq 2^{l_2-1}\alpha_{\min}$. By an argument similar to the proof of \Cref{high-prob-lemma}, for any $t>0$ and $\gamma \in (0,1)$, we have that
    \begin{align}&\mathbb{P}\left(\max_{j\in[\lceil 2r_{l_1,l_2}^{-1}\rceil]}\left|\frac{1}{|\mathcal{B}|}\sum_{k\in \mathcal{B}}\{Z_{k,j} - \E[Z_{k,j}]\}\right|\geq t\right)\nonumber\\
    \stackrel{(a)}{\leq} & \sum_{j\in [\lceil 2r_{l_1,l_2}^{-1}\rceil]}\mathbb{P}\left(\left|\frac{1}{|\mathcal{B}|}\sum_{k\in \mathcal{B}}\{Z_{k,j} - \E[Z_{k,j}]\}\right|\geq t\right)\nonumber\\
    \stackrel{(b)}{\leq} &\sum_{j\in [\lceil 2r_{l_1,l_2}^{-1}\rceil]}\mathbb{P}\left(\left|\frac{1}{|\mathcal{B}|}\sum_{k\in \mathcal{B}}\left\{\left(2r_{l_1,l_2}-\left|\frac{1}{n_k}\sum_{i=1}^{n_k}[X_{k,i}]_{-T_{l_1,l_2}}^{T_{l_1,l_2}}-a_j\right|\right)_+ - \E[Z_{k,j}]\right\}\right|\geq (1-\gamma) t\right)\nonumber\\
    &\hspace{2cm}+\sum_{j\in [\lceil 2r_{l_1,l_2}^{-1}\rceil]}\mathbb{P}\left(\left|\frac{1}{|\mathcal{B}|}\sum_{k\in\mathcal{B}}\frac{16T_{l_1,l_2}}{\alpha_k n_k}W_{k,j}\right|\geq \gamma t\right)\nonumber\\
    \stackrel{(c)}{\leq} &\frac{6}{r_{l_1,l_2}}\exp\left(-\frac{|\mathcal{B}|(1-\gamma)^2t^2}{2r_{l_1,l_2}^2}\right)\nonumber\\
    &\label{eq:hetero-smoothed-vote-concentration-bound}\hspace{2cm}+ \frac{6}{r_{l_1,l_2}}\exp\left(-\frac{|\mathcal{B}|}{256}\min\left\{\frac{2^{2(l_1+l_2)}n_{\min}^2\alpha_{\min}^2\gamma^2 t^2}{128T_{l_1,l_2}^2},\frac{2^{l_1+l_2} n_{\min}\alpha_{\min} \gamma t}{T_{l_1,l_2}}\right\}\right),\end{align}
    where (a) follows from a union bound, (b) follows from the definition of $Z_{k,j}$ in \eqref{Z_kj-def-hetero} and a further union bound, and (c) follows from Hoeffding's inequality \citep[e.g.~Proposition 2.5 in][]{Wainwright2019high-dimensional} and Bernstein's inequality \citep[e.g.~Proposition 2.9 in][]{Wainwright2019high-dimensional}. Given $\delta \in (0,1)$, choose $\gamma = \{1+2^{l_1+l_2}n_{\min}\alpha_{\min}r_{l_1,l_2}/(128T_{l_1,l_2})\}^{-1}$ and

    \begin{align}t &= \sqrt{\frac{2\log(12/(r_{l_1,l_2}\delta))}{|\mathcal{B}|}}\max\left\{\frac{r_{l_1,l_2}}{(1-\gamma)},\frac{128T_{l_1,l_2}}{2^{l_1+l_2}n_{\min}\alpha_{\min}\gamma},\frac{128T_{l_1,l_2}}{2^{l_1+l_2}n_{\min}\alpha_{\min}\gamma}\sqrt{\frac{2\log(12/(r_{l_1,l_2}\delta))}{|\mathcal{B}|}}\right\}\nonumber\\
    \label{hetero-step-1-lemma-t-def}&= \left(r_{l_1,l_2}+\frac{128 T_{l_1,l_2}}{2^{l_1+l_2}n_{\min}\alpha_{\min}}\right)\sqrt{\frac{2\log(12/(r_{l_1,l_2}\delta))}{|\mathcal{B}|}},\end{align}
    where the final inequality follows from the assumption $2\log(12/(r_{l_1,l_2}\delta))\leq |\mathcal{B}|$ that follows from \eqref{step1-heterogeneous-delta-condition}. With this choice of $t$, it follows from \eqref{eq:hetero-smoothed-vote-concentration-bound} that the event 
    \begin{equation}\label{Omega_t-def-hetero}\Omega_t= \left\{\max_{j\in [\lceil 2r_{l_1,l_2}^{-1}\rceil]}\left|\frac{1}{|\mathcal{B}|}\sum_{k\in \mathcal{B}}\{Z_{k,j} - \E[Z_{k,j}]\}\right|< t\right\}\end{equation}
    occurs with probability at least $1-\delta$.

    \medskip
    \noindent \textbf{Step 2. }In this step, we upper bound the expected value of smoothed votes from bins far away from $\theta$. For $j$ such that $|\theta-a_j|>5r_{l_1,l_2}/2$, we can upper bound $\E Z_{k,j}$, for each $k\in \mathcal{I}_{l_1,l_2}$, by
    \begin{align}&\E Z_{k,j}\nonumber\\
    \leq& 2r_{l_1,l_2}\mathbb{P}\left(\left|\frac{1}{n_k}\sum_{i=1}^{n_k}[X_{k,i}]_{-T_{l_1,l_2}}^{T_{l_1,l_2}}-a_j\right|\leq 2r_{l_1,l_2}\right)\nonumber\\
    \leq& 2r_{l_1,l_2}\left\{\varepsilon +(1-\varepsilon)\mathbb{P}_{X_k\sim \mathcal{N}^{\otimes n_k}(\theta,1)}\left(\left|\frac{1}{n_k}\sum_{i=1}^{n_k}[X_{k,i}]_{-T_{l_1,l_2}}^{T_{l_1,l_2}}-\theta\right|\geq \frac{1}{2}r_{l_1,l_2}\right) \right\}\nonumber\\
    \stackrel{(a)}{\leq}& 2r_{l_1,l_2}\varepsilon +2r_{l_1,l_2}\nonumber\\
    &\cdot\mathbb{P}_{X_k\sim \mathcal{N}^{\otimes n_k}(\theta,1)}\left(\left|\frac{1}{n_k}\sum_{i=1}^{n_k}[X_{k,i}]_{-T_{l_1,l_2}}^{T_{l_1,l_2}}-\E_{X_k\sim \mathcal{N}^{\otimes n_k}(\theta,1)}\left\{[X_{k,1}]_{-T_{l_1,l_2}}^{T_{l_1,l_2}}\right\}\right|\geq \frac{r_{l_1,l_2}}{2}- \sqrt{2}e^{-T_{l_1,l_2}^2/16}\right) \nonumber\\
    \stackrel{(b)}{\leq}& 2r_{l_1,l_2}\left\{\varepsilon +\mathbb{P}_{X_k\sim \mathcal{N}^{\otimes n_k}(\theta,1)}\left(\left|\frac{1}{n_k}\sum_{i=1}^{n_k}[X_{k,i}]_{-T_{l_1,l_2}}^{T_{l_1,l_2}}-\E_{X_k\sim \mathcal{N}^{\otimes n_k}(\theta,1)}\left\{[X_{k,1}]_{-T_{l_1,l_2}}^{T_{l_1,l_2}}\right\}\right|\geq C_1'r_{l_1,l_2}\right) \right\}\nonumber\\
    \stackrel{(c)}{\leq}& 2r_{l_1,l_2}\left\{\varepsilon +2\exp({-C_1'^2n_k r_{l_1,l_2}^2/16})\right\}\nonumber\\
    \label{hetero-smoothed-vote-UB}\stackrel{(d)}{\leq}&2r_{l_1,l_2} \left(\varepsilon + \frac{2}{\sqrt{|\mathcal{I}_{l_1,l_2}|(2^{l_1}n_{\min}\wedge2^{2(l_1+l_2)}n_{\min}^2\alpha_{\min}^2)}}\right),\end{align}
    where $C_1' \in (0,1/2)$ is a positive constant, (a) follows from \eqref{Gaussian-trunc-lemma-bias}, (b) follows from the choices of $T_{l_1,l_2}$ and $r_{l_1,l_2}$ in \eqref{eq:hetero-UB-thm-T-def} and \eqref{eq:hetero-UB-thm-r-def}, (c) follows from \eqref{Gaussian-trunc-lemma-concentration}, and (d) follows from noting that $n_k\geq 2^{l_1-1}n_{\min}$ and the choices of $T_{l_1,l_2}$ and $r_{l_1,l_2}$ in \eqref{eq:hetero-UB-thm-T-def} and \eqref{eq:hetero-UB-thm-r-def}.

    \medskip
    \noindent\textbf{Step 3.} For $j$ such that $|\theta-a_j|\leq r_{l_1,l_2}/2$, we can lower bound $\E Z_{k,j}$, for each $k \in \mathcal{I}_{l_1,l_2}$, as follows: 
    \begin{align}
    \E Z_{k,j} &\geq (1-\varepsilon)\E_{X_k\sim\mathcal{N}^{\otimes n_k}(\theta,1)}\left[2r_{l_1,l_2} - \left|\frac{1}{n_k}\sum_{i=1}^{n_k}[X_{k,i}]_{-T_{l_1,l_2}}^{T_{l_1,l_2}} - \theta\right|-|\theta-a_j|\right]\nonumber+\varepsilon\cdot 0\\
        &\geq (1-\varepsilon) \left(\frac{3}{2}r_{l_1,l_2}-\left\{\E_{X_k\sim \mathcal{N}^{\otimes n_k}(\theta,1)}\left[\left|\frac{1}{n_k}\sum_{i=1}^{n_k}[X_{k,i}]_{-T_{l_1,l_2}}^{T_{l_1,l_2}}-\theta\right|^2\right]\right\}^{1/2}\right)\nonumber\\
        &\stackrel{(a)}{\geq} (1-\varepsilon)\left(\frac{3}{2}r_{l_1,l_2}-\left\{\frac{1}{n_k}+ 2e^{-T_{l_1,l_2}^2/8}\right\}^{1/2}\right)\nonumber\\
    &\label{hetero-smoothed-vote-LB}\stackrel{(b)}{\geq} (1-\varepsilon)r_{l_1,l_2},\end{align}
    where (a) follows from \eqref{Gaussian-trunc-lemma-bias}, \eqref{Gaussian-trunc-lemma-var} and the bias-variance decomposition, and (b) follows from the choices of $T_{l_1,l_2}$ and $r_{l_1,l_2}$ in \eqref{eq:hetero-UB-thm-T-def} and \eqref{eq:hetero-UB-thm-r-def}.

    \medskip
    \noindent\textbf{Step 4. }On the event $\Omega_t$ defined in \eqref{Omega_t-def-hetero}, we have, for $j$ satisfying $|\theta-a_j|\leq r_{l_1,l_2}/2$ and $l$ satisfying $|\theta-a_l|>5r_{l_1,l_2}/2$,
    \begin{align}\frac{1}{|\mathcal{B}|}\sum_{k\in \mathcal{B}} (Z_{k,j}- Z_{k,l})
    >&r_{l_1,l_2}\left(1-3\varepsilon-\frac{4}{\sqrt{|\mathcal{I}_{l_1,l_2}|(2^{l_1}n_{\min}\wedge2^{2(l_1+l_2)}n_{\min}^2\alpha_{\min}^2)}}\right)-2t\nonumber\\
    =&r_{l_1,l_2}\left(1-3\varepsilon-\frac{4}{\sqrt{|\mathcal{I}_{l_1,l_2}|(2^{l_1}n_{\min}\wedge2^{2(l_1+l_2)}n_{\min}^2\alpha_{\min}^2)}}\right)\nonumber\\
    \label{smoothed-vote-gap}&\hspace{4cm}-2\left(r_{l_1,l_2}+\frac{128 T_{l_1,l_2}}{2^{l_1+l_2}n_{\min}\alpha_{\min}}\right)\sqrt{\frac{2\log(12/(r_{l_1,l_2}\delta))}{|\mathcal{B}|}}\end{align}
    due to \eqref{hetero-smoothed-vote-UB}, \eqref{hetero-smoothed-vote-LB} and the choice of $t$ in \eqref{hetero-step-1-lemma-t-def}.
    Hence, we have $|\mathfrak{A}^*_{l_1,l_2}-\theta|\leq 5r_{l_1,l_2}/2$ on the event $\Omega_t$ provided that the RHS of \eqref{smoothed-vote-gap} is positive.
    Noting that $|\mathcal{B}|=|\mathcal{I}_{l_1,l_2}|/2$ and 
    \[\frac{ T_{l_1,l_2}}{2^{l_1+l_2}n_{\min}\alpha_{\min}r_{l_1,l_2}}\sqrt{\frac{\log(12/(r_{l_1,l_2}\delta))}{|\mathcal{I}_{l_1,l_2}|}}\leq \frac{1}{C_2}\sqrt{\frac{\log(12/(r_{l_1,l_2}\delta))}{T_{l_1,l_2}^2+\log(4/\delta)}}\leq\frac{1}{C_2}\] by the definition \eqref{eq:hetero-UB-thm-r-def} of $r$, we conclude that \eqref{step1-heterogeneous-delta-condition} is sufficient for the RHS of \eqref{smoothed-vote-gap} to be positive, and hence sufficient for
    \[\mathbb{P}(|\mathfrak{A}^*_{l_1,l_2}-\theta|\leq 5r_{l_1,l_2}/2)\geq 1-\delta\]
    to hold.
\end{proof}
\begin{lemma}\label{hetero-refinement-high-prob-bound} 

     Consider the contamination model with heterogeneous server sizes detailed in~\eqref{Gaussian-distribution} and~\eqref{heterogeneous-model}. Fix any $(l_1,l_2)$ satisfying that $\mathcal{I}_{l_1,l_2} \neq \emptyset$. Suppose that the event $|\mathfrak{A}^*_{l_1,l_2}-\theta|\leq 5r_{l_1,l_2}/2$ from \Cref{heterogeneity-step1-lemma} holds. Given $\delta \in [4e^{-|\mathcal{I}_{l_1,l_2}|/4},1]$, the estimator $\hat{\theta}_{l_1,l_2}$ defined in \eqref{theta-hat_l1l2-def} satisfies
    \[\mathbb{P}\left(|\hat{\theta}_{l_1,l_2} - \theta|\geq L(\varepsilon, l_1,l_2, \delta)\big| \mathfrak{A}^*_{l_1,l_2}\right) \leq  \delta\]
with tuning parameters $T_{l_1,l_2}$, $r_{l_1,l_2}$ and $M_{l_1,l_2}$ defined in \eqref{eq:hetero-UB-thm-T-def} and \eqref{eq:hetero-UB-thm-r-def}, 
where the function $L(\cdot,\cdot,\cdot,\cdot)$ is defined in \eqref{hetero-L-def}
and $C_1,C_2>0$ are sufficiently large absolute constants.
\end{lemma}
\begin{proof} In \textbf{Step 1}, we adapt the proof of \Cref{E-S-given-J-lemma} to bound the bias of $Z_k$ as an estimator of $\theta$, assuming that the event $|\mathfrak{A}^*_{l_1,l_2}-\theta|\leq 5r_{l_1,l_2}/2$ holds. In \textbf{Step 2}, we prove a concentration bound for $\hat{\theta}_{l_1,l_2}$ using Bernstein's inequality and Hoeffding's inequality, similar to the arguments in \textbf{Step 1} of the proof of \Cref{heterogeneity-step1-lemma}.  In \textbf{Step 3}, we combine the bias bound from \textbf{Step 1} and the concentration bound from \textbf{Step 2} to show the high probability guarantee for $\hat{\theta}_{l_1,l_2}$.

\medskip
\noindent\textbf{Step 1. }Recall that $\mathcal{I}_{l_1,l_2}$ is split into two halves, $\mathcal{B}$ and $\mathcal{B}'$, such that only the first half $\mathcal{B}$ is used in the construction of the initial estimator $\mathfrak{A}^*_{l_1,l_2}$. Following the proof of \Cref{E-S-given-J-lemma} with $n = n_k$, $r= r_{l_1,l_2}$, $T= T_{l_1,l_2}$ and $M= M_{l_1,l_2}$ instead, and observing that $r_{l_1,l_2}\vee n_k^{-1/2}\lesssim M_{l_1,l_2}$ for each $k \in \mathcal{I}_{l_1,l_2}$, we obtain the following version of \eqref{Gaussian-three-cases-truncation-2}. For each $k \in \mathcal{B}'$, we have that
\begin{align}
    &|\E_{X_k\sim \mathcal{N}^{\otimes n_k}(\theta,1)} [Z_k-\theta|\mathfrak{A}^*_{l_1,l_2}]|\nonumber\\
    =&\left|\E_{X_k\sim \mathcal{N}^{\otimes n_k}(\theta,1)}\left[\mathfrak{A}^*_{l_1,l_2}-\theta+\left[\frac{1}{n_k}\sum_{i=1}^{n_k}[X_{k,i}]_{-T_{l_1,l_2}}^{T_{l_1,l_2}}-\mathfrak{A}^*_{l_1,l_2}\right]_{-M_{l_1,l_2}}^{M_{l_1,l_2}}\Bigg|\mathfrak{A}^*_{l_1,l_2}\right]\right|\nonumber\\
    \lesssim &e^{-T_{l_1,l_2}^2/16}+((M_{l_1,l_2}+1)\mathbf{1}\{|\mathfrak{A}^*_{l_1,l_2}-\theta|>5r_{l_1,l_2}/2\}+M_{l_1,l_2}\mathbf{1}\{|\mathfrak{A}^*_{l_1,l_2}-\theta|\leq5r_{l_1,l_2}/2\})\nonumber\\
    &\label{hetero-lemma-expectation bound}\cdot\left\{\mathbb{P}_{X_k\sim \mathcal{N}^{\otimes n_k}(\theta,1)}\left(\left.\left|\frac{1}{n_k}\sum_{i = 1}^{n_k} [X_{k,i}]_{-T_{l_1,l_2}}^{T_{l_1,l_2}} - \mathfrak{A}^*_{l_1,l_2}\right|>M_{l_1,l_2}\right|\mathfrak{A}^*_{l_1,l_2}\right)\right\}^{1/2}.
\end{align}
Continuing to follow the proof of \Cref{E-S-given-J-lemma}, we have the following version of \eqref{eq:gaussian-E-S-given-J-prob-bound}. For each $k\in \mathcal{B}'$, we have that 
\begin{align}
    &\mathbb{P}_{X_k\sim \mathcal{N}^{\otimes n_k}(\theta,1)}\left(\left.\left|\frac{1}{n_k}\sum_{i = 1}^{n_k} [X_{k,i}]_{-T_{l_1,l_2}}^{T_{l_1,l_2}} - \mathfrak{A}^*_{l_1,l_2}\right|>M_{l_1,l_2}\right|\mathfrak{A}^*_{l_1,l_2}\right)\nonumber\\
    &\stackrel{}{\leq} \mathbf{1}\{|\mathfrak{A}^*_{l_1,l_2}-\theta|\leq5r_{l_1,l_2}/2\}\exp\left(-n_k(M_{l_1,l_2}-5r_{l_1,l_2}/2-\sqrt{2}e^{-T_{l_1,l_2}^2/16})^2/16\right)\nonumber\\
        &\hspace{1cm}+\mathbf{1}\{|\mathfrak{A}^*_{l_1,l_2}-\theta|>5r_{l_1,l_2}/2\}\nonumber\\
        &\hspace{2cm}\cdot\left[\exp\left(-n_k(M_{l_1,l_2}-2-\sqrt{2}e^{-T_{l_1,l_2}^2/16})^2/16\right)\mathbf{1}\{M_{l_1,l_2}>3\}+\mathbf{1}\{M_{l_1,l_2}\leq 3\}\right]\nonumber\\
        &\label{eq:hetero-lemma-prob-bound}\leq e^{-n_k(M_{l_1,l_2}/12)^2/16}+\mathbf{1}\{|\mathfrak{A}^*_{l_1,l_2}-\theta|>5r_{l_1,l_2}/2\}\mathbf{1}\{M_{l_1,l_2}\leq 3\}.
\end{align}
Substituting \eqref{eq:hetero-lemma-prob-bound} into \eqref{hetero-lemma-expectation bound}, we have that, for each $k \in \mathcal{B}'$,
\begin{align}&|\E_{X_k \sim \mathcal{N}^{\otimes n_k}(\theta,1)}[Z_k-\theta|\mathfrak{A}^*_{l_1,l_2}]|\nonumber\\&\lesssim e^{-T_{l_1,l_2}^2/16}+M_{l_1,l_2}e^{-n_kM_{l_1,l_2}^2/4608}+(e^{-n_kM_{l_1,l_2}^2/4608}+1)\mathbf{1}\{|\mathfrak{A}^*_{l_1,l_2}-\theta|>5r_{l_1,l_2}/2\}\nonumber\\
&\stackrel{(a)}{\lesssim} e^{-T_{l_1,l_2}^2/16}+e^{-n_kM_{l_1,l_2}^2/4609}+\mathbf{1}\{|\mathfrak{A}^*_{l_1,l_2}-\theta|> 5r_{l_1,l_2}/2\}\nonumber\\
&\label{hetero-smoothed-vote-bias-bound}\lesssim \frac{1}{\sqrt{|\mathcal{I}_{l_1,l_2}|}(\sqrt{2^{l_1}n_{\min}}\wedge 2^{l_1+l_2}n_{\min}\alpha_{\min})}+\mathbf{1}\{|\mathfrak{A}^*_{l_1,l_2}-\theta|> 5r_{l_1,l_2}/2\},\end{align}
where (a) follows from noting that $M_{l_1,l_2}\lesssim n_k M_{l_1,l_2}^2\lesssim e^{(1/4608-1/4609)n_kM_{l_1,l_2}^2}$ by the definition of $M_{l_1,l_2}$ in \eqref{eq:hetero-UB-thm-r-def}. Hence, on the event $|\mathfrak{A}^*_{l_1,l_2}-\theta|\leq 5r_{l_1,l_2}/2$, we have that
\begin{align}|\E [Z_k|\mathfrak{A}^*_{l_1,l_2}]-\theta| &\leq (1-\varepsilon)|\E_{X_k \sim \mathcal{N}^{\otimes n_k}(\theta,1)}[Z_k-\theta|\mathfrak{A}^*_{l_1,l_2}]| + \varepsilon(|\mathfrak{A}^*_{l_1,l_2}-\theta|+M_{l_1,l_2})\nonumber\\
&\label{hetero-refinement-bias-bound}\lesssim  \frac{1}{\sqrt{|\mathcal{I}_{l_1,l_2}|}(\sqrt{2^{l_1}n_{\min}}\wedge 2^{l_1+l_2}n_{\min}\alpha_{\min})}+\varepsilon M_{l_1,l_2},
\end{align}
where the last inequality follows from \eqref{hetero-smoothed-vote-bias-bound} and that $r_{l_1,l_2}\lesssim M_{l_1,l_2}$.

\medskip
\noindent\textbf{Step 2. }Note that, from \eqref{eq-def-Z_k-hetero} and \eqref{theta-hat_l1l2-def2}, $\hat{\theta}_{l_1,l_2}$ is the sum of private noise and a part that is bounded between $\mathfrak{A}^*_{l_1,l_2}-M_{l_1,l_2}$ and $\mathfrak{A}^*_{l_1,l_2}+M_{l_1,l_2}$. The private noise $\sum_{k\in B'_{l_1,l_2}}2T_{l_1,l_2} L_k/(\alpha_k n_k |\mathcal{I}_{l_1,l_2}|/2)$ is sub-exponential with parameters \[(4T_{l_1,l_2} /(2^{l_1+l_2-2}n_{\min}\alpha_{\min}|\mathcal{B}'|^{1/2}), 4T_{l_1,l_2}/(2^{l_1+l_2-2}n_{\min}\alpha_{\min}|\mathcal{B}'|)).\]
Using Hoeffding's inequality \citep[e.g.~Proposition 2.5 in][]{Wainwright2019high-dimensional} on the bounded part and Bernstein's inequality \citep[e.g.~Proposition 2.9 in][]{Wainwright2019high-dimensional} on the private noise, we have that, for any $\gamma \in (0,1)$,
    \begin{multline*}
        \mathbb{P}\left(\left|\hat{\theta}_{l_1,l_2}- \E \left[\frac{1}{|\mathcal{B}'|}\sum_{k\in \mathcal{B}'}Z_k \bigg|\mathfrak{A}^*_{l_1,l_2}\right]\right|\geq t\bigg|\mathfrak{A}^*_{l_1,l_2}\right) \leq 2\exp\left(-\frac{|\mathcal{B}'|(1-\gamma)^2t^2}{2M_{l_1,l_2}^2}\right)+ \\
        2\exp\left(-\left(\frac{2^{2(l_1+l_2-2)}n_{\min}^2\alpha_{\min}^2|\mathcal{B}'| \gamma ^2t^2}{32T_{l_1,l_2}^2}\wedge \frac{2^{l_1+l_2-2}n_{\min}\alpha_{\min}|\mathcal{B}'|\gamma t}{8T_{l_1,l_2}}\right)\right).
    \end{multline*}
    Since $4\log(4/\delta) \leq |\mathcal{I}_{l_1,l_2}|$, choosing $\gamma = (1+2^{l_1+l_2}n_{\min}\alpha_{\min}M_{l_1,l_2}/(16T_{l_1,l_2}))^{-1}$ and
    \[t = \sqrt{\frac{2\log(4/\delta)}{|\mathcal{B}'|}}\left(M_{l_1,l_2}+\frac{16T_{l_1,l_2}}{2^{l_1+l_2}n_{\min}\alpha_{\min}}\right)\]
gives us that 
\begin{equation}\label{hetero-refinement-concentration-bound}
    \mathbb{P}\left(\left|\hat{\theta}_{l_1,l_2}- \E \left[\frac{1}{|\mathcal{B}'|}\sum_{k\in \mathcal{B}'}Z_k \bigg|\mathfrak{A}^*_{l_1,l_2}\right]\right|\geq \sqrt{\frac{2\log(4/\delta)}{|\mathcal{B}'|}}\left(M_{l_1,l_2}+\frac{16T_{l_1,l_2}}{2^{l_1+l_2}n_{\min}\alpha_{\min}}\right)\bigg|\mathfrak{A}^*_{l_1,l_2}\right)\leq\delta.
\end{equation}
    
\medskip
    \noindent\textbf{Step 3.} Combining the results of \textbf{Step 1} and \textbf{Step 2}, we have that, on the event $|\mathfrak{A}^*_{l_1,l_2}-\theta|\leq 5r_{l_1,l_2}/2$,
    \begin{align*}
        &\mathbb{P}\left(\left|\hat{\theta}_{l_1,l_2}- \theta\right|\leq C\left\{\varepsilon M_{l_1,l_2}\vee\sqrt{\frac{\log(4/\delta)}{|\mathcal{I}_{l_1,l_2}|}}\left(M_{l_1,l_2}\vee\frac{T_{l_1,l_2}}{2^{l_1+l_2}n_{\min}\alpha_{\min}}\right)\right\}\bigg|\mathfrak{A}^*_{l_1,l_2}\right)\\
        \geq& \mathbb{P}\left(\left|\hat{\theta}_{l_1,l_2}- \E \left[\frac{1}{|\mathcal{B}'|}\sum_{k\in \mathcal{B}'}Z_k \bigg|\mathfrak{A}^*_{l_1,l_2}\right]\right|\leq \sqrt{\frac{2\log(4/\delta)}{|\mathcal{B}'|}}\left\{M_{l_1,l_2}+\frac{16T_{l_1,l_2}}{2^{l_1+l_2}n_{\min}\alpha_{\min}}\right\}\bigg|\mathfrak{A}^*_{l_1,l_2}\right)\\
        \geq& 1-\delta,
    \end{align*}
    where the first inequality follows from \eqref{hetero-refinement-bias-bound} and noting that $|\mathcal{B}'| = |\mathcal{I}_{l_1,l_2}|/2$, $M_{l_1,l_2}\gtrsim 1/\sqrt{2^{l_1}n_{\min}}$ and $T_{l_1,l_2}\gtrsim 1$, and the final inequality follows from \eqref{hetero-refinement-concentration-bound}.
\end{proof}
\end{document}